\documentclass[10pt,a4paper]{article}
\pdfoutput=1 %this serves for arxiv

\usepackage{lineno}

\modulolinenumbers[5]

\usepackage{amssymb}
\usepackage{amsmath}
\usepackage{algorithm}     %serve per scrivere lo pseudo codice
\usepackage{algpseudocode} %serve per scrivere lo pseudo codice
\usepackage{subfig}
\usepackage{csquotes}
\usepackage{mathtools}
\usepackage{bbm}
\usepackage{amsthm}
\usepackage{float}
 \usepackage{dutchcal}

\usepackage{hyperref}
\usepackage{cleveref}
\usepackage{wrapfig}

 \usepackage{tikz}
\usetikzlibrary{patterns}
\usetikzlibrary{arrows,shapes,backgrounds}
\usepackage{pgfplots}
\usepgfplotslibrary{external}
\usetikzlibrary{calc}

\usepackage{geometry}
\newcommand{\email}[1]{\hspace*{\stretch{1}}\emph{\texttt{#1}}}

\makeatletter
\def\blfootnote{\xdef\@thefnmark{$\star$}\@footnotetext}
\makeatother
\newenvironment{Authors}%
  {\begin{center}\begin{bfseries}}%
  {\end{bfseries}\end{center}}
\newenvironment{Addresses}%
  {\begin{flushleft}\begin{itshape}}%
  {\end{itshape}\end{flushleft}}

 \usepackage{fancyhdr}  
  
\fancypagestyle{plain}{
	\fancyhead{}
	\fancyhead[C]{\hfill submitted to Siam  Journal on Scientific Computing, February 2022}
}

\newtheorem{theorem}{Theorem}[section]
\newtheorem{theorem1}{Theorem}[section]
\newtheorem{theorem2}{Theorem}[section]
\newtheorem{proposition}[theorem]{Proposition}
\newtheorem{lemma}[theorem1]{Lemma}
\newtheorem{remark}[theorem2]{Remark}

  \newcommand{\vertiii}[1]{{\left\vert\kern-0.25ex\left\vert\kern-0.25ex\left\vert #1 
    \right\vert\kern-0.25ex\right\vert\kern-0.25ex\right\vert}}

\begin{document}

%\begin{frontmatter}
\thispagestyle{plain}

\title{Localized model reduction for nonlinear elliptic partial differential equations: localized training, 
partition of unity, and adaptive enrichment.}
 \date{}
 
 \maketitle
\vspace{-50pt} 
 
\begin{Authors}
Kathrin Smetana$^{1}$
Tommaso Taddei$^{2}$
\end{Authors}

\begin{Addresses}
$^1$
Stevens Institute of Technology, 
Department of Mathematical Sciences,
1 Castle Point Terrace, Hoboken, NJ 07030, USA \email{ksmetana@stevens.edu}\\
$^2$
IMB, UMR 5251, Univ. Bordeaux;  33400, Talence, France.
Inria Bordeaux Sud-Ouest, Team MEMPHIS;  33400, Talence, France, \email{tommaso.taddei@inria.fr} \\
\end{Addresses}

\begin{abstract}
We  {propose} a component-based (CB) parametric model order reduction (pMOR) formulation for parameterized  {nonlinear} elliptic partial differential equations (PDEs). 
CB-pMOR is designed to deal with large-scale problems for which full-order solves  are not affordable in a reasonable time frame  or parameters' variations induce topology changes that prevent the application of monolithic pMOR techniques.
We rely on the partition-of-unity method (PUM) to devise global approximation spaces from local reduced spaces, and on Galerkin projection to compute the global state estimate.
We propose a randomized data compression algorithm based on oversampling for the construction of the components' reduced spaces: the approach exploits random boundary conditions of controlled smoothness on the oversampling boundary.
We further propose an adaptive residual-based enrichment algorithm that exploits global reduced-order solves on representative systems to update the local reduced spaces.
We prove exponential convergence of the enrichment procedure for linear coercive problems; we further present numerical results for a two-dimensional nonlinear diffusion problem to illustrate the  many  features of our proposal and demonstrate its effectiveness.
\end{abstract}

% REQUIRED

\section{Introduction}
\label{sec:introduction}

%Acronyms to be introduced somewhere in the introduction.
%\begin{itemize}
%\item
%pMOR = parametric model order reduction.
%\item
%PDE = partial differential equation
%\item
%\texttt{hf} = high-fidelity.
%\item
%CB-pMOR = component-based parametric model order reduction.
%\end{itemize}

\subsection{Component-based model reduction for  parameterized PDEs}

Numerical modeling and simulation is of paramount importance to predict the response, improve the design, and monitor the structural health of engineering systems.
Several problems of interest involve repeated solutions of a partial differential equation   (PDE) for many values of the model parameters or require real-time responses: 
these tasks are prohibitively expensive for standard (e.g., finite element) methods.
Parametric model order reduction (pMOR, 
\cite{Haa17,HeRoSt16,QuMaNe16}) aims to reduce the marginal cost associated with the solution to parameterized systems  over a range of parameters.
The goal of this paper is to develop a pMOR procedure for large-scale nonlinear elliptic PDEs with parameter-induced topology changes.

pMOR techniques may rely on  an offline/online decomposition  to reduce marginal costs.
During the offline phase, we rely  on several high-fidelity (HF) solves to generate a 
reduced-order model (ROM) for the solution field.
During the online phase, given a new value of the parameter, we query the  ROM to estimate the solution field  and relevant quantities of interest.
Monolithic pMOR methods rely on HF solves at the training stage, which  might be unaffordable for very large-scale problems.
Furthermore, pMOR methods rely  on the assumption  that the solution field is defined over a parameter-independent domain or over a family of diffeomorphic domains: they thus cannot deal with  problems for which parametric variations induce topology changes.

To address these issues, several authors have proposed component-based pMOR procedures (cf. \cite{HuKnPa13} and the review \cite{Buhr2020localized}).
During  the offline stage, we define a  library of \emph{archetype components} and we build local reduced-order bases (ROBs) and local ROMs; then, during the online stage, we instantiate components to form the global system and we estimate the global solution by coupling local ROMs. CB-pMOR strategies consist of two distinct  building blocks:
(i) a rapid and reliable domain decomposition (DD)  strategy for online global predictions, and
(ii) a localized training strategy exclusively based on local solves for the construction of the local approximations.

{CB-pMOR shares important features with multiscale methods  
\cite{BabLip11,EfeHou07,le2014msfem,strouboulis2000design,Owh15,
Owh17,OwZ11,OZB13,MalPet14,MalPet21,SmePat16,MaScDo21,GrGrSa12, Ver19,Liu21,CHOS21,CLLW19,CEGL16,Chenetal20,HouWu97,Hug95}. Similarly to CB-pMOR, multiscale methods rely on local solves to build suitable approximation  spaces  that are tailored to the problem of interest.} The  emphasis in CB-pMOR is to devise  and  then exploit a library of inter-operable archetype components and associated ROMs that can be used for a broad range of {potentially parameter dependent} problems in a specific domain of interest.
%This work is a contribution in CB-pMOR: nevertheless, in the future, we wish to extend our proposals to the multiscale setting.

%\todo[inline]{KS: People from the DD community (I am collaborating with them) will not appreciate if we call that subsection DD methods and then provide mainly Localized MOR references. Therefore I changed.}

\subsection{Domain decomposition strategies within CB-pMOR}
Since the seminal work by Maday and R{\o}nquist
\cite{MadRon02} ---  that proposed a non-overlapping non-conforming reduced basis element method based on   Mortar DD --- 
several authors have combined  DD methods with model reduction methods to devise effective CB-ROMs.
As discussed in detail in the review \cite{Buhr2020localized}, 
we can distinguish between
 conforming non-overlapping  approaches {\cite{HuKnPa13,EftPat13,SmePat16}},
 non-conforming non-overlapping approaches based on Lagrange multipliers
 \cite{hoang2021domain,IaQuRo12,MadRon02,pegolotti2021model},  non-conforming non-overlapping approaches based on discontinuous Galerkin (DG) coupling
\cite{AHOK12,AnPaQu16,OhlSch15,riffaud2021dgdd}, and 
overlapping methods  
 \cite{bergmann2018zonal,buhr2018randomized}.
{The vast majority of existing contributions
(with few recent exceptions 
\cite{baiges2013domain,bergmann2018zonal,hoang2021domain,pegolotti2021model})
is restricted to linear PDEs.} 
%, and does not address the issues of hyper-reduction 
%--- which is crucial to achieve significant speedups ---
%and a posteriori error estimation.
%\todo[inline]{KS: Many papers actually do address a posteriori error estimation and if people do only look at linear PDEs they do not need to look at hyper-reduction?}

In this work, we rely on the partition-of-unity method (PUM) to devise global approximation spaces from local reduced spaces, and on Galerkin projection to compute the global state estimate.
PUM was proposed by Babu{\v{s}}ka and Melenk in 
\cite{BabMel97,MelBab96} and further developed and analyzed  in the framework of generalized finite element methods for multiscale problems 
(cf. \cite{babuvska2020multiscale});
PUM was also considered in the pMOR literature for linear elliptic and parabolic problems
\cite{buhr2018randomized,schleuss2020optimal}.
In  the CB-pMOR framework, PUM offers 
{ a general (i.e., independent of the underlying PDE) framework with} strong theoretical guarantees.

\subsection{Localized training based on oversampling and randomization}
Given the domain $\widehat{\Omega}$ associated with a given  archetype component, oversampling methods consist in
(i) defining a patch $\widehat{U} \supset \widehat{\Omega}$ and a suitable local PDE problem in $\widehat{U}$,  
(ii) solving the local PDE for several choices of the boundary conditions on $\partial \widehat{U}$ and then restricting the solution to $\widehat{\Omega}$, and finally
(iii) exploiting the results to build a local approximation space for the solution in $\widehat{\Omega}$.
Randomized methods rely on independent identically distributed (iid) samples  of the boundary conditions on (a subset  of) $\partial \widehat{U}$: they thus require the introduction of a probability density function for the functions defined on $\partial \widehat{U}$.

Oversampling methods exploit low-pass filtering properties of the differential operator to identify low-dimensional structures:
we refer to  \cite[Chapter 5]{Tad17} and
\cite[Remark 3.3]{SmePat16} for two representative working examples.
In detail, Caccioppoli type inequalities (see e.g. \cite{GiaSou85}) provide the theoretical foundations for the  application of oversampling methods to a particular class of PDEs.
{Oversampling methods have been suggested and used extensively in the context of multiscale methods (see e.g. 
\cite{HouWu97,BabLip11,HenPet13,MalPet14} and references therein) and then used as well in CB-pMOR
\cite{EftPat13,SmePat16} for linear PDEs.}

As {suggested} in \cite{buhr2018randomized}, randomized oversampling methods for linear parameter- independent PDEs can be linked to randomized singular value decomposition (SVD) techniques developed and analyzed in randomized numerical linear algebra 
 \cite{HaMaTr11,MaRoTy11,Mah11,Sar06,RoSzTy09,DriMah16,MahDri09,Woo14}: this link allows to extend methodological and theoretical contributions in randomized linear algebra to CB-pMOR.
{In particular, we can exploit concentration inequalities
(see e.g. \cite{BoLuMa13})
to analyze the error of randomized techniques}, and inform the choice of the sampling distribution.
The influence of the choice of the sampling distribution for nonlinear PDEs  remains an open question  in CB-pMOR.

%\begin{itemize}
%\item
%Two major elements:
%(i) evanescence of high-frequency modes and (ii) performance of randomized methods for  operators with rapidly-decaying singular values;
%Emphasize that they are two separate factors.
%\item
%Evanescence of high-frequency modes: acoustic waveguide (TT: PhD thesis TT Chapter 5) or diffusion equation (KS,ATP: SIAM 2016, Remark 3.3). 
%\textbf{References:} 
%working examples: 
%
%oversampling:
%\cite{BabLip11,HenPet13};
%in MOR \cite{EftPat13,SmePat16} 
%Caccioppoli inequalities \cite{GiaSou85}.
%\item
%Randomized methods:
%explain in words why randomized methods are good to identify structures; mention measure concentration effect.
%Randomized linear algebra \cite{HaMaTr11}.
%non-asymptotic concentration inequalities:
%\cite{BoLuMa13}.
%A few results that matter (in the linear case)
%\cite[Theorem 10.8]{HaMaTr11};
%\cite[Proposition 10.6]{HaMaTr11};
%results for Gaussians
%\cite[Propositions 10.1, 10.2]{HaMaTr11}.
%\end{itemize}

\subsection{Contributions of the paper and outline}
In this work, we propose a CB-pMOR procedure based on PUM for parametric nonlinear elliptic PDEs; {we do not require the nonlinear operator to be monotone.} The contributions of the paper are twofold.
First, we propose a randomized data compression algorithm based on oversampling: the approach relies on random samples of local parameters and boundary conditions on the oversampling boundary.
We propose a new sampler that controls the smoothness of the boundary condition,  and we empirically demonstrate its effectiveness for a nonlinear diffusion problem.
Second, we propose a basis enrichment algorithm that relies on global reduced solves to  enrich the local reduced spaces.
The algorithm relies on a {local} residual-based error indicator to identify boundary conditions for which the local ROM is inaccurate {and a rigorous global a posteriori error bound as a termination criterion}. We prove in-sample a priori exponential convergence of the enrichment algorithm for linear coercive problems; we further investigate performance for a nonlinear diffusion problem.

Our randomized  algorithm reads as a randomized 
proper orthogonal decomposition \cite{yu2015randomized} 
with respect to parameter and boundary  conditions.
On the other hand, the enrichment algorithm is closely related 
to the online enrichment strategy proposed in \cite{OhlSch15}
for non-overlapping  DG DD, and to the residual-based online enrichment algorithm considered  in 
\cite{buhr2018exponential} {for linear problems}.
The major difference is that the enrichment is  performed at training stage and aims to update the local approximation spaces associated with the archetype components, rather than during the online stage on the ``instantiated components''.

The outline of the paper is as follows.
In  \cref{sec:model_problem}, we introduce the model problem considered throughout the paper to illustrate the main definitions and to numerically validate our proposal: the model problem involves a high-dimensional ($\mathcal{O}(10^2)$) parameterization and topology  changes.
In   \cref{sec:domain_decomposition}, we discuss  the DD strategy based on the PUM and we introduce local and global discrete approximation spaces;
in  \cref{sec:data_compression},  we discuss the randomized localized data compression;
in  \cref{sec:local_global}, we present the enrichment strategy; and  in   \cref{sec:numerics}  we present thorough numerical investigations for the model problem.
 \Cref{sec:conclusions} {concludes} the paper.
 
\section{Model problem: nonlinear diffusion}
\label{sec:model_problem}
 
\begin{figure}[t]
\centering

\subfloat[ ]{
\includegraphics[width=0.44\textwidth]
{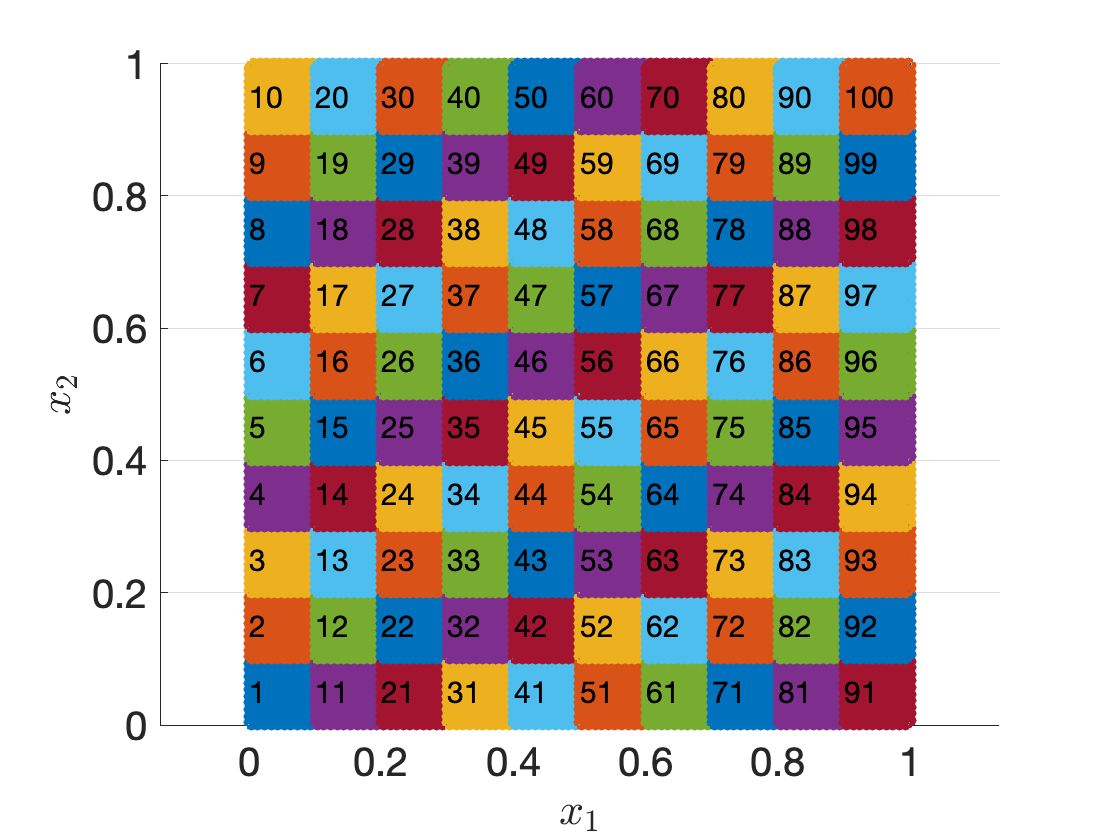}
}
~~~~~~~~~~
\subfloat[ ]{
\begin{tikzpicture}[scale=0.9]
\draw [black,thick] (-3,-3) -- (3,-3) -- (3,3) --
(-3,3) --(-3,-3) ; 

\draw [blue,dashed,ultra thick] (-3,-3) -- (-2,-3) -- (-2,-2) -- (-3,-2) --(-3,-3) ; 

\coordinate [label={right:  {\Large{\color{blue} co}}}] (E) at (-3.05,-2.6) ;

\draw [red,dashed,ultra thick] (-3,1) -- (-2,1) -- (-2,2) -- (-3,2) --(-3,1) ; 

\coordinate [label={right:  {\Large{\color{red} ed}}}] (E) at (-3.05,1.5) ;

\draw [purple,dashed,ultra thick] (0,0) -- (1,0) -- (1,1) -- (0,1) -- (0,0);

\coordinate [label={right:  {\Large{\color{purple} int}}}] (E) at (-0.12,0.5) ;

\coordinate [label={left:  {\Huge{$\Omega$}}}] (E) at (3,2.3) ;
 
\end{tikzpicture}
}

\subfloat[ ]{
\includegraphics[width=0.44\textwidth]
{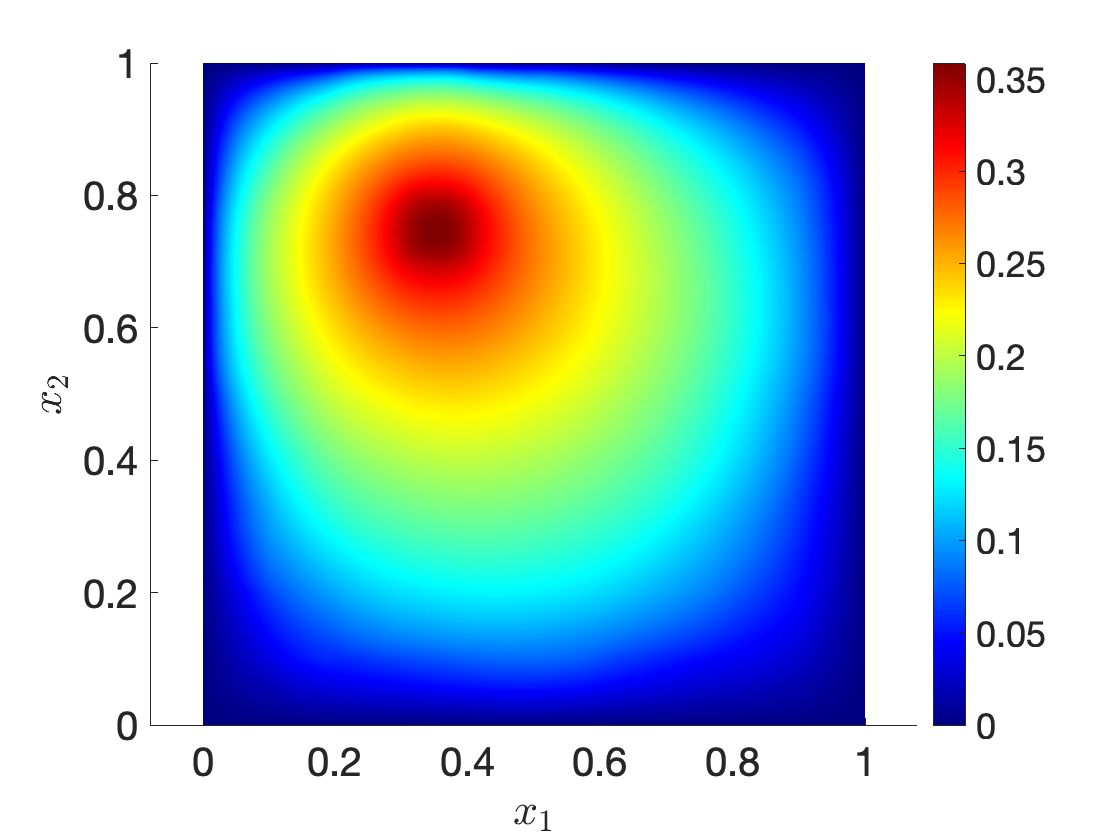}
}
~~
\subfloat[ ]{
\includegraphics[width=0.44\textwidth]
{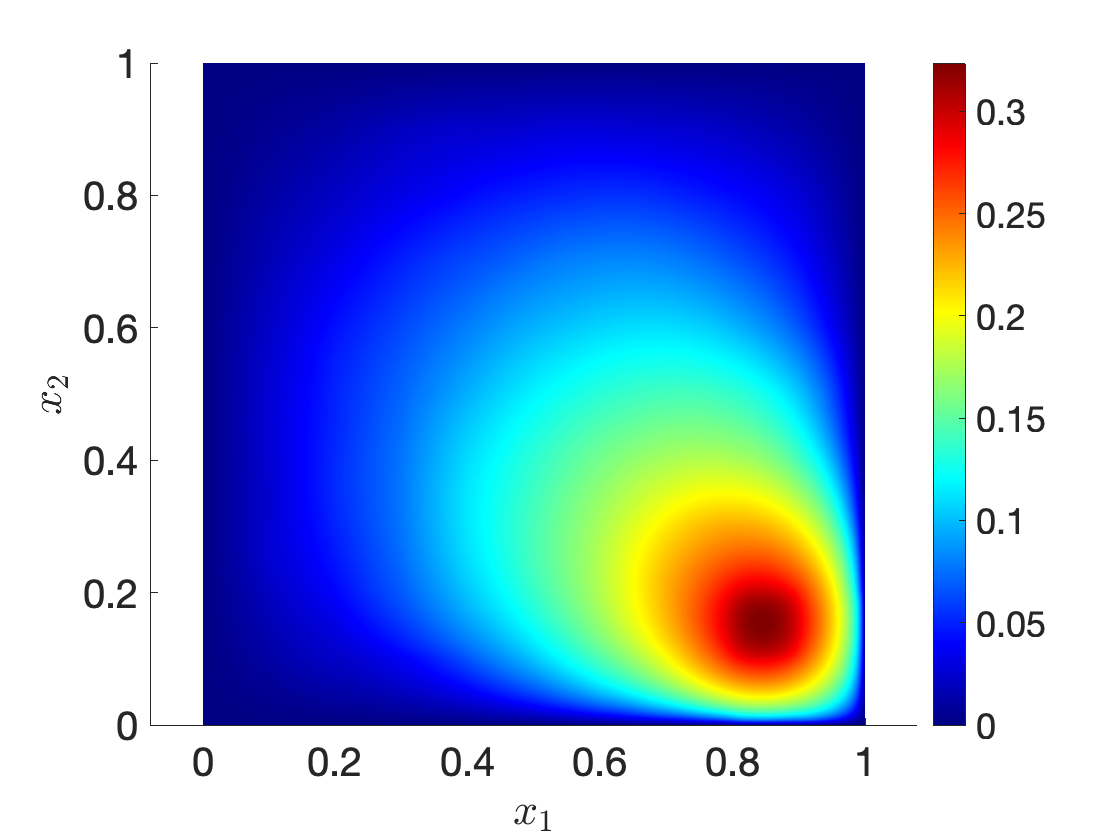}
}

\caption{nonlinear diffusion. (a) domain $\Omega$ and subdomains $\{ \Omega_i \}_i$ for $N_{\rm dd}=100$. 
(b) instantiated archetype components.
(c)-(d) solution fields for $N_{\rm dd}=100$ and two choices of the parameters.
}
\label{fig:nonlineardiff_vis}
\vspace{-15pt}
\end{figure}
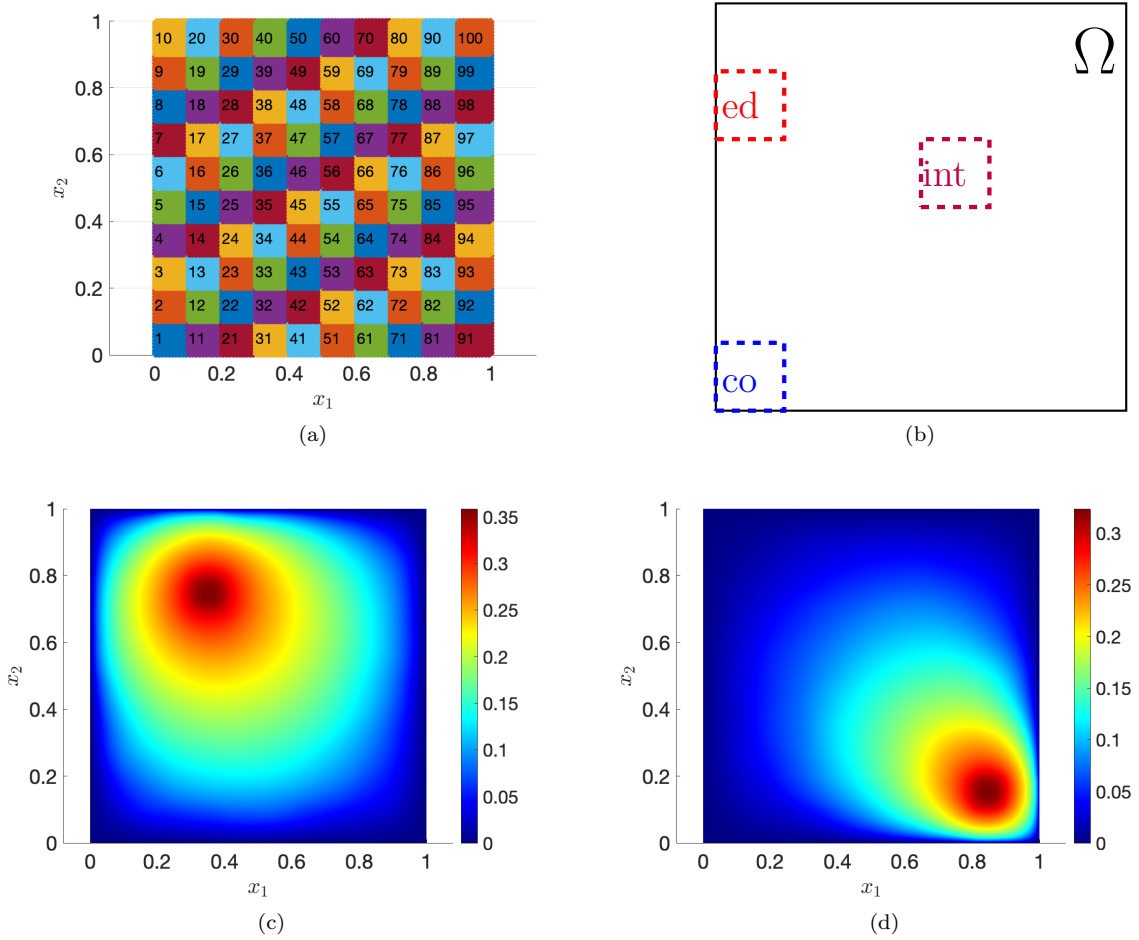

 Given $n_{\rm dd}\in \mathbb{N}$ and $H=0.1$, we define the domains
 \begin{equation}
 \label{eq:local_domains}
\Omega_{i + (j-1) n_{\rm dd} }
\, = \,
 \left\{
 [x_1 + H (i-1), x_2 + H(j-1)] \, : \,
 x_1,x_2\in (0,H)
  \right\},
 \end{equation}
 for $i,j=1,\ldots,n_{\rm dd}$, 
and $\Omega = \bigcup_{k=1}^{N_{\rm dd}} \Omega_k$ with $N_{\rm dd} = n_{\rm dd}^2$.
We further introduce 
$\widehat{\mathcal{P}} = [0.1,0.2]\times [30,40]$, the permeability coefficient
$\kappa: \Omega \times \mathbb{R} \times \bigotimes_{i=1}^{N_{\rm dd}} \widehat{\mathcal{P}} \to \mathbb{R}_+$ such that
{ $
\kappa\big|_{\Omega_i} = 
\kappa \left(x; u , \mu^{(1)},\ldots,\mu^{(N_{\rm dd})}  \right)  \big|_{\Omega_i}$ satisfies}
\begin{subequations}
\label{eq:nonlinear_diffusion}
\begin{equation}
\kappa\big|_{\Omega_i}\, = \,
\frac{36}{\mu_2^{(i)}}  \left(
\frac{u(1 - u)}{   u^3 + \frac{12}{\mu_2^{(i)}}(1-u)^3  }
\right)^2 \, + \, \mu_1^{(i)},
\quad
i=1,\ldots,N_{\rm dd};
\end{equation}
and the source term
\begin{equation}
f(x;  i^{\star}) \, = \,
100 \, e^{  - 50 \| x - x_{\rm c, i^{\star}}  \|_2^2 } \mathbbm{1}_{\Omega_i^{\star}}(x).
\end{equation}
Then, we introduce the PDE model in strong form:
given 
$\mu = [ \mu^{(1)},\ldots,\mu^{(N_{\rm dd})} , i^{\star}   ] \in 
\mathcal{P}_{\rm glo}(n_{\rm dd}): = 
\bigotimes_{i=1}^{N_{\rm dd}} \widehat{\mathcal{P}} \times \{1,\ldots,N_{\rm dd} \}$, 
find $u_{\mu}$ such that
\begin{equation}
\label{eq:PDE_model_diff}
\left\{
\begin{array}{ll}
- \nabla \cdot \left(  
\kappa_{\mu}(u_{\mu})  \, 
\nabla u_{\mu}
 \right) \, = \, f_{\mu} & {\rm in} \; \Omega, \\[3mm]
u_{\mu} = 0 & {\rm on} \; \partial  \Omega. \\
\end{array}
\right.
\end{equation}
In   \cref{fig:nonlineardiff_vis}(a)-(c)-(d), we show the domain $\Omega$ and selected snapshots for $N_{\rm dd} = 100$.
We discretize the problem using a Q3 spectral element method based on a structured grid with $961$ degrees of freedom in each subdomain $\Omega_i$.
The PDE model \eqref{eq:PDE_model_diff} has been previously considered in \cite{SO17}, and   is inspired by the model for immiscible two-phase flows in porous media studied in \cite{Michel2003}. {While we focus on \cref{eq:PDE_model_diff} in this paper to ease the exposition of ideas, we emphasize that the proposed methods can be readily applied to other nonlinear elliptic PDEs.}
\end{subequations}

In the following, we devise a component-based reduced-order model (CB-ROM) for \eqref{eq:PDE_model_diff}: the ROM should handle arbitrary choices of $N_{\rm dd}$ and of the parameters
{ $\mu^{(1)},\ldots,\mu^{(N_{\rm dd})} \in \widehat{\mathcal{P}}$ and $i^{\star}$}. 
In view of the presentation of the methodology, we introduce the overlapping partition that will be used to define the partition of unity in  \cref{sec:domain_decomposition},  
\begin{equation}
\label{eq:overlapping_partition}
 \left\{  \omega_i \right\}_{i=1}^{N_{\rm dd}},
\quad
\omega_i = \left\{
x \in \Omega: \min_{y \in \Omega_i} \| x - y \|_{\infty} < \delta_{\rm over}
\right\},
\; i=1,\ldots,N_{\rm dd}.
\end{equation}
Note that $\bigcup_i \omega_i = \Omega$. As explained in   \cref{sec:domain_decomposition}, for non-homogeneous Dirichlet or Neumann boundary conditions we shall require  $\bigcup_i \omega_i \supset \supset  \Omega$.
Here, $\delta_{\rm over}$ is the overlapping size; in the numerical experiments, we consider 
$\delta_{\rm over} = 0.1 H$.
We further introduce the archetype components that we use to describe the system: the ``corner'' (\texttt{co}) component is associated with the corner elements of the partition $\left\{  \omega_i \right\}_i$; the 
``edge'' (\texttt{ed}) component is associated with the edge elements of   
$\left\{  \omega_i \right\}_i$;
 the 
``internal'' (\texttt{int}) component is associated with the internal elements of   $\left\{  \omega_i \right\}_i$
(see   \cref{fig:nonlineardiff_vis}(b)).

We denote by $\widehat{\Omega}^{\rm co},\widehat{\Omega}^{\rm int}, \widehat{\Omega}^{\rm ed}$ the master elements associated with the three
archetype components. For edge and corner components, we denote by $\widehat{\Gamma}_{\rm dir}^{\rm ed}, \widehat{\Gamma}_{\rm dir}^{\rm co}$ the Dirichlet boundaries; furthermore, we introduce 
the local HF spaces associated with the spectral element discretization
\begin{equation}
\label{eq:local_hilbert}
\mathcal{Y}^{\rm int} \subset  H^1(\widehat{\Omega}^{\rm int} ),
\;\;
\mathcal{Y}^{\rm ed}  \subset  H_{0, \widehat{\Gamma}_{\rm dir}^{\rm ed} }^1(\widehat{\Omega}^{\rm ed} ),
\;\;
\mathcal{Y}^{\rm co}  \subset  H_{0, \widehat{\Gamma}_{\rm dir}^{\rm co} }^1(\widehat{\Omega}^{\rm co} )
\end{equation} 
and the corresponding (semi-) norms 
$\| \cdot \|_{\bullet} $ with 
$\bullet \in \{\texttt{co}, \texttt{ed}, \texttt{int}\}$
that are introduced in   \cref{sec:domain_decomposition}.
We denote by $\texttt{L} \in \{\texttt{co}, \texttt{ed}, \texttt{int}\}^{N_{\rm dd}}$ the set of labels  that link the   elements of $\left\{  \omega_i \right\}_i$ to the corresponding component; we further denote by 
$\Phi_i: \widehat{\Omega}^{\texttt{L}_i} \to \omega_i$ the mapping from the (appropriate) component to the $i$-th element of the partition.
We remark that the mappings $\Phi_i$ are simple translations for all internal components, while they are the composition of a rigid translation and a rotation for boundary (edge and corner) components 
to ensure that $\Phi_i( \widehat{\Gamma}_{\rm dir}^{\texttt{L}_i} ) \subset \partial \Omega$ and thus 
to facilitate the imposition of Dirichlet conditions. 
 
\section{Partition of unity method for localized model reduction}
\label{sec:domain_decomposition}

\subsection{Partition of unity}
In \cite{BabMel97,MelBab96}, Babu{\v{s}}ka and Melenk proposed the partition of unity method (PUM) to construct ansatz spaces with local properties. As discussed in 
\cite{BabMel97}, the PUM is designed to include a priori knowledge about the PDE in the ansatz spaces, and it provides a systematic way to construct ansatz spaces of any desired regularity. In the framework of CB-pMOR, the PUM provides a systematic  framework to construct global ansatz spaces and  offers strong theoretical guarantess concerning approximation and robustness.

Given the overlapping cover of $\Omega$,
$ \{ \omega_i \}_{i=1}^{N_{\rm dd}}$, we denote by $M$ the minimum constant such that
\begin{subequations}
\label{eq:PoU}
\begin{equation}
\label{eq:PoU_cardinality}
\forall \, x\in \Omega, 
\quad
{\rm card} \left\{  i\in \{1,\ldots,N_{\rm dd}\}: x\in \omega_i   \right\} \leq M,
\end{equation}
where ${\rm card} (A)$ denotes the cardinality of the discrete set $A$. Then, we define the partition of unity (PoU) $\{ \phi_i \}_{i=1}^{N_{\rm dd}}$ such that
\begin{equation}
\label{eq:PoU_properties}
\left\{
\begin{array}{l}
\displaystyle{
{\rm supp} \left( \phi_i \right) \subset \overline{\omega}_i , 
\;\;
0\leq \phi_i (x) \leq 1,  \;\;
\|  \nabla \phi_i \|_{L^{\infty}(\Omega)}  \leq C_i,
}
\\[3mm]
\displaystyle{\sum_{j=1}^{N_{\rm dd}} \phi_j(x) = 1,
\quad\quad\quad\quad\quad\quad
x \in \Omega, \;\; 
i=1,\ldots, N_{\rm dd}.
}
\\
\end{array}
\right.
\end{equation}
\end{subequations}
We say that  $\{ \phi_i \}_{i=1}^{N_{\rm dd}}$ is of degree $m$ if  $\{ \phi_i \}_{i=1}^{N_{\rm dd}} \subset C^m(\Omega; \mathbb{R})$. Then, we define the PUM spaces
 \begin{equation}
 \label{eq:pum_HF_space}
 \mathcal{X}_{\rm pum} :=
 \left\{
\sum_{i=1}^{N_{\rm dd}}  \phi_i   \psi_i \; : \; 
\psi_i \in \mathcal{X}_{i}
 \right\}   \subset H_0^1(\Omega),
 \end{equation}
  where $\mathcal{X}_{i} = \{  
 \zeta \circ \Phi_{i}^{-1} \, : \,
\zeta \in \mathcal{Y}^{\texttt{L}_i} 
  \}$.
 Note that by construction $
\phi_i \zeta \circ \Phi_{i}^{-1} 
  \in H_0^1(\omega_i)$ and can thus be trivially extended to $\mathbb{R}^d$. 
Next, given the reduced spaces
 $\{ \mathcal{Z}^{\bullet} \}_{\bullet} $ such that $  \mathcal{Z}^{\bullet} \subset \mathcal{Y}^{\bullet}$, we define the 
 global reduced  space
 \begin{equation}
 \label{eq:pum_space}
 \mathcal{Z}_{\rm gfem} :=
 \left\{
\sum_{i=1}^{N_{\rm dd}} \phi_i \zeta_i\circ \Phi_i^{-1} \; : \; 
\zeta_i \in \mathcal{Z}^{\texttt{L}_i}
 \right\}
 \subset  \mathcal{X}_{\rm pum}.
 \end{equation}

\Cref{thm:PUM_approximation} provides a rigorous upper bound for the approximation properties of the PUM space in $\Omega$ {---  the local approximation condition \eqref{eq:local_approximation_pum} provides the foundations for the localized data compression strategy proposed in   \cref{sec:data_compression}.} 
 
\begin{theorem}
\label{thm:PUM_approximation}
(\cite[Theorem 1]{BabMel97})
Let $u\in H_0^1(\Omega)$. Assume that there exist $\zeta_1,\ldots$, $\zeta_{N_{\rm dd}}$ such that
$\zeta_i\circ \Phi_i \in \mathcal{Z}^{\texttt{L}_i}$ and
\begin{equation}
\label{eq:local_approximation_pum}
\| u - \zeta_i  \|_{L^2(\Omega \cap  \omega_i)} \leq \epsilon_i,
\;\;
\| \nabla  u - \nabla \zeta_i  \|_{L^2(\Omega \cap  \omega_i)} \leq \epsilon_{\nabla, i},
\;\;
i=1,\ldots,N_{\rm dd},
\end{equation}
for some positive constants $\{ \epsilon_i \}_i$ and $\{  \epsilon_{\nabla, i}\}_i$.
Then, the function $u_{\rm gfem} = \sum_{i=1}^{N_{\rm dd}} \phi_i \zeta_i \in \mathcal{Z}_{\rm gfem}$ satisfies
\begin{equation}
\label{eq:global_approximations}
\left\{
\begin{array}{l}
\displaystyle{ 
\| u - u_{\rm gfem}  \|_{L^2(\Omega)} \leq \sqrt{M} \sqrt{  \sum_{i=1}^{N_{\rm dd}} \epsilon_i^2 }
};
\\[3mm]
\displaystyle{ 
\| \nabla u - \nabla u_{\rm gfem}  \|_{L^2(\Omega)} \leq \sqrt{2 M} \sqrt{  \sum_{i=1}^{N_{\rm dd}} C_i^2 \epsilon_i^2
+ \epsilon_{\nabla,i}^2
 }
}.
\\ 
\end{array}
\right.
\end{equation}
\end{theorem}

In this work, we consider a piecewise tensorized bilinear PoU
 $\{ \phi_{i+ (j-1)n_{\rm dd}} (x)= \phi_i^{\rm 1d}(x_1)\phi_j^{\rm 1d}(x_2) \}_{i,j=1}^{n_{\rm dd}}$ where
 $\{   \phi_i^{\rm 1d} \}_i$ is a PoU subordinate to the cover
 $$
 \{ \omega_i^{\rm 1d} = ((i-1) H-\delta_{\rm over}/2, i H + \delta_{\rm over}/2) \}_{i=1}^{n_{\rm dd}}.
 $$ 
 For this choice of the PoU, we have that
 $\| \frac{d}{dx} \phi_i^{\rm 1d} \|_{L^{\infty}(\Omega)}  =  \frac{1}{\delta_{\rm over}}$ and thus the constants
 $C_i$ in \eqref{eq:PoU_properties} are given by
    $C_i  =  \frac{\sqrt{2}}{\delta_{\rm over}}$ for $i=1,\ldots,N_{\rm dd}$. Note that, since we impose Dirichlet conditions on $\partial \Omega$, we can consider 
    $\bigcup_{i=1}^{N_{\rm dd}} \omega_i = \Omega$.
    Note also that the constant $M$ in 
     \eqref{eq:PoU_cardinality} is equal to four.
     
We observe that there exist 
$\widehat{\phi}^{\rm int},\widehat{\phi}^{\rm co},\widehat{\phi}^{\rm ed}$ such that
$$
\phi_i  = \widehat{\phi}^{\texttt{L}_i} \circ \Phi_i^{-1} :  i=1, \ldots, N_{\rm dd},
$$
 for any choice of $n_{\rm dd}\in \mathbb{N}$:
this is due to the particular choice of the mappings
 $\{\Phi_i \}_i$ and of the archetype components.
  For arbitrary partitions $\{\omega_i \}_i$ 
   and arbitrary mappings $\{\Phi_i \}_i$,  given $\{ \widehat{\phi}^{\bullet} \}_{\bullet}$ such that
$0\leq \widehat{\phi}^{\bullet}(x) \leq 1$ in $\mathbb{R}^d$,
$ \widehat{\phi}^{\bullet}(x) = 0$ if $x\notin \widehat{\Omega}^{\bullet}$, $\| \nabla \widehat{\phi}^{\bullet}\|_{L^2(\mathbb{R}^2)} \leq C_{\bullet}$, we can show that the set $\{ \phi_i \}_i$ such that
$$
\phi_i  =
\frac{1}{ \sum_{j=1}^{N_{\rm dd}}   \widehat{\phi}^{\texttt{L}_i} \circ \Phi_j^{-1}    } \, 
 \widehat{\phi}^{\texttt{L}_i} \circ \Phi_i^{-1},
\quad
\forall \, i=1, \ldots, N_{\rm dd},
$$
is a PoU in $\Omega$.

%\todo[inline]{KS: I found these last two sentences confusing. Maybe rephrase? I did not do it, as, as I said, I got confused when reading. ;-) }

\subsection{Discrete variational formulation and functional norms}
We introduce the local semi-norms
\begin{equation}
\label{eq:local_norm}
\|  w  \|_{\bullet}
\, = \,
\|   \widehat{\phi}^{\bullet} \, w  \|_{H^1(\widehat{\Omega}^{\bullet})},
\quad
\bullet \in \{\texttt{int}, \texttt{co}, \texttt{ed} \}.
\end{equation}
Note that for this choice of the local norms, 
since the mappings $\{\Phi_i \}_i$ are roto-translations,
if $\{  \zeta_j^{\bullet} \}_{j=1}^n$ are orthonormal bases  with respect to  $\|  \cdot  \|_{\bullet}$, then 
$\{  \phi_i \zeta_j^{\texttt{L}_i} \}_{j=1}^n$ is orthonormal in $H^1(\omega_i)$, for $i=1,\ldots,N_{\rm dd}$.
Given the spaces
$\mathcal{X}_{i,0}:=  \{  
 \phi_i \, \zeta \circ \Phi_{i}^{-1} \, : \,
\zeta \in \mathcal{Y}^{\texttt{L}_i} 
  \}$ for  $i=1,\ldots,N_{\rm dd}$, we further introduce the inner products and induced norms
\begin{equation}
\label{eq:deployed_norms}
(w,v)_{1,\omega_i}
\, = \,
\int_{\omega_i} \, \nabla w \cdot \nabla v  + w v \, dx,
\;\;
\| w \|_{1,\omega_i}
=
\sqrt{ (w,w)_{1,\omega_i}  },
\quad
w,v\in \mathcal{X}_{i,0};
\end{equation}
the global norm 
$\| w \|_{1,\Omega}
=
\sqrt{  \int_{\Omega} \, \| \nabla w \|_2^2  + w^2 \, dx   }$; and the dual norms 
\begin{equation}
\label{eq:dual_deployed_norms}
\|  f \|_{-1,\omega_i}
\, = \,
\sup_{v\in \mathcal{X}_{i,0} }
\frac{f(v)}{\| v \|_{1,\omega_i}},
\quad
 \| F \|_{-1,\Omega}
\, = \,
\sup_{v\in \mathcal{X}_{\rm pum} }
\frac{f(v)}{\| v \|_{1,\Omega}},
\quad
i=1,\ldots,N_{\rm dd},
\end{equation}
for $f\in \mathcal{X}_{i,0}'$ and
$F\in \mathcal{X}_{\rm pum}'$.

Then, we introduce the HF problem: given $\mu \in \mathcal{P}_{\rm glo}$, find 
$u_{\mu} \in \mathcal{X}_{\rm pum}$ such that
\begin{subequations}
\label{eq:variational_form}
\begin{equation}
\mathcal{G}_{\mu}(u_{\mu} ,v)  = 0 
\quad
\forall \, v\in \mathcal{X}_{\rm pum},
\end{equation}
 where
\begin{equation}
\mathcal{G}_{\mu}(w,v) := 
\int_{\Omega} \; \eta_{\mu}(x; \, w,v) \, dx 
\quad
{\rm with} \quad
\eta_{\mu}(x; w,v) \, =\,
\kappa_{\mu}(x;w) \nabla w \cdot \nabla v - f_{\mu} v
\end{equation}
\end{subequations}
and $\mathcal{G}_{\mu}: \mathcal{X}_{\rm pum} \rightarrow \mathcal{X}_{\rm pum}'$.
%\todo[inline]{KS: $\eta$ depends on $x$ but the right-hand side in the definition does not -> Added $x$ to $\kappa$. }
%{ 
%We denote by $\mathcal{G}: \mathcal{X}_{\rm pum} \times \mathcal{X}_{\rm pum}  \times \mathcal{P}_{\rm glo} \to \mathbb{R}$ 
%the variational form associated with \eqref{eq:nonlinear_diffusion},
%}

\subsection{Residual assembly and algebraic formulation of the reduced-order model}
\label{sec:residual_computation}

We omit dependence of $\Omega$ and $\mathcal{P}_{\rm glo}$ 
on $n_{\rm dd}$
to shorten notation.
%{
We consider the Galerkin ROM:
\begin{equation}
\label{eq:CB_PUM}
{\rm find} \; \;
\widehat{u}_{\mu} \in \mathcal{Z}_{\rm gfem} \, : \,
\mathcal{G}_{\mu} (\widehat{u}_{\mu} , v) = 0\quad
\forall \, v\in \mathcal{Z}_{\rm gfem}.
\end{equation}
%}
%\todo[inline]{KS: Below we use $u_{\rm gfem}$ for the GFEM Approximation?}
Given the local approximation spaces 
$\mathcal{Z}^{\rm int},\mathcal{Z}^{\rm ed},\mathcal{Z}^{\rm co}$ with bases\footnote{
Here, we choose $n = {\rm dim} ( \mathcal{Z}^{\rm int}   )= {\rm dim} ( \mathcal{Z}^{\rm ed}   )= {\rm dim} ( \mathcal{Z}^{\rm co} )$. This choice simplifies notation and is also convenient for code vectorization.
The extension to reduced spaces of arbitrary size is straightforward.
} $\{  \zeta_i^{\bullet} \}_{i=1}^n$, we define the basis of $\mathcal{Z}_{\rm gfem}$ 
$\{ \zeta_{i,j}  \}_{i,j}$ such that
\begin{subequations}
\begin{equation}
\label{eq:pum_basis}
\zeta_{i,j} = \zeta_i^{\texttt{L}_j}  \circ \Phi_i^{-1}  \, \phi_j,
\quad
i=1,\ldots,n, \;\;
j=1,\ldots,N_{\rm dd}.
\end{equation}
Given $u \in \mathcal{Z}_{\rm gfem}$, we set $N:= n N_{\rm dd}$ and we denote by $\mathbf{u} \in \mathbb{R}^N$ the vector of coefficients such that
\begin{equation}
\label{eq:field2vector}
u = \sum_{j=1}^{N_{\rm dd}} \sum_{i=1}^n \left( \mathbf{u}    \right)_{i + (j-1) n} \zeta_{i,j}.
\end{equation} 
\end{subequations}
Then, we introduce the discrete residual $\mathbf{R}: \mathbb{R}^N \times 
\mathcal{P}_{\rm glo} \to \mathbb{R}^N$ such that
\begin{subequations}
\begin{equation}
\label{eq:residual_vector} 
\left(
\mathbf{R}_{\mu}( \mathbf{u})
\right)_{i + (j-1) n}
\, = \,
\mathcal{G}_{\mu}(u, \zeta_{i,j} )
\end{equation}
and the algebraic nonlinear problem associated to \eqref{eq:CB_PUM}:
\begin{equation}
{\rm find} \; \widehat{\mathbf{u}}_{\mu}\in \mathbb{R}^N \;\; {\rm such \;\; that }
\;\;
\mathbf{R}_{\mu} \left( 
\widehat{\mathbf{u}}_{\mu} \right) = \mathbf{0}.
\end{equation}
\end{subequations}

In order to discuss the practical evaluation of the discrete residual $\mathbf{R}_{\mu}$ in  \eqref{eq:residual_vector}, we define
$\{  
\widehat{u}_{i, \mu} = \widehat{u}_{\mu} \phi_i  \}_i$ and ${\rm Neigh}_i = \{j : \omega_i \cap \omega_j \neq \emptyset \}$. Then, we observe that
\begin{subequations}
\label{eq:residual_assembly}
\begin{equation}
\label{eq:residual_assembly_a} 
\begin{array}{rl}
\displaystyle{\mathcal{G}_{\mu}(  \widehat{u}_{\mu}  ,  \zeta_{i,j} ) =    }
&
\displaystyle{
\int_{\omega_i} \; \eta_{\mu} \left(x; 
\sum_{k\in {\rm Neigh}_i} \, \widehat{u}_{k, \mu}
\, , \zeta_{i,j} \right) \, dx
} \\[3mm]
=& 
\displaystyle{
\int_{\widehat{\Omega}^{\texttt{L}_i} }    \; \widehat{\eta}_{\mu}^{(i)}
\left(x; 
\sum_{k\in {\rm Neigh}_i} \, \widehat{u}_{k, \mu} \circ \Phi_i
\, , \zeta_i^{\texttt{L}_j}  \widehat{\phi}^{\texttt{L}_j} \right)
 \, dx,
}
\\
\end{array}
\end{equation}
where
\begin{equation}
\label{eq:residual_assembly_b}
\widehat{\eta}_{\mu}^{(i)}
\left(x; 
w, \,  v \right)
=
\left(
\kappa_{\mu}(\Phi_i(x); w) \nabla \Phi_i^{-1} 
\nabla \Phi_i^{-T} \nabla w \cdot
\nabla v
\, - \, \widetilde{f}_{\mu} v
\right) \, {\rm det }(\nabla \Phi_i),
\end{equation}
with  
$\widetilde{f}_{\mu} = f_{\mu} \circ \Phi_i$.
Since $\{\Phi_i \}_i$ are roto-translations,   \eqref{eq:residual_assembly_b}  reduces to:
\begin{equation}
\label{eq:residual_assembly_c}
\widehat{\eta}_{\mu}^{(i)}
\left(x; 
w, v \right)
\, =
\,
\kappa_{\mu}(\Phi_i(x); w)  
 \nabla w \cdot \nabla v 
\, - \, \widetilde{f}_{\mu} v.
\end{equation}
\end{subequations}

We observe that the Jacobian $\mathbf{J}_{\mu}(\cdot)$ of the algebraic residual 
 $\mathbf{R}_{\mu}(\cdot)$ is sparse for large values of $N_{\rm dd}$. More precisely, exploiting \eqref{eq:residual_assembly}, it is easy to verify that the number of non-zero elements of 
$\mathbf{J}_{\mu}(\cdot)$ is bounded by
$$
\texttt{nnz} \left( \mathbf{J}_{\mu}(\mathbf{u}) \right)
\leq \sum_{i=1}^{N_{\rm dd}} n^2 {\rm card}  ( {\rm Neigh}_i )
= \mathcal{O} \left(  n^2 N_{\rm dd} \right),
\quad
\forall \, \mathbf{u} \in \mathbb{R}^N.
$$
For the model problem considered in this work we have 
 ${\rm card} \left( {\rm Neigh}_i \right) \leq 9$ for $i=1,\ldots,N_{\rm dd}$. 
 
 Assembly of the residual in \eqref{eq:residual_assembly} is extremely expensive due to the need to integrate over all instantiated components $\{   \widehat{\Omega}^{L_{i}}\}$. To speed up computations, we should thus resort to hyper-reduction techniques
 \cite{BMNP04,ChaSor10,farhat2015structure,ryckelynck2005priori,YanPat19}.
The choice of the hyper-reduction procedure strongly depends on the PDE model of interest, on the underlying high-fidelity numerical scheme,     and on the geometrical parameterization: we refer to \cite{taddei2021discretize} for a discussion on the treatment of geometry parameterizations. We further observe that evaluation of \eqref{eq:residual_assembly_a} involves evaluation of $\widehat{u}_{j,\mu}$ in the mapped quadrature points of the mesh $\widehat{\Omega}^{\texttt{L}_i}$: this evaluation is extremely expensive for unstructured meshes { and thus requires a specialized treatment}.
The development of specialized   hyper-reduction techniques for  CB-pMOR is part of ongoing research and is not addressed in the present work.

\section{Data compression: randomized localized training}
\label{sec:data_compression}
 
 To highlight the main features of our methodology without unnecessary notation, we assume here that the system is described by a single archetype component. We denote by $\mathcal{Z} \subset \mathcal{Y}$ the local approximation space. 
For the model problem in  \cref{sec:model_problem}, this would correspond to the case of Neumann or Robin boundary conditions on $\partial \Omega$. 
 In the supplementary materials, we discuss the extension to the case of multiple components and we provide further details for the particular test case considered.
 
The aim of this section is to devise an actionable procedure to build  a local approximation space 
$\mathcal{Z} \subset \mathcal{Y}$
such that
\begin{equation}
\label{eq:goal_data_compression}
\min_{\zeta \in \mathcal{Z} }
\; \| u_{\mu} \big|_{\omega_i} - \zeta\circ \Phi_i^{-1}  \|_{1,\omega_i}
\leq \varepsilon_{\rm tol}
\;\;
{\rm for} \; i=1,\ldots,N_{\rm dd},
\;\;
\mu \in \mathcal{P}_{\rm glo}(n_{\rm dd}),
\end{equation} 
where $\varepsilon_{\rm tol}>0$ is a prescribed tolerance. 
Condition \eqref{eq:goal_data_compression} implies that the local space $\mathcal{Z}$ should approximate the manifold
\begin{equation}
\label{eq:true_local_manifolds}
\mathcal{M} = \left\{
 u_{\mu} \big|_{\omega_i}  \circ \Phi_i \, : \,
 i=1,\ldots,N_{\rm dd}, \; 
 \mu\in 
\mathcal{P}_{\rm glo}(n_{\rm dd}),
\;\;
n_{\rm dd} \in \mathbb{N}
\right\} \subset \mathcal{Y}.
\end{equation}
The computation of snapshots {that belong to} the manifold
$\mathcal{M}$ requires to solve global problems and is thus unfeasible in our framework.
Instead, 
in \cref{sec:oversampling}, 
we propose to rely on oversampling to identify an actionable localized manifold  
$\widetilde{\mathcal{M}}$ for which we can compute snapshots;
then, in   \cref{sec:random_training}, we propose a randomized training algorithm to construct local approximation spaces.
  
\subsection{Oversampling}
\label{sec:oversampling}
\begin{wrapfigure}{r}{0.6\textwidth}
\centering
%\vspace{-15pt}
\includegraphics[scale=0.3]{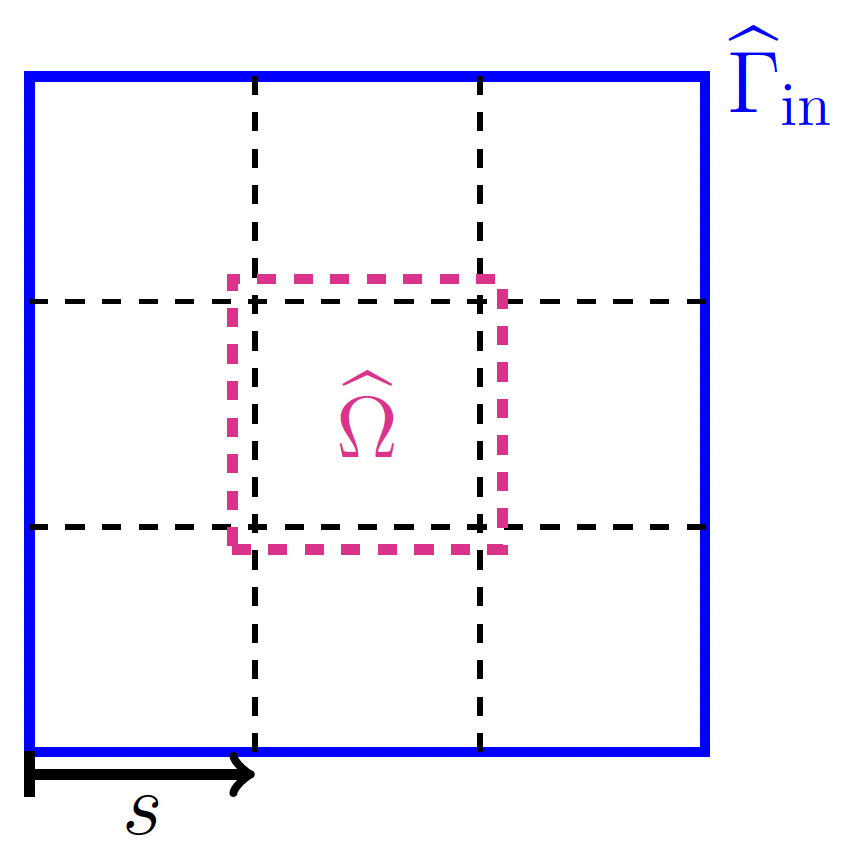}
\caption{Nonlinear diffusion. Archetype component  with corresponding oversampling domain $\widehat{U}$.}
%\vspace{-15pt}
\label{fig:nonlineardiff_vis_components_main}
\end{wrapfigure}

We fix $i\in \{1,\ldots,N_{\rm dd} \}$, and we define the patch $\widehat{U} \subset \mathbb{R}^2$ with 
input boundary $\widehat{\Gamma}_{\rm in}  \subset \partial \widehat{U}$. 
We extend the mapping $\Phi_i$ to $\widehat{U}$ and we define $U_i : = \Phi_i( \widehat{U})$ --- for the considered model problem, the mappings $\{\Phi_i \}_i$ are linear maps that can  be trivially extended to $\mathbb{R}^2$.
As depicted in   \cref{fig:nonlineardiff_vis_components_main}, we consider $U_i = \bigcup_{j \in {\rm Neigh}_i} \Omega_j$ where 
${\rm Neigh}_i$ is introduced in   \cref{sec:domain_decomposition}.  We denote by $u_{\mu,i}$ the restriction of the solution 
$u_{\mu}$ to $U_i$ and we define
$\widetilde{u}_{\mu,i} : =u_{\mu,i} \circ \Phi_i$. We observe that $\widetilde{u}_{\mu,i}$ solves the problem  (cf.  \eqref{eq:residual_assembly_c}):
\begin{subequations}
\begin{equation}
\label{eq:transfer_operator}
\int_{  \widehat{U}  }
\kappa_{\mu}(\widetilde{u}_{\mu,i}) \nabla \widetilde{u}_{\mu,i} \cdot \nabla v  \, dx \, = \,
\int_{  \widehat{U}    }
\, \tilde{f}_{\mu} \,  v  \, dx,
\quad
\forall \, v \in 
\mathcal{Y}_{i,0}^{\rm ovr},
\end{equation}
with $\widetilde{u}_{\mu,i}|_{  \widehat{\Gamma}_{\rm in}  } = u_{\mu,i} \circ \Phi_i$ and
$\mathcal{Y}_{i,0}^{\rm ovr} = 
\{ v \circ \Phi_i :
v|_{U_i} \in \mathcal{X}_{\rm pum}, 
v|_{\partial U} = 0   \} \subset H_0^1( \widehat{U} )$. We observe that  $\widetilde{u}_{\mu,i}$ is a function of the subset of parameters $\{  \mu^{(j)}\}_{j \in {\rm Neigh}_i}$, and of the index $i^{\star}$: we can then define the  active set of parameters 
$\mathcal{P}^{\rm co} =  
\bigotimes_{i=1}^{N_{\rm dd}^{\rm co}} \widehat{\mathcal{P}} \times \{1,\ldots,N_{\rm dd}^{\rm co} , 0\}$ where 
$N_{\rm dd}^{\rm co} = {\rm card} ( {\rm Neigh}_i )$ and 
$i^{\star} = 0$ means that the source term is outside the patch. The parameterization
$\mathcal{P}^{\rm co}$ is associated to the archetype component of interest and is independent of the size of the {system} (i.e., the number of subdomains $N_{\rm dd}$).
Exploiting \eqref{eq:transfer_operator}, we define the transfer operator
$T: G \times \mathcal{P}^{\rm co} \to \mathcal{Y}$ such that  
$T_{\mu}(g) = u|_{  \widehat{\Omega}} $ where 
$G\subset H^{1/2}(  \widehat{\Gamma}_{\rm in}  )$,
$u$ satisfies \eqref{eq:transfer_operator} with 
$u|_{\partial \widehat{U}   \setminus \widehat{\Gamma}_{\rm in}} = 0$ and 
$u|_{\widehat{\Gamma}_{\rm in}} = g$.
In the implementation, we replace 
$\mathcal{Y}_i^{\rm ovr}$ with a standard discretization of $ H_0^1( \widehat{U} )$.

We define the (unknown) set  $G^{\rm true} \subset H^{1/2}(  \widehat{\Gamma}_{\rm in}  )$ 
that  contains all possible restrictions of the solution field to the input boundary for all instantiated  components $\omega_i$, all parameters, and all choices of $n_{\rm dd}$; clearly, we have 
$\mathcal{M} = \{T_{\mu}(g) \, : \, g \in  G^{\rm true}, \mu \in \mathcal{P}^{\rm co} \}$. If we introduce the ``approximation''  $G$ of $G^{\rm true}$, we obtain the localized manifold:
\begin{equation}
\label{eq:fake_localized_manifold}
\widetilde{\mathcal{M}} = \left\{
T_{\mu}  (g)  \; : \;
g \in G, \;
\mu \in   \mathcal{P}^{\rm co} 
\right\}.
\end{equation}
\end{subequations}

We observe that snapshots of $\widetilde{\mathcal{M}}$ can be computed by solving local problems  in the patch $\widehat{U}$ for prescribed choices of the active parameters $\mu \in \mathcal{P}^{\rm co}$ and the boundary conditions. 
The patch $\widehat{U}$ should be significantly smaller than $\Omega$ to ensure rapid computations; at the same time, $\widehat{U}$ should be large enough to ensure decay of high-frequency modes on $\widehat{\Gamma}_{\rm in}$. 

The choice of the set of boundary conditions $G$ is of paramount importance; clearly, $G$ should  be rich enough to ensure that $\sup_{w\in \mathcal{M}} {\rm dist} (w, \widetilde{\mathcal{M}}) \leq \varepsilon_{\rm tol}$. Since the problem is nonlinear, generating a discrete representative approximation of  the high-dimensional set $G$ is also particularly challenging. In the next section, we directly 
prescribe  a probability density function (pdf) $p_{\rm bc}$ of the space of boundary conditions: the set $G$ is thus defined as the support of the pdf $p_{\rm bc}$.
  
\subsection{Randomized training}
\label{sec:random_training}
 
We introduce notation
$$
\mathcal{Z} = {\rm POD} \left(  
\{  u^{(i)} \}_{i=1}^{n_{\rm train}}, \, (\cdot, \cdot), \; n
\right)
$$ 
 to refer to the application of proper orthogonal decomposition (POD, \cite{volkwein2011model}) to the snapshot set $\{  u^{(i)} \}_{i=1}^{n_{\rm train}}$ with inner product $(\cdot,\cdot)$; here, $n$ denotes the number of POD modes in the output space 
 $\mathcal{Z}$. 
We further denote by $p_{\mu}$ and $p_{\rm bc}$ the pdfs for  parameter and boundary conditions. 

\begin{algorithm}[t]                      
\caption{Randomized localized training}     
\label{alg:localized_training}     

\small
\begin{flushleft}
\emph{Inputs:}  $n_{\rm train}$ size of training   set, 
$p_{\mu} , p_{\rm bc}$ pdfs.
\smallskip

\emph{Output:} 
$\mathcal{Z}$ local approximation spaces.
\end{flushleft}                      

 \normalsize

 \begin{algorithmic}[1]

%\For {$\bullet \in \{\texttt{co}, \texttt{ed}, \texttt{int}\}$}

\State
Generate $\mu^{(k)} \overset{\rm iid}{\sim} p_{\mu}$,
$g^{(k)} \overset{\rm iid}{\sim} p_{\rm bc}$, $k=1,\ldots,n_{\rm train}$.
\smallskip

\State
Compute $u^{k} = T_{\mu^{(k)}}  ( g^{(k)} )$  for  $k=1,\ldots,n_{\rm train}$.
\smallskip

\State
$ \mathcal{Z}  = {\rm POD} \left(  \left\{ u^{k}  \right\}_{k=1}^{n_{\rm train}} ,  (\cdot, \cdot),  n
 \right)$.
\smallskip
\end{algorithmic}
\end{algorithm}
 
  \Cref{alg:localized_training}  illustrates the randomized training procedure. The algorithm reads as a randomized POD \cite{yu2015randomized} with respect to parameter and boundary conditions: the inputs of the algorithm are the number of training  points $n_{\rm train}$, the size of the sought reduced spaces $n$, and the pdfs $\{p_{\mu}, p_{\rm bc} \}$ for the archetype component; the output is the reduced space   $  \mathcal{Z}$.
    
It is well-known that POD is optimal in $L^2( p_{\mu}  \times p_{\rm bc}   )$    in the limit $n_{\rm train}  \to \infty$; however, since the pdfs 
$p_{\mu}, p_{\rm bc}$ are chosen \emph{a priori}, they might not be representative of the 
true distributions for  the global systems. 
%Similarly, the error estimator provides a probabilistic error bound for the expected relative error under the assumption   $\mu \sim p_{\mu}$ and $g \sim p_{\rm bc}$.
Provided that additional information on the class of global systems of interest is available,  these observations motivate the enrichment strategy proposed in   \cref{sec:local_global}.

\begin{remark}
\label{remark:error_probabilistic}
\emph{Probabilistic a posteriori error estimation.}
Given 
$n_{\rm test}$ additional simulations
$\{  u^{(i)} \}_{i=1}^{n_{\rm test}} \subset \widetilde{\mathcal{M}}$ and the space 
 ${\mathcal{Z}}$, we  introduce the error indicator 
 \begin{equation}
 \label{eq:error_indicator}
\widehat{E} : =
\frac{1}{n_{\rm test}} \sum_{i=1}^{n_{\rm test}} \; \frac{ \| u^{(i)} - \Pi_{\mathcal{Z}} u^{(i)}  \|  }{\| u^{(i)}  \|},
 \end{equation}
 which  measures the average relative projection error on the test set $\{  u^{(i)}  \}_i$. Here,
 $\Pi_{\mathcal{Z}} : \mathcal{Y}  \to \mathcal{Z}$ is the projection operator on $\mathcal{Z}$.
 Provided that $u^{(i)} = T_{\mu^{(i)}} ( g^{(i)} ) $
 with $\mu^{(i)} \overset{\rm iid}{\sim} p_{\mu}$ and $g^{(i)} \overset{\rm iid}{\sim} p_{\rm bc}$, then \eqref{eq:error_indicator} is an unbiased estimator of the expected relative projection error
\begin{equation}
\label{eq:estimate_modelreduction}
E  : =
\mathbb{E}_{\mu \sim  p_{\mu}, g \sim      
p_{\rm bc}} 
\left[
\frac{\|   T_{\mu} ( g  ) - \Pi_{\mathcal{Z}} T_{\mu} ( g  )      \|      }{\| T_{\mu} ( g  )   \|}
\right].
\end{equation}
 Note that the error indicator provides a measure of the performance of $\mathcal{Z}$ for the particular choice of the sampling distribution.
 \end{remark}
 
\subsubsection{Random boundary conditions}
The oversampling domain $\widehat{U}$ in  \cref{fig:nonlineardiff_vis_components_main} 
contains $N_{\rm dd}^{\rm co} $ subdomains (cf.   \cref{fig:nonlineardiff_vis_components_main}): in absence of prior information, we propose to set 
\begin{equation}
\label{eq:parameter_sampling}
\begin{array}{l}
\displaystyle{
\mu = \left[\mu^{(1)}, \ldots, \mu^{(N_{\rm dd}^{\rm co})}, i^{\star} \right] ,
\quad
\mu^{(i)} \overset{\rm iid}{\sim} {\rm Uniform} \left(  \widehat{\mathcal{P}} 
\right),
}
\\[3mm]
\displaystyle{
{\rm Pr} \left( i^{\star} = t \right)
=\left\{
\begin{array}{ll}
\frac{p_{\rm src}}{N_{\rm dd}^{\rm co}}
&
t = 1,\ldots, N_{\rm dd}^{\rm co} \\
1- p_{\rm src} & t = 0 \\
\end{array}
\right.
}
\\
\end{array}
\end{equation}
 where $p_{\rm src}$ is the probability that a source term is present in the patch. If $N_{\rm dd}$ is known \emph{a priori}, we might set $p_{\rm s} = \frac{ N_{\rm dd}^{\rm co}  }{N_{\rm dd}}$. In this work, however, we consider $p_{\rm src} = 0.5$.

In view of the definition of   $p_{\rm bc}$, we introduce the curvilinear coordinate 
$s\in  [0,1]$ (cf.   \cref{fig:nonlineardiff_vis_components_main}); then, given $N_{\rm f} \in \mathbb{N}$ and $\alpha\in \mathbb{R}_+$, we define the complex-valued random field $\widetilde{g}$ such that
\begin{equation}
\label{eq:randomBC_basic}
\widetilde{g}(s; \mathbf{c}^{\rm re}, \mathbf{c}^{\rm im}) 
\, = \,
\sum_{k=0}^{N_{\rm f}-1} \frac{ c_{k+1}^{\rm re} + {\texttt{i}}  c_{k+1}^{\rm im}   }{\sqrt{1  + (2\pi k)^{2\alpha}}}
\;
e^{ 2 \pi k s \texttt{i} },
\quad
c_k^{\rm re}, c_k^{\rm im} \overset{\rm iid}{\sim} 
 \mathcal{N}(0,1).
\end{equation}
Recalling that for any $k,k'=0,\ldots,N_{\rm f}-1$ and $\alpha\in \mathbb{N}$, we have
$$
\begin{array}{l}
\displaystyle{
\int_0^1 \, e^{ 2 \pi k s \texttt{i} } \, e^{ - 2 \pi k' s \texttt{i} } \, ds = \delta_{k,k'},
\quad
\frac{d^{\alpha}}{ds^{\alpha}} e^{ 2 \pi k s \texttt{i}  }
=
(2 \pi k   \texttt{i})^{\alpha} e^{ 2 \pi k s \texttt{i}  };
}
\\[3mm]
\displaystyle{
(2 \pi k   \texttt{i})^{\alpha}  (-2 \pi k'   \texttt{i})^{\alpha} 
=
(4 \pi^2 k  k' )^{\alpha};
}\\
\end{array}
$$
 we find that 
 $$
 \begin{array}{l}
 \displaystyle{
 \|  \widetilde{g}(\cdot ; \mathbf{c}^{\rm re}, \mathbf{c}^{\rm im})  \|_{H^{\alpha}(0,1)}^2 
= }  
 \displaystyle{
 \|  \widetilde{g}(\cdot ; \mathbf{c}^{\rm re}, \mathbf{c}^{\rm im})  \|_{L^2(0,1)}^2 
\, + \,
\|  \widetilde{g}^{(\alpha)}(\cdot ; \mathbf{c}^{\rm re}, \mathbf{c}^{\rm im})  \|_{L^2(0,1)}^2
 } \\[3mm]
 = 
  \displaystyle{
\sum_{k=0}^{N_{\rm f}-1 }
\left(  c_{k+1}^{\rm re} + \texttt{i} c_{k+1}^{\rm im}   \right)
\left(  c_{k+1}^{\rm re} - \texttt{i} c_{k+1}^{\rm im}   \right)
=
\sum_{k=1}^{N_{\rm f}}
\left(  c_k^{\rm re}  \right)^2 + 
\left(   c_k^{\rm im}   \right)^ 2.
 } \\ 
 \end{array}
$$
 The latter implies that  the random variable 
$X: =  \|  \widetilde{g}(\cdot ; \mathbf{c}^{\rm re}, \mathbf{c}^{\rm im})  \|_{H^{\alpha}(0,1)}^2$ is distributed as  a $\chi^2$ distribution with $2 N_{\rm f}$ degrees of freedom: therefore, the parameter $\alpha$ in \eqref{eq:randomBC_basic} controls (in a {probabilistic} sense) the Sobolev regularity of the datum $\widetilde{g}$.

In order to choose $p_{\rm bc}$, exploiting a physical argument --- the solution $u_{\mu}$ to \eqref{eq:PDE_model_diff} represents water  saturation --- we anticipate that 
 $u_{\mu} \in [0, \bar{u}_{\rm max}]$ for some $\bar{u}_{\rm max}<1$.
Furthermore, we wish to devise samplers that reflect  the Sobolev regularity of the datum $g$.
   For these reasons, we propose to consider the  procedure in   \cref{alg:randomBC} to generate random samples of the boundary condition.
We first generate a sample of the random field $\widetilde{g}$ in \eqref{eq:randomBC_basic} and we extract its real part (cf. Line 3).
then, we simply rescale the datum 
to ensure that the image of $g$, ${\rm Im}[g]$, is contained in $[0,\bar{u}_{\rm max}]$.
{
The strategy in   \cref{alg:randomBC} is not well-suited for sampling of boundary conditions for the corner and edge components due to the presence of Dirichlet boundaries; we postpone the description of the full sampling strategy to \cref{suppsec:mult_comp}.
}

\begin{algorithm}[t]                      
\caption{Random sample generator of boundary conditions}     
\label{alg:randomBC}     

\small
\begin{flushleft}
\emph{Inputs:}  $N_{\rm f}, \alpha$  (cf. 
\eqref{eq:randomBC_basic}),
$\bar{u}_{\rm max} \in (0,1]$.
\smallskip

\emph{Output:} 
$g: [0,1]\to [0,1)$
boundary condition.
\end{flushleft}                      

 \normalsize 

\begin{algorithmic}[1]

 \State
Draw $\mathbf{c}^{\rm re}, \mathbf{c}^{\rm im} \in \mathbb{R}^{N_{\rm f}}$ s.t.
$c_k^{\rm re}, c_k^{\rm im} \overset{\rm iid}{\sim} 
 \mathcal{N}(0,1)$.
 \smallskip
 
\State
Draw $X_1,X_2  
\overset{\rm iid}{\sim} 
 {\rm Uniform}(0,\bar{u}_{\rm max}) 
 $, set $a = \min\{X_1,X_2\}$,  $b = \max\{X_1,X_2\}$.
 \smallskip

\State
Set  $\displaystyle{
g^{(1)} 
\, = \, {\texttt{Real}} 
\left[ 
\widetilde{g}(\cdot; \mathbf{c}^{\rm re}, \mathbf{c}^{\rm im}) \right]}.
$
 \smallskip

\State $g = a + \frac{b-a}{\max g^{(1)}  - \min g^{(1)} } \left( g^{(1)} -  \min g^{(1)}  \right)$.

\end{algorithmic}
\end{algorithm}

In the numerical experiments, we provide samples of the boundary conditions for various values of $\alpha\in \mathbb{R}_+$ and we investigate performance for the model problem considered. In particular, we discuss the impact of the choice of $\alpha$.
Note that the sampling strategy proposed in this section depends on several parameters
--- $N_{\rm f}, \alpha, \bar{u}_{\rm max}$ 
in   \cref{alg:randomBC} 
and $p_{\rm src}$ in \eqref{eq:parameter_sampling}
 ---
 that might be  difficult to tune.
  This observation justifies the use of few global reduced solves at training stage to improve performance of the CB-ROM.

\section{Basis enrichment  based on reduced global solves}
\label{sec:local_global}

In several contexts, it is possible to identify at the training stage a class of global configurations of interest. To provide a concrete reference for the model problem of   \cref{sec:model_problem}, we might be interested in solving the global PDE for 
(i) any choice of $n_{\rm dd} \in \{n_{\rm dd,LB}, \ldots,n_{\rm dd,UB}\}$ with $n_{\rm dd,LB},n_{\rm dd,UB}\in \mathbb{N}$,
(ii) any $\mu^{(i)} \in \widehat{\mathcal{P}}$,
(iii) up to $n_{\rm src}$ distinct sources.
The aim of this section is to devise a localized training procedure with adaptive global enrichment that exploits prior knowledge about the global system to enrich the local spaces. In   \cref{sec:error_indicator}, we present a residual-based error estimator that will be used to drive the enrichment strategy; in   \cref{sec:algorithm_adaptive}, we present the training procedure; in  
\cref{sec:coercive_linear}, we present an \emph{a priori} convergence result for linear coercive problems. As in   \cref{sec:data_compression}, we assume that the system is described by a single archetype component to shorten notation.

\subsection{Residual-based error estimation}
\label{sec:error_indicator}

Exploiting notation introduced in  
  \cref{sec:residual_computation}, given 
$i\in \{1,\ldots,N_{\rm dd}\}$, and $u\in H^1(\Omega)$, we define the local Riesz elements $\psi_{\mu}[u] \in  \mathcal{X}_{i,0}$ as
\begin{subequations}
\label{eq:residual_localized}
\begin{equation}
\label{eq:residual_localized_a} 
\left( \psi_{\mu}[u], v \right)_{1, \omega_i}
\, =\,
\int_{\omega_i} \widehat{\eta}_{\mu}^{(i)}
\left(
x; u , v \right) \, dx,
\quad
\forall \, v\in  \mathcal{X}_{i,0},
\end{equation}
and the dual residual
\begin{equation}
\label{eq:residual_localized_b}
\mathcal{r}_{\mu}^{(i)} [u] : =
\|  \psi_{\mu}[u]  \|_{1, \omega_i}.
\end{equation}
\end{subequations}
Next, \cref{lem:global_local_red} provides an upper bound for the global dual residual in terms of the localized dual residuals $\{ \mathcal{r}^{(i)}[\cdot]  \}_i$.
An analogous result for linear elliptic problems was proved in \cite[Proposition 5.1]{buhr2017arbilomod}.
The proof of \cref{lem:global_local_red} can be found in \cref{sec:proofs}.

\begin{lemma}\label{lem:global_local_red}
\label{th:residual_estimator}
Let $\{ \phi_i \}_{i}$ be a PoU that satisfies \eqref{eq:PoU}. Then, given $u\in \mathcal{X}_{\rm pum}$, we have
\begin{equation}
\label{eq:residual_estimator}
\|  \mathcal{G}_{\mu}(u, \cdot  )  \|_{-1, \Omega}
\leq
\sqrt{M}
\left(
\max_{i=1,\ldots,N_{\rm dd}} C_i^{\rm r} \,
\right)
\sqrt{  \sum_{i=1}^{N_{\rm  dd}} \left( \mathcal{r}_{\mu}^{(i)} [u] \right)^2    }
=:
\mathcal{R}_{\mu} [u],
\end{equation}
with 
$C_i^{\rm r} : =  \sqrt{ \max\{ C_i + C_i^2 + 1, 2 \}   } $.
\end{lemma}
We will employ the local residuals \eqref{eq:residual_localized_b} to mark instantiated components of the partitions where the error is large; see \cref{sec:algorithm_adaptive}. {Let us also note that as the inifinite-dimensional analogon of $\mathcal{G}_{\mu}$ as a map from $H^{1}_{0}(\Omega)$ to $H^{-1}(\Omega)$ is not in $C^{1}$, one cannot expect that the $\| \cdot \|_{-1,\Omega}$-norm of the residual (see \cref{eq:dual_deployed_norms} for the definitions) stays bounded if the mesh size goes to zero. As a remedy one may consider $\mathcal{G}_{\mu}$ as a mapping from $W^{1,p}_{0}(\Omega)$ to $W^{-1,p}(\Omega)$, $p>2$; see \cite{CalRap1997,PouRap94,SO17} and \cref{sec:proofs}. As this significantly complicates the calculations of the dual norms, we opt here for assuming that the dimension of the HF space is fixed and considering the $\| \cdot \|_{-1,\Omega}$-norm. We may then define the error indicator 
\begin{equation}
\label{eq:residual_adaptive}
\Delta_{\mu} 
=
\sqrt{
\sum_{i=1}^{N_{\rm dd,\mu}}
\, \left( \mathcal{r}_{\mu}^i \right)^2
}.
\end{equation}
For linear problems it is straightforward to derive a rigorous a posteriori bound based on $\mathcal{R}_{\mu} [\cdot]$ (see e.g. \cite{BabMel97,Buhr2020localized}). Here, we combine \cref{th:residual_estimator} with the Brezzi-Rappaz-Raviart (BRR) theory \cite{BrRaRa80,CalRap1997} to derive a rigorous residual-based error bound for the global error; see in particular \cite{CaToUr09,VerPat05} for the application of the BRR theory in the context of model order reduction. To that end, we require that
\begin{align}
\label{inf-sup-apost-H1} 0 < \beta_{2,p} &:= \underset{|w|_{1,\Omega} \neq 0}{\underset{w \in \mathcal{X}_{\rm pum}}{\inf}} \underset{|v|_{1,\Omega} \neq 0}{\underset{v \in \mathcal{X}_{\rm pum}}{\sup}} \frac{\langle \mathcal{G}_{\mu}^{\prime}(\widehat{u}_{\mu}) w, v \rangle}{| w|_{1,\Omega} | v|_{1,\Omega}},
\end{align}
and that there exist constants $\gamma_{2,p}$ and $L_{2,p}$ such that
\begin{align}
\label{continuity-apost-H1} \langle \mathcal{G}_{\mu}^{\prime}(\widehat{u}_{\mu}) w, v \rangle &\leq \gamma_{2,p} | w|_{1,\Omega} | v|_{1,\Omega}, \\
\label{Lipschitz-apost-H1}
\| \mathcal{G}_{\mu}^{\prime}(\widehat{u}_{\mu}) - \mathcal{G}_{\mu}^{\prime}(w) \| &\leq L_{2,p}\, |\,\widehat{u}_{\mu} - w\,|_{W^{1,p}(\Omega)}
\end{align}
for $w \in B(\widehat{u}_{\mu},R) \subset \mathcal{X}_{\rm pum}$ and $v \in \mathcal{X}_{\rm pum}$. Here, $R$ is supposed to be sufficiently large and $| w |_{1,\Omega}:=\|\nabla w \|_{L^{2}(\Omega)}$ and $|w|_{W^{1,p}(\Omega)}:=\|\nabla w \|_{L^{p}(\Omega)}$. To obtain a proximity indicator \cite{VerPat05,CaToUr09}, which is based on localized and easily computable residuals via the Riesz representation, we employ, as in \cite{SO17}, the finite dimensionality of $\mathcal{X}_{\rm pum}$ and define $c_{h}:=\sup_{v \in \mathcal{X}_{\rm pum}}\frac{|v|_{W^{1,p}(\Omega)}}{|v|_{1,\Omega}}$ and
\begin{equation}\label{proximity indicator}
\tau_{\mu,p} := \frac{2L_{2,p}c_{h}}{\beta_{2,p}^{2}} \sqrt{M}\left(\max_{i=1,\ldots,N_{\rm dd}} C_i^{\rm r} \,\right) \sqrt{  \sum_{i=1}^{N_{\rm  dd}} \left( \mathcal{r}_{\mu}^{(i)} [u] \right)^2    }.
\end{equation}
The proximity indicator $\tau_{\mu,p}$ will be used to validate whether $\widehat{u}_{\mu}$ is close enough to $u_{\mu}$ within the adaptive \cref{alg:localized_enrichment}. We obtain the following result, which is proved in \cref{sec:proofs}.
\begin{proposition}[Global a posteriori error bound]\label{prop: a post est}
Let $\tau_{\mu,p} < 1$ and \eqref{inf-sup-apost-H1}, \eqref{continuity-apost-H1} and \eqref{Lipschitz-apost-H1} be fulfilled. Then there exists a unique solution $u_{\mu} \in B(\widehat{u}_{\mu},\frac{\beta_{2,p}}{L_{2,p}c_{h}}) \subset \mathcal{X}_{\rm pum}$ of \eqref{eq:variational_form} and the error estimator
\begin{equation}\label{delta_BBR}
\Delta_{\mu,p} := \frac{\beta_{2,p}}{L_{2,p}c_{h}}(1 - \sqrt{1 - \tau_{\mu,p}})
\end{equation}
satisfies
\begin{equation}\label{eq:err_est}
|\widehat{u}_{\mu} - u_{\mu} |_{1,\Omega} \leq \Delta_{\mu, p}.
\end{equation}
\end{proposition}
\begin{remark}[Discussion of result]
It is well-known that for nonlinear PDEs the dual norm of the residual can only be used as an a posteriori error estimator if the approximation is already close to the high-fidelity solution (see e.g. \cite{CalRap1997,Verfuerth94,Ortner09}). Relying solely on the dual norm of the residual can therefore be problematic as it may seem that the approximation error is acceptable even though that might not be the case. The proximity indicator $\tau_{\mu,p}$ \cref{proximity indicator}, which only relies on computable constants, can be used to assess, whether indeed the approximation $\widehat{u}_{\mu}$ is close enough to $u_{\mu}$ such that the error estimation \cref{eq:err_est} is valid. While the proximity indicator $\tau_{\mu,p}$ \cref{proximity indicator} and thus the a posteriori error estimator \cref{delta_BBR} solely rely on the dual norms of local residuals that can be computed on the components and therefore do not require any global solutions, the constants $L_{2,p}$ and $\beta_{2,p}$ are global constants. We will discuss some strategies on how to estimate these constants in the next \cref{rem:est_const}. To the best of our knowledge even for linear elliptic PDEs there are no results in the conforming setting that solely rely on local constants (the a posteriori error estimators in \cite{buhr2017arbilomod,Sme15} e.g. both contain the global coercivity constant). A fully localizable a posteriori error estimator for nonlinear non-monotone PDEs would therefore be at least a paper on its own and is thus beyond the scope of this paper. 
\end{remark}
\begin{remark}[Estimation of constants]\label{rem:est_const}
Regarding the estimation of the constant $c_{h}$ in the inverse inequality, we refer to classical results e.g. in \cite{Ern2004} noting that the global inverse inequality only requires the measure of $\Omega$. Estimating the constant $L_{2,p}$ relies on estimates of the constant in the Poincar\'{e} inequality for $L^{p}$, $W^{1,p}$ and the Sobolev embedding inequality $\|v\|_{C^{0}(\Omega)} \leq c_{E} |v|_{W^{1,p}(\Omega)}$ (see e.g. \cite[Subsection 3.1.2]{Sme13}). The estimation of $c_{E}$ can be easily localized. An estimate of the constant in the Poincar\'{e} inequality involving the measure of $\Omega$ can be found in \cite[(7.44)]{GilTru77} for functions that are zero on $\partial\Omega$. We hope that if the local reduced bases contain the constant function it is maybe possible to obtain localized and more precise estimates of the Poincar\'{e} constant. Finally, similar as in \cite{Sme15}, we propose to use a localized model order reduction approximation of $\beta_{2,p}$. In detail, we suggest using the following heuristic and hierarchical estimator
\begin{equation*}
\beta_{2,p}^{app} := \underset{|w|_{1,\Omega} \neq 0}{\underset{w \in \widetilde{\mathcal{Z}}_{\rm gfem}}{\inf}} \underset{|v|_{1,\Omega} \neq 0}{\underset{v \in \widetilde{\mathcal{Z}}_{\rm gfem}}{\sup}} \frac{\langle \mathcal{G}_{\mu}^{\prime}(\widehat{u}_{\mu}) w, v \rangle}{| w|_{1,\Omega} | v|_{1,\Omega}},
\end{equation*}
where $\mathcal{Z}_{\rm gfem} \subsetneq \widetilde{\mathcal{Z}}_{\rm gfem} \subset \mathcal{X}_{\rm pum}$. We conjecture that using a certain number of additional local basis functions per component might already yield an acceptable estimate of $\beta_{2,p}$. 
\end{remark}
}

\subsection{Adaptive algorithm}
\label{sec:algorithm_adaptive}

We define the pdfs $ p_{\mu}, p_{\rm bc} $ for localized sampling and the pdf  $p_{\mu}^{\rm glo}$ that is used to generate global problems.
In the numerical examples, we consider $n_{\rm dd} \sim {\rm Uniform} (\{4,\ldots,12\} )$, $\mu^{(i)} \overset{\rm iid}{\sim}  {\rm Uniform} (\widehat{\mathcal{P}})$ and we assume that exactly one source term is active in $\Omega$ (that is, $n_{\rm src}=1$).
Given the partition 
$\{ \omega_i \}_i$, 
we define the local solution operators
%\footnote{
%We observe that 
%$\mathcal{X}_i \cap H_0^1(\omega_i) \neq  \mathcal{X}_{i,0}$.
%}
\begin{equation}
\label{eq:local_solution}
T_{\mu}^{(i)}: \mathcal{X}_{i} \to 
\mathcal{X}_{i,0}
 \;\;\;\;
{\rm s.t.}
\;\;
\mathcal{G}_{\mu} \left( u + T_{\mu}^{(i)}(u), v \right)
= 0 \quad\forall \, v\in \mathcal{X}_{i,0},
\end{equation}
for $i=1,\ldots,N_{\rm dd, \mu}$.
The particular choice of the operators
$\{ T_{\mu}^{(i)} \}_i$ is motivated by the convergence analysis in   \cref{sec:coercive_linear}.

\Cref{alg:localized_enrichment} contains the data compression procedure.
First, we initialize the local spaces using   \cref{alg:localized_training}. 
Then, 
we sample $n_{\rm train}^{\rm glo}$ configurations
$
\mathcal{P}_{\rm train}  = \{  \mu^{(k)} \}_{k=1}^{n_{\rm train}^{\rm glo}}$ with 
$\mu^{(k)} \overset{\rm iid}{\sim} p_{\mu}^{\rm glo}
$, and we  proceed with the enrichment iterations.
At the $\ell$-th iteration, 
for each $\mu \in \mathcal{P}_{\rm train}$, we resort to the CB-ROM proposed in   \cref{sec:domain_decomposition} to estimate the solution $\widehat{u}_{\mu}$; then,  we compute the local residuals
\eqref{eq:residual_localized}
$$
 \mathfrak{r}_{\mu}^{i} =  
\mathfrak{r}_{\mu}^{(i)} [  \widehat{u}_{\mu} ] ,
\quad
i=1,\ldots,N_{\rm dd,\mu} ,
$$ 
and we mark the $m_{\rm r} \%$ instantiated components with the largest residual,
   $\texttt{I}_{\rm mark}^{\mu}\subset \{1,\ldots,N_{\rm dd,\mu} \}$ . Then, we solve \eqref{eq:local_solution} to obtain 
$$
u_{i,\mu} =
\frac{1}{\phi_i}
 T_{\mu}^{(i)} (  \widehat{u}_{\mu} |_{\omega_i}),
\;\;
\forall \, i \in  \texttt{I}_{\rm mark}^{\mu},
$$
  and we update the dataset of  simulations associated with the marked elements
$\mathcal{D}$,
$$
\mathcal{D}  = \mathcal{D}  \cup \{
u_{i,\mu} \circ \Phi_{i} : 
i \in \texttt{I}_{\rm mark}^{\mu}   \} .
$$
Note that 
$u_{i,\mu}$ is not well-defined on $\partial \omega_i$ (i.e., division of $0$ by $0$): however, since we are ultimately interested in the PUM space $\mathcal{Z}_{\rm gfem}$ \eqref{eq:pum_space} and due to the choice of the local norm $\|\cdot \|$
(cf. \eqref{eq:local_norm}), this issue does not affect our procedure.
{In view of the termination condition, we further compute the error estimator $\Delta_{\mu,p}$ \cref{delta_BBR} with approximate constants.} At the end of the loop over the parameters, we update the local space  using POD (cf.  \cref{sec:random_training})
 $$
 \mathcal{Z}  = 
\mathcal{Z} 
\cup
{\rm POD}\left( \{ w - \Pi_{\mathcal{Z}} w : w\in \mathcal{D} \}, (\cdot,\cdot), n^{\rm glo}       \right),
 $$
and we check if {
$\max_{ \mu \in \mathcal{P}_{\rm train} } \Delta_{\mu,p}$} is below a user-defined tolerance.

\begin{algorithm}[t]                      
\caption{randomized localized training with global enrichment}     
\label{alg:localized_enrichment}     

\small
\begin{flushleft}
\emph{Inputs (localized training):}  
$n_{\rm train}^{\rm loc}$ = number of  solves,
 $n^{\rm loc}$ = size of the POD spaces,
 $ p_{\mu}, p_{\rm bc} $
 sampling pdfs.
 \medskip
 
 \emph{Inputs (enrichment):}  
$n_{\rm train}^{\rm glo}$ = number of global simulations per iteration,
 $n^{\rm glo}$ = number of   modes added at each iteration,
 $\texttt{maxit}$ = maximum number of outer loop iterations,
 $tol$ = tolerance for termination criterion,
 $p_{\mu}^{\rm glo}$ = global configuration sampler,
 $m_{\rm r}$ = percentage of marked components at each iteration.
 \medskip
 
\emph{Outputs:} 
$  \mathcal{Z}  $ local approximation space.
\end{flushleft}                      

 \normalsize 

\begin{center}
\textbf{Localized training}
\end{center}

\begin{algorithmic}[1]
\State
Apply   \cref{alg:localized_training} to obtain the local space 
$\mathcal{Z}$.
\end{algorithmic} 

 \begin{center}
\textbf{Enrichment}
\end{center}

\begin{algorithmic}[1]

\State Sample $n_{\rm train}^{\rm glo}$ configurations
 $\mu^{(k)} \overset{\rm iid}{\sim} p_{\mu}^{\rm glo}$,
$\mathcal{P}_{\rm train} : =\{ \mu^{(k)}  \}_k$
\smallskip

\For{$\ell=1,\ldots,\texttt{maxit}$}

\State Initialize the dataset  
$  \mathcal{D} = \emptyset$.
\smallskip

\For{ $\mu \in \mathcal{P}_{\rm train}$}

\State
Compute $\widehat{u}_{\mu}$ using the PUM-CB-ROM (cf.  \cref{sec:domain_decomposition}).
\smallskip

\State
Compute local residuals 
\eqref{eq:residual_localized}
$ \mathfrak{r}_{\mu}^{i} =  
\mathfrak{r}_{\mu}^{(i)} [  \widehat{u}_{\mu} ] $ 
for 
$i=1,\ldots,N_{\rm dd,\mu}$.
\smallskip

\State
Mark the $m_{\rm r}$ $\%$ instantiated components 
with the largest  residuals,
$\{  \omega_i \}_{i \in \texttt{I}_{\rm mark}^{\mu}   }$.

\smallskip

\State
Solve the local problems
\eqref{eq:local_solution}
in
$\{  \omega_i \}_{i \in \texttt{I}_{\rm mark}^{\mu}  }$,
$u_{i,\mu}  = \frac{1}{\phi_i} T_{\mu}^{(i)} (  \widehat{u}_{\mu} |_{ \omega_i})$.
\smallskip

\State
Augment the dataset $\mathcal{D} = \mathcal{D}  \cup \{
u_{i,\mu} \circ \Phi_{i} : 
i \in \texttt{I}_{\rm mark}^{\mu}   \} $.
\smallskip

\State
{Compute $\Delta_{\mu,p}$ \cref{delta_BBR} with approximate constants}.

\EndFor 

\State
Update the POD space 
$\mathcal{Z}  = 
\mathcal{Z} 
\cup
{\rm POD}\left( \{ w - \Pi_{\mathcal{Z}} w : w\in \mathcal{D} \}, (\cdot,\cdot), n^{\rm glo}       \right)$.

\If{ $\max_{ \mu \in \mathcal{P}_{\rm train}} \Delta_{\mu,p}  
 < tol$ }

\State  \texttt{BREAK}
 \smallskip
 
\EndIf

\EndFor 
\end{algorithmic}
\end{algorithm}

Several steps of the Algorithm are embarrassingly parallelizable:
the loop over the configurations (cf. Lines 4 to 11),
the computation of the residuals (cf. Line 6),
the solution to the local problems (cf. Line 8).
Note also that the solution to \eqref{eq:local_solution} is performed over the domain $\omega_i$ (or equivalently $\widehat{\Omega} \subset \widehat{U}$), and the Newton solver can be initialized with
the null solution: for this reason, it is significantly cheaper than the solution to \eqref{eq:transfer_operator}.
As discussed in the introduction, the enrichment algorithm
is closely linked to the online enrichment strategy proposed in \cite{OhlSch15} and also to related approaches in the multiscale FE literature
\cite{chung2018fast}.

\subsection{A priori convergence analysis for coercive linear problems}
\label{sec:coercive_linear}

We study the in-sample \emph{a priori} convergence of the enrichment procedure in  \cref{alg:localized_enrichment} : we consider the case of  linear coercive problems, and we apply the  simplified randomized procedure contained in   \cref{alg:localized_enrichment_simplified}; the proof follows the argument of \cite[Theorem 1]{buhr2018exponential}.
We assume $\mathcal{P}_{\rm train}=\{\mu\}$ and we omit dependency on $\mu$; in \ref{sec:proofs}, we discuss the extension to multiple configurations.
Then, we define the model problem:
\begin{equation}
\label{eq:coercivePDE}
{\rm find} \; u\in \mathcal{X} \, : \,
\mathcal{G}(u, v) = f(v) - a(u,v) = 0 \quad
\forall \, v\in \mathcal{X},
\end{equation}
where $H_0^1(\Omega) \subset \mathcal{X} \subset H^1(\Omega) $ is a suitable Hilbert space on $\Omega$.
We also  introduce the energy norm and the associated dual norm:
\begin{equation}
\label{eq:energy_norm}
\| w  \|_a
=
\sqrt{a(w,w)} 
\quad \forall \, w\in \mathcal{X},
\quad
\| f  \|_{\mathcal{X}'}
=
\sup_{ v\in   \mathcal{X}}
\frac{f(v)}{\| v \|_a}
\quad \forall \, f \in \mathcal{X}'.
\end{equation}
Given the partition $\{ \omega_i \}_{i=1}^{N_{\rm dd}}$, we further define the associated mappings
$\{ \Phi_i \}_{i=1}^{N_{\rm dd}}$, 
 the associated PoU  $\{ \phi_i \}_{i=1}^{N_{\rm dd}}$, 
and the local spaces
$\mathcal{X}_i = H^1(\omega_i) \cap \mathcal{X}$ and 
 $\mathcal{X}_{i,0} = H_0^1(\omega_i)$.
Then,   we define the local dual residual norms 
such that
\begin{equation}
\label{eq:loc_residuals}
\mathcal{r}^{(i)}[u]
= 
\sup_{v \in \mathcal{X}_{i,0}} \frac{\mathcal{G}(u,v)}{\| v \|_a},
\quad
i=1,\ldots,N_{\rm dd}.
\end{equation}  
Finally, we denote by $c_{\rm pu}$ the constant such that (see \eqref{eq:residual_estimator}):
\begin{equation}
\label{eq:cpu_coercive}
\|  \mathcal{G}_{\mu}(u, \cdot  )  \|_{\mathcal{X}'}
\leq
c_{\rm pu}
 \, 
\sqrt{  \sum_{i=1}^{N_{\rm  dd}} \left( \mathcal{r}^{(i)} [u] \right)^2}.
\end{equation}

\begin{algorithm}[t]                      
\caption{simplified randomized localized training with global enrichment}     
\label{alg:localized_enrichment_simplified}     
 
\begin{algorithmic}[1]

\State Initialize $\mathcal{Z} = \mathcal{Z}_0$.

\State Sample $n_{\rm train}^{\rm glo}=1$ configurations
 $\mu  \sim  p_{\mu}^{\rm glo}$,
$\mathcal{P}_{\rm train} : =\{ \mu  \}_k$
\smallskip

\For{$\ell=0,\ldots,\texttt{maxit}$}

\State
Compute $\widehat{u}_{\ell}$ using the PUM-CB-ROM (cf.   \cref{sec:domain_decomposition}).
\smallskip

\State
Find $k = {\rm arg} \max_{i=1,\ldots,N_{\rm dd}} \mathfrak{r}^{(i)} [\widehat{u}_{\ell}]$.
\smallskip

\State
Solve the local problem: find 
$u_{k}^{\rm loc}\in \mathcal{X}_{k,0}$ such that
$\mathcal{G}(\widehat{u}_{\ell}+u_{k}^{\rm loc} , v  ) = 0$ for all $v\in \mathcal{X}_{k,0}$.
\smallskip

\State
Define $u^{\star} = \frac{u_{k}^{\rm loc}}{\phi_k} $ and update the local space 
$\mathcal{Z}  = 
\mathcal{Z}  \cup {\rm span} \{  u^{\star}   \circ \Phi_k    \} $.
 
 \EndFor 
\end{algorithmic}
\end{algorithm}

 \Cref{th:exponential_convergence} shows that the reconstruction error decreases exponentially with respect to the iteration count $\ell$ for any choice of the initial reduced space. 
%The result can be  extended to the case of multiple configurations ($n_{\rm train}^{\rm glo}>1$): we omit the details.

\begin{proposition}
\label{th:exponential_convergence}
The sequence of PUM-CB-ROM solutions $\{ \widehat{u}_{\ell} \}_{\ell=1,2,\ldots}$ satisfies
$
\| u - \widehat{u}_{\ell}  \|_a \leq
\left({1 - \frac{1}{N_{\rm dd} c_{\rm pu}^2  }}
\right)^{\ell/2} \; 
 \| u - \widehat{u}_{0}  \|_a.
$
\end{proposition}

Next Lemma summarizes two standard results that will be used in the proof of    \cref{th:exponential_convergence}.
\begin{lemma}
\label{th:textbook_coercive}
Let $\mathcal{Z}_{\rm gfem} \subset \mathcal{X}$ and let $\widehat{u} \in \mathcal{Z}_{\rm gfem}$ satisfy
$\mathcal{G}(\widehat{u} , v) = 0$ $\forall \; v\in \mathcal{Z}_{\rm gfem}$.
Then, we have
\begin{subequations}
\label{eq:galerkin_optimality}
\begin{equation}
\label{eq:galerkin_optimality_a}
\|\widehat{u}  - u  \|_a = \inf_{\varphi\in \mathcal{Z}_{\rm gfem}} \| \varphi - u    \|_a;
\end{equation}
\vspace{-10pt}
\begin{equation}
\label{eq:galerkin_optimality_b}
\|\widehat{u}  - u  \|_a =   \| \mathcal{G}(u, \cdot)    \|_{\mathcal{X}'}.
\end{equation}
\end{subequations}
\end{lemma}

\begin{proof}
(\emph{\cref{th:exponential_convergence}}).
Exploiting 
\eqref{eq:cpu_coercive}
 and  then \eqref{eq:galerkin_optimality_b}, we find
\begin{equation}
\label{eq:aux_a}
 \left(
 \mathcal{r}^{(k)}[ \widehat{u}_{\ell}  ]
  \right)^2
 \geq
 \frac{1}{N_{\rm dd}}
 \sum_{j=1}^{ N_{\rm dd}   }
  \left(
 \mathcal{r}^{(j)}[ \widehat{u}_{\ell}  ]
  \right)^2\geq
  \frac{1}{ N_{\rm dd} c_{\rm pu}^2 }
  \; 
  \|   \mathcal{G}(\widehat{u}_{\ell} , \cdot)   \|_{\mathcal{X}'}^2
  =
   \frac{1}{ N_{\rm dd} c_{\rm pu}^2 }
  \; 
  \|u - \widehat{u}_{\ell}  \|_a^2.
\end{equation}
By construction, $u_k^{\rm loc}$ in   \cref{alg:localized_enrichment_simplified} belongs to $\mathcal{Z}_{\rm gfem}^{(\ell+1)}$. As a result, if we consider $\varphi =\widehat{u}_{\ell} +  u_k^{\rm loc}$ in \eqref{eq:galerkin_optimality_a}, we find
$$
\|u - \widehat{u}_{\ell+1}   \|_a^2
\leq
\|u - \widehat{u}_{\ell} - u_k^{\rm loc}  \|_a^2
=
\|u - \widehat{u}_{\ell}    \|_a^2
- 2 
a \left(  u - \widehat{u}_{\ell} ,  u_k^{\rm loc}  \right)
+
\| u_k  \|_a^2.
$$
Since
$$
a \left(  u - \widehat{u}_{\ell} ,  u_k^{\rm loc}  \right)
%=
%a \left(  u - \widehat{u}_{\ell} - u_k^{\rm loc} ,  u_k^{\rm loc}  \right)
%+
%a \left(   u_k^{\rm loc},  u_k^{\rm loc}  \right)
=
\underbrace{  
 \mathcal{G} \left(    \widehat{u}_{\ell}  + u_k^{\rm loc}  ,  u_k^{\rm loc}  \right)}_{ = 0 }
+
\|  u_k^{\rm loc}  \|_a^2
=
\|  u_k^{\rm loc}  \|_a^2,
$$
and
$$
\|  u_k^{\rm loc}  \|_a 
=
\sup_{v\in \mathcal{X}_{k,0}   }
\frac{a( u_k^{\rm loc} , v  )   }{\|  v \|_a}
=
\sup_{v\in \mathcal{X}_{k,0}   }
\frac{  \mathcal{G}( u - \widehat{u}_{\ell} , v  )   }{\|  v \|_a}
=
\mathcal{r}^{(k)}[ \widehat{u}_{\ell}  ],
$$
we obtain
\begin{equation}
\label{eq:aux_b}
\|u - \widehat{u}_{\ell+1}   \|_a^2
\leq
\|u - \widehat{u}_{\ell}   \|_a^2
-      \left(  \mathcal{r}^{(k)}[ \widehat{u}_{\ell}  ] \right)^2.
\end{equation}
Thesis follows by combining \eqref{eq:aux_a} and \eqref{eq:aux_b}.
\end{proof}

\section{Numerical results}
\label{sec:numerics}

\subsection{Performance of randomized training for a linear problem}
\label{sec:linear}
We first provide numerical investigations for the linear advection-diffusion-reaction problem
$$
\left\{
\begin{array}{ll}
-\nabla \cdot \left(
\mu_1 \kappa \nabla u_{\mu,g} + [\mu_2,\mu_3]^T u_{\mu,g}
\right) \, + \mu_4 u_{\mu,g} \, = \, 0
& {\rm in} \; U = (0,0.3)^2, \\[3mm]
u_{\mu,g} = g 
& {\rm on} \; \partial U =: \Gamma_{\rm in}, \\
\end{array}
\right.
$$
where
$\kappa (x) = \frac{1}{1 + \| x \|_2^2}$ and
$\mu  = [\mu_1,\mu_2,\mu_3,\mu_4] \in \mathcal{P} = [0.2,1]\times [-1,1]^2 \times [0,1]$.
We consider the extracted domain $\widehat{\Omega} = (0.1,0.2)^2$. The linear problem allows us to compare our randomized method with a previously-developed data compression algorithm. Note that the transfer operator
$T: (\mu, g) \mapsto u_{\mu,g}|_{ \widehat{\Omega}}$ is nonlinear due to the presence of parameters.
We discretize the problem using the finite element method based on cubic (P3) polynomials, with $N_{\rm in} = 360$ degrees of freedom on the boundary $\Gamma_{\rm in}$.

We compare performance of our randomized algorithm with the approach in \cite{TadPat18} (TE+POD): 
given the training set $\mathcal{P}_{\rm train} = \{ \mu^k \}_{k=1}^{n_{\rm train}} \subset \mathcal{P}$, we first solve $n_{\rm train}$ independent transfer eigenproblems 
\cite{BabLip11}
for each value of the parameter and then we use POD to combine the resulting spaces.
We refer to \cite{TadPat18}  for further details and analysis, and we refer to
\cite{SmePat16} for a similar data compression algorithm.
In the numerical experiments, we set $n_{\rm train}=100$: this implies that TE+POD requires to solve $n_{\rm train} \cdot N_{\rm in} = 36000$ PDEs. We envision that the total number of PDE solves can be reduced up to 
$\mathcal{O}( n \cdot n_{\rm train} ) $ by resorting to Krylov methods   to solve the transfer eigenproblem: we refer to the above-mentioned literature for further details.

\begin{figure}[t]
\centering

\subfloat[$\alpha=1$]{
\includegraphics[width=0.3\textwidth]
{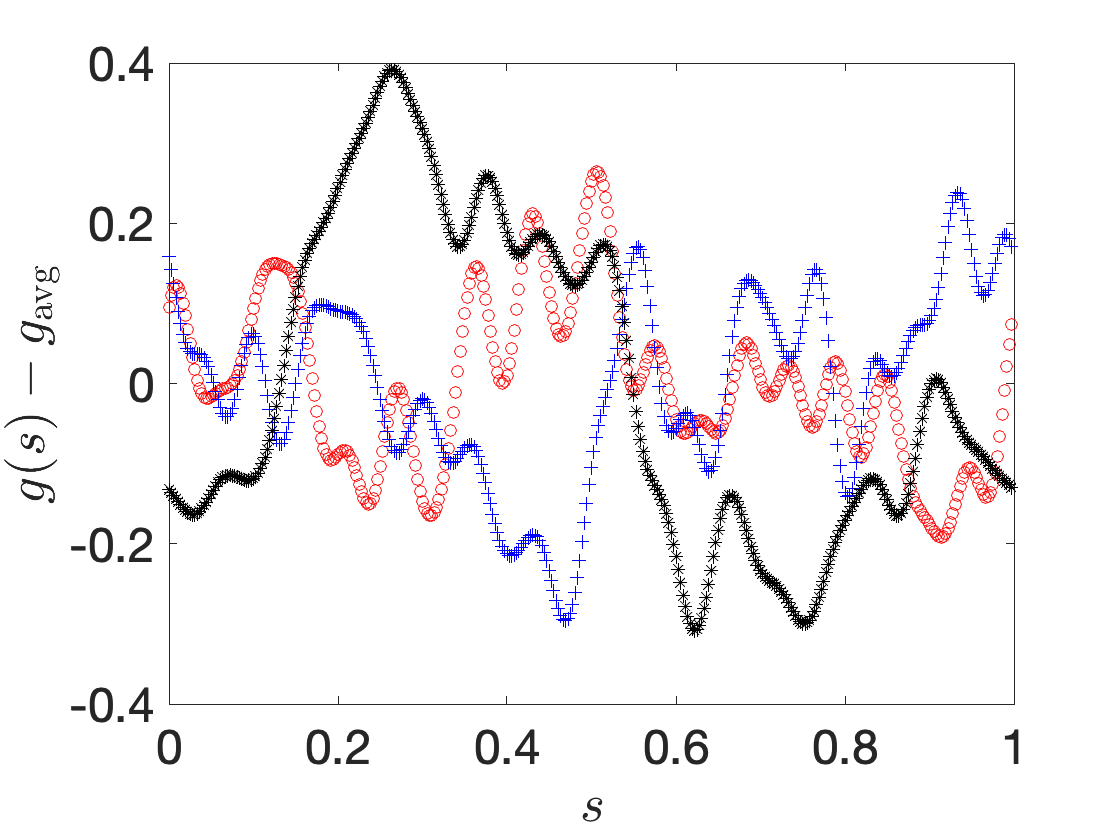}
}
~
\subfloat[$\alpha=2$]{
\includegraphics[width=0.3\textwidth]
 {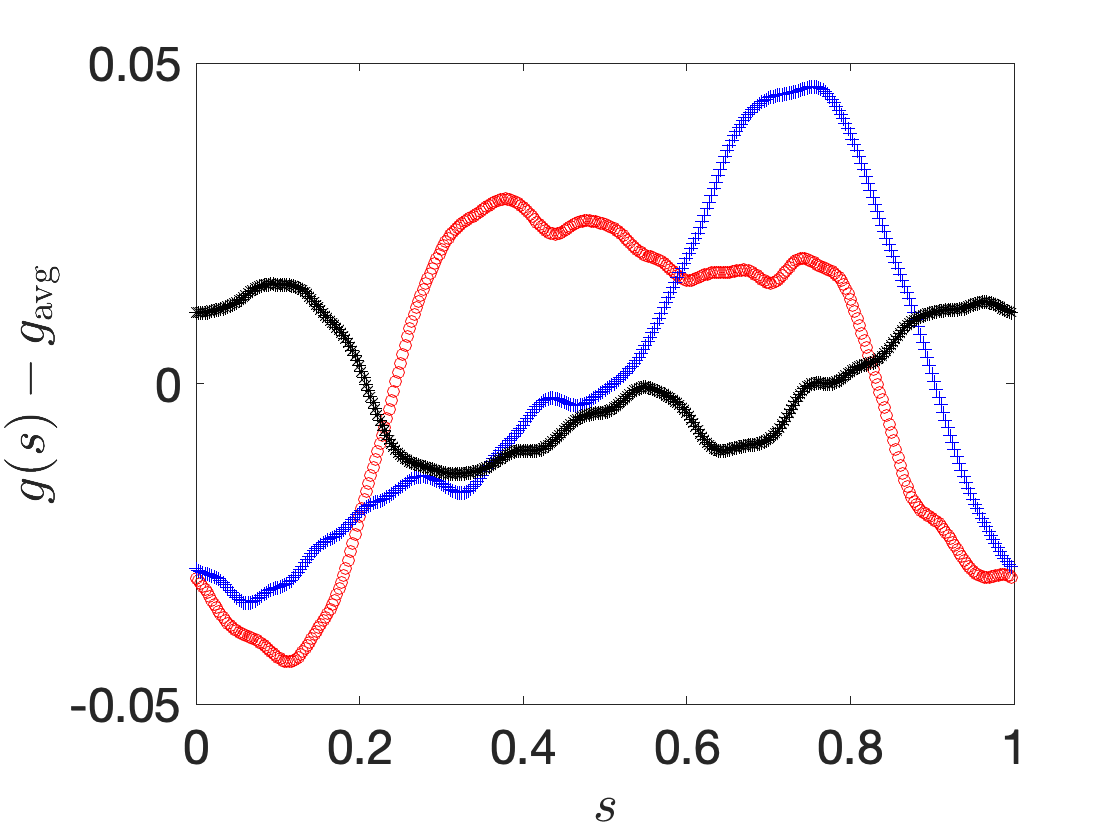}
}
~
\subfloat[$\alpha=3$]{
\includegraphics[width=0.3\textwidth]
 {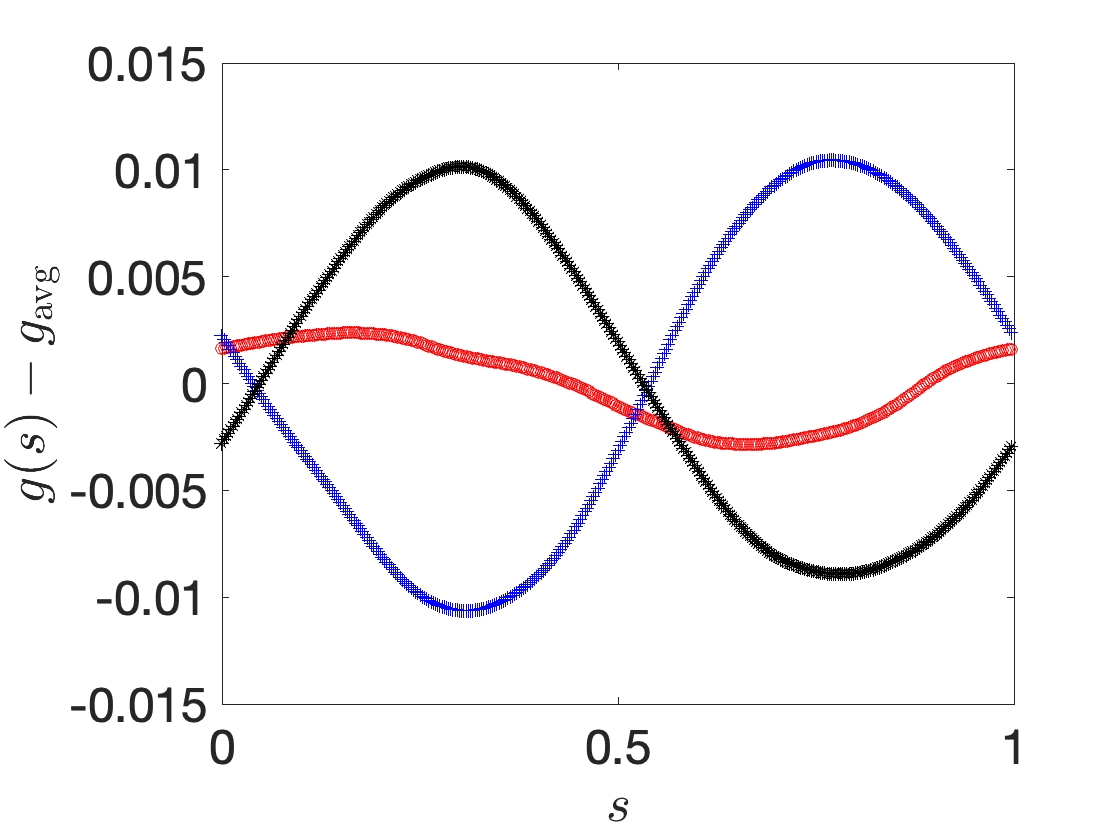}
}
 \caption{linear problem. Samples of random boundary conditions for $3$ choices of $\alpha$.
}
\label{fig:linear_BCs}
\vspace{-15pt}
\end{figure}

We set $p_{\mu} = {\rm Uniform}(\mathcal{P})$ and we consider samples of the random field $g = {\rm Real} [\widetilde{g}(\cdot; \mathbf{c}^{\rm re}, \mathbf{c}^{\rm im}) ]$  (cf. \eqref{eq:randomBC_basic}) with $N_{\rm f}  = 20$. In   \cref{fig:linear_BCs}, we show random samples of $g - g_{\rm avg}$ with $g_{\rm avg} = \int_0^1 g(s) \, ds$ for various choices of $\alpha$: we observe that, as $\alpha$ increases, the samples become increasingly smooth. Given the restriction of the finite element  Lagrangian basis  to the input boundary $\{ \phi_i^{\rm fe} \}_{i\in \texttt{I}_{\rm dir}}$, we further define the random field
\begin{equation}
\label{eq:gaussian_field}
g(x; \mathbf{c}) : = 
\sum_{i\in \texttt{I}_{\rm dir}} \, c_i \phi_i^{\rm fe}(x),
\quad
{\rm with} \;\;
c_i \overset{\rm iid}{\sim} \mathcal{N}(0,1),
\end{equation}
which is used below for comparison. To assess performance, we compare the maximum relative projection error
\begin{equation}
\label{eq:max_relative_error}
E_{\rm max,rel}(\mathcal{Z})
: =
\max_{j=1,\ldots,n_{\rm test}} \;
\frac{\| \Pi_{\mathcal{Z}^{\perp}} u_{\mu^{(j)}, g^{(j)}} |_{\widehat{\Omega}}  \|_{  H^1(\widehat{\Omega}) }    }{
\|  u_{\mu^{(j)}, g^{(j)}} |_{\widehat{\Omega}}  \|_{  H^1(\widehat{\Omega}) }
},
\quad
\mu^{(j)}  \overset{\rm iid}{\sim}  {\rm Uniform} (\mathcal{P}), \; \;
g^{(j)}  \overset{\rm iid}{\sim}  p_{\rm bc},
\end{equation}
for the two choices of $p_{\rm bc}$ ---   ``smooth'' (with $\alpha=1$) and ``Gaussian'' \eqref{eq:gaussian_field} --- and $n_{\rm test}=100$.

\begin{figure}[t]
\centering

\subfloat[][smooth training

smooth test
]{
\includegraphics[width=0.4\textwidth]
{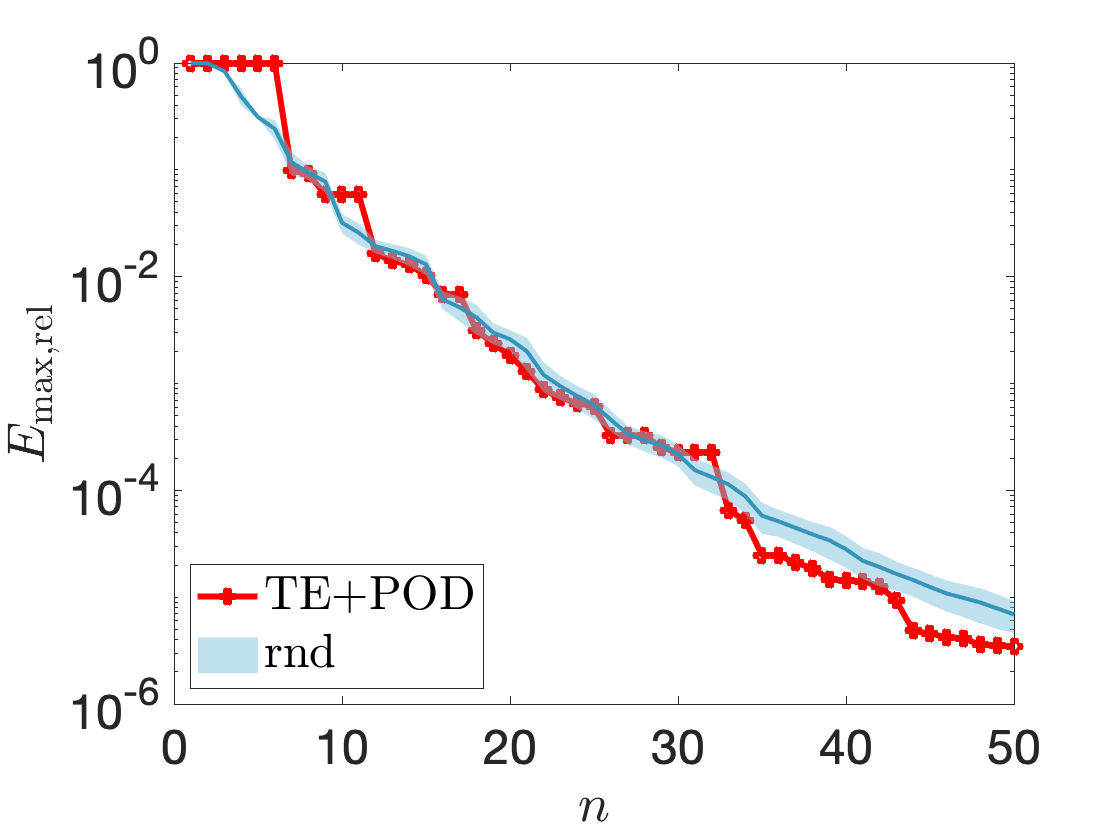}
}
~
\subfloat[][Gaussian training

smooth test]{
\includegraphics[width=0.4\textwidth]
 {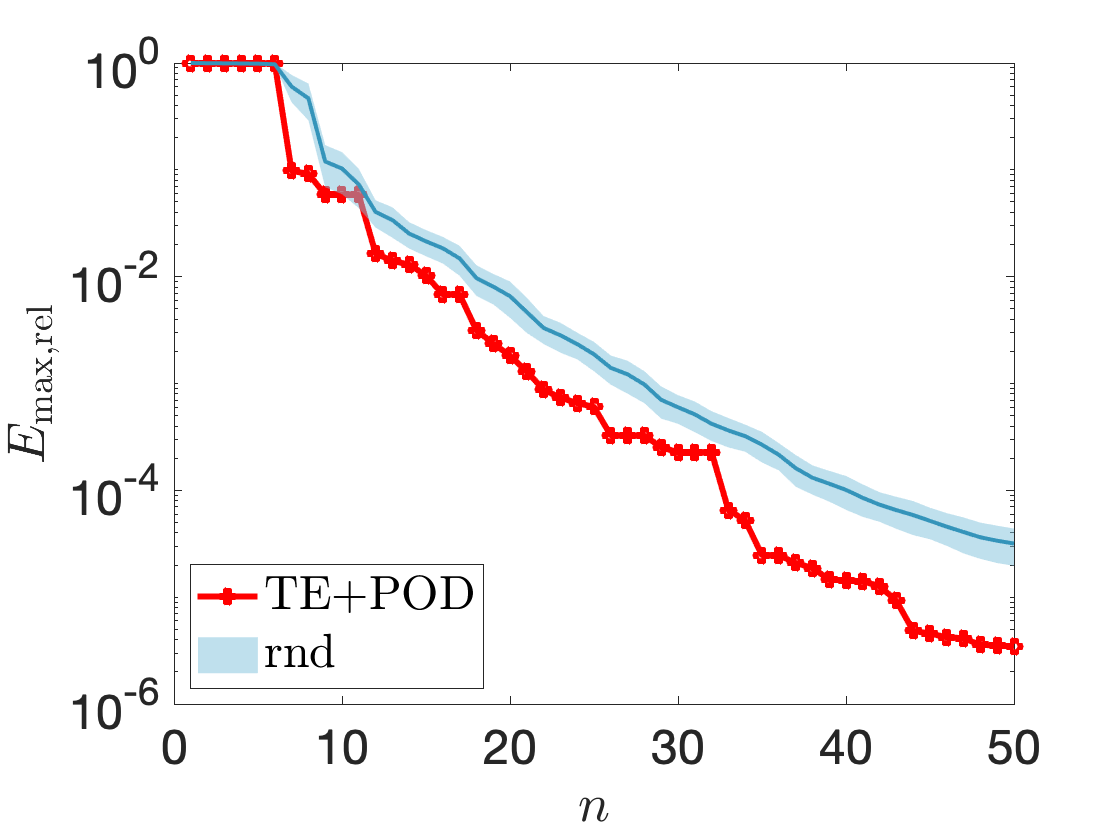}
}

\subfloat[][smooth training

Gaussian test
]{
\includegraphics[width=0.4\textwidth]
{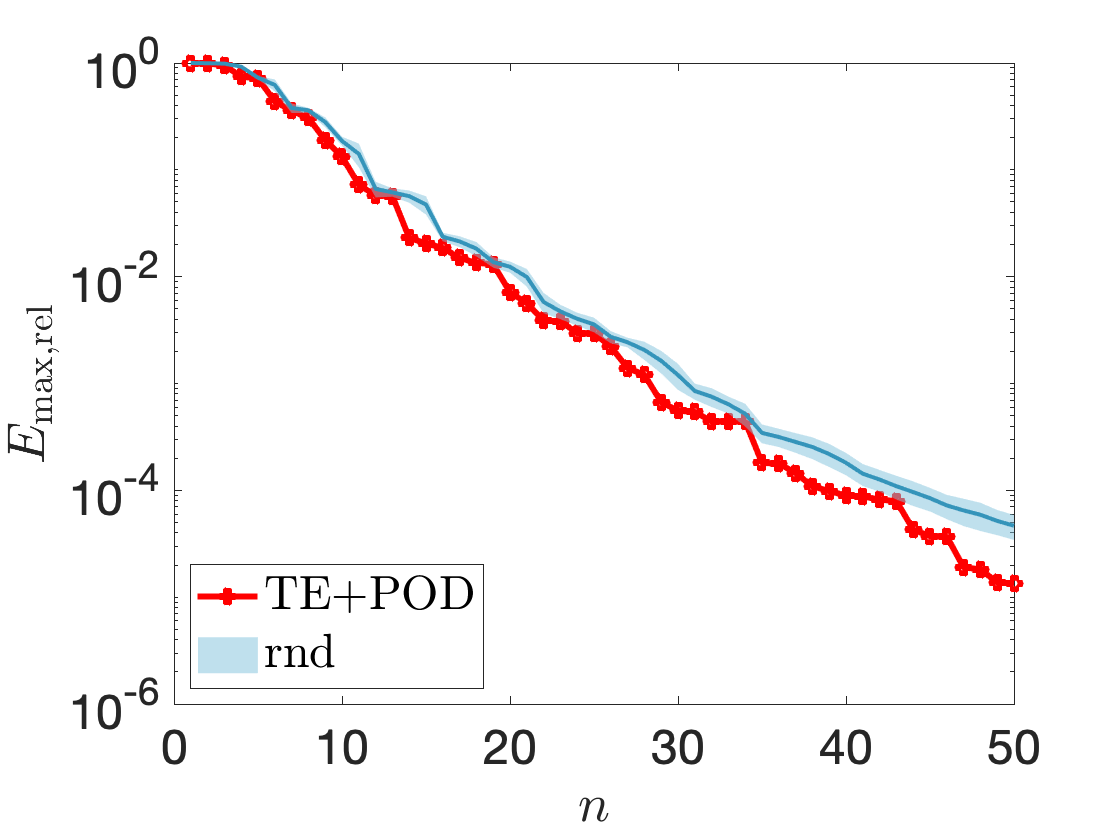}
}
~
\subfloat[][Gaussian training

Gaussian test
]{
\includegraphics[width=0.4\textwidth]
 {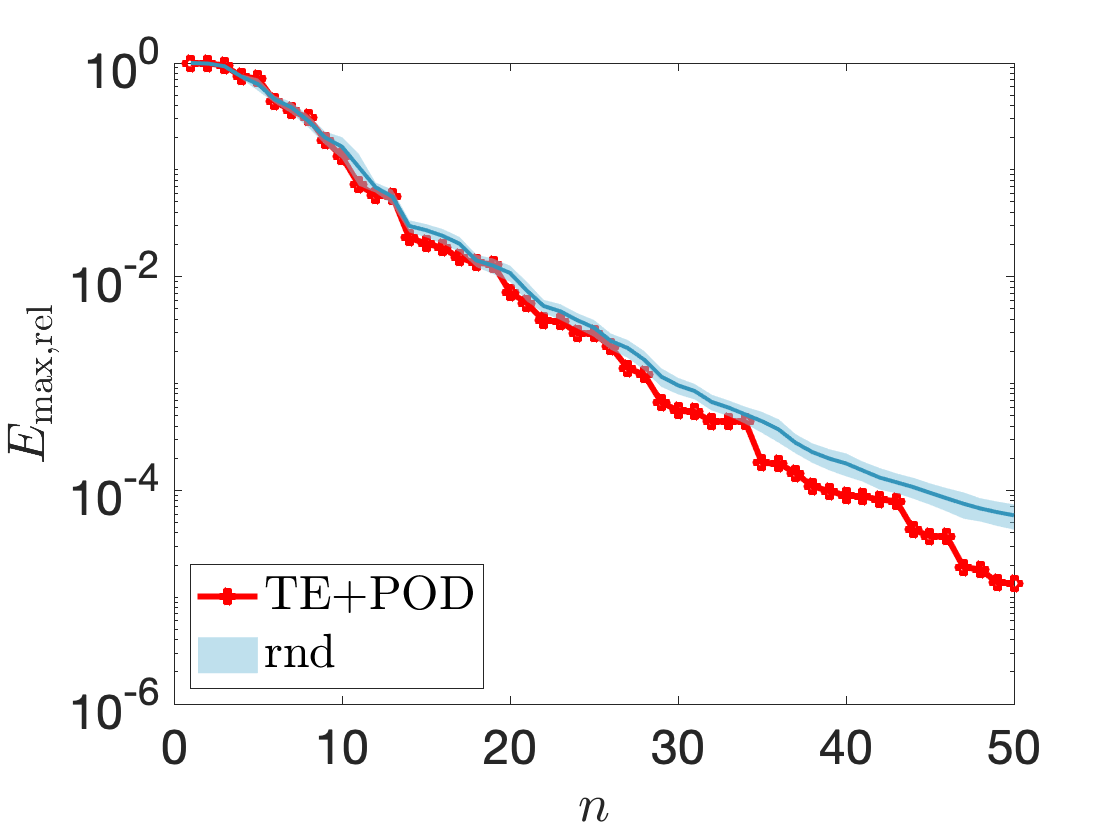}
}
 \caption{linear problem. Out-of-sample performance; comparison with deterministic training for $N_{\rm r}=100$ choices of the random samples and for two fixed test sets.
(a) - (b) smooth test set.
(c) - (d) Gaussian test set. 
}
\vspace{-15pt}
\label{fig:linear_BCs_performance}
\end{figure}

 \Cref{fig:linear_BCs_performance} shows the results for smooth and Gaussian training and test sets.
Here, we consider training sets of size $n_{\rm train}=50$ in   \cref{alg:randomBC};
furthermore, we compare errorbar plots based on $100$ independent choices of the training set.
We observe that our smooth sampling strategy is nearly as effective as TE+POD for $n \lesssim 40$ for both smooth and Gaussian test sets. This result empirically demonstrates that randomized methods are extremely effective to identify dominant POD modes even for nonlinear transfer  operators. 
We further observe that Gaussian sampling is clearly inferior when tested on smooth data, while it performs as accurately as smooth sampling on the Gaussian test set: we conjecture that this behavior is due to the low-pass filtering properties of the differential operator.

\begin{figure}[t]
\centering

\includegraphics[width=0.4\textwidth]
{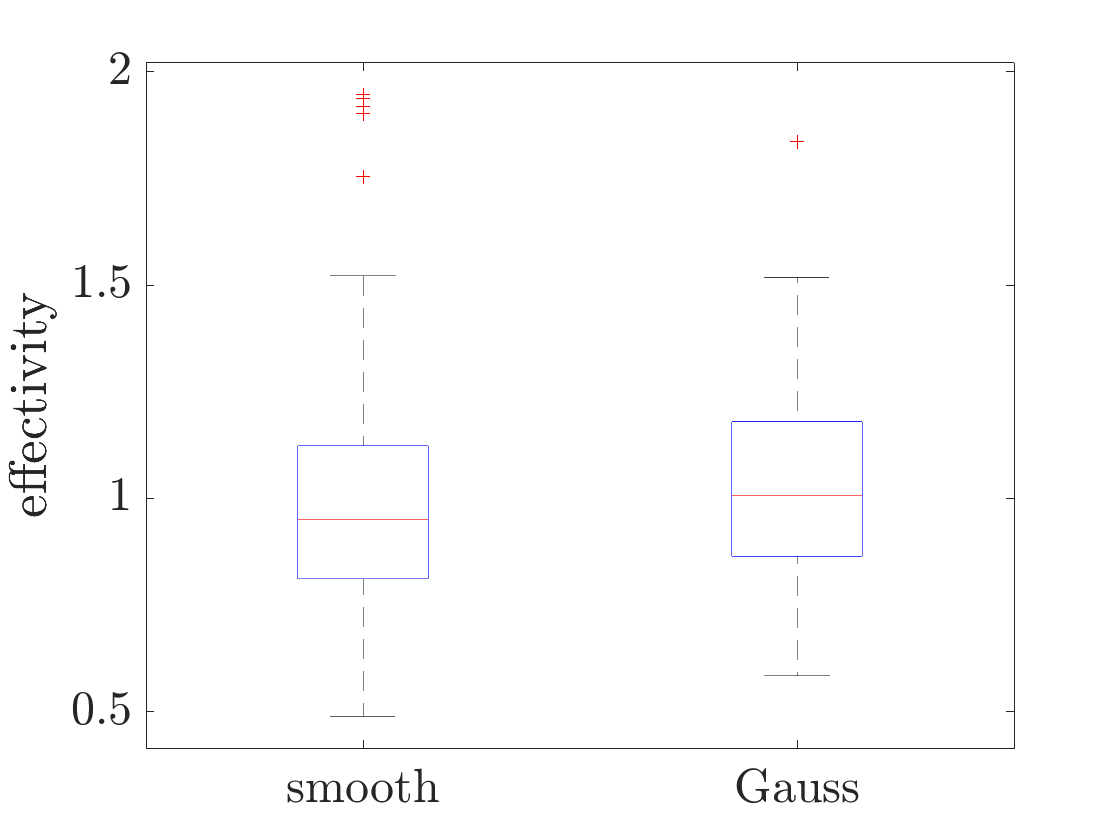}
 
 \caption{linear problem. Effectivity of the error indicator
$\eta = \frac{\widehat{E}( n_{\rm test} = 10 )}{\widehat{E}( n_{\rm test} = 100 )}$
for $100$ independent runs and for both  Gaussian and smooth training. 
}
\label{fig:linear_BCs_effectivity}
\vspace{-15pt}
\end{figure}

In   \cref{fig:linear_BCs_effectivity}, we show the behavior of the error indicator $\widehat{E}$ in   \cref{remark:error_probabilistic}: more precisely, we show boxplots of the approximate effectivity
{$
\eta = \frac{\widehat{E}( n_{\rm test} = 10 )}{\widehat{E}( n_{\rm test} = 100 )}
$}
for $100$ independent runs and for both  Gaussian and smooth training. Note that, with very high-probability,  $\eta \in [0.5, 1.5]$. Note also, however, that 
$\widehat{E}$ strongly depends on the choice of the sampling distribution, which in practice is largely unknown.

\subsection{Application to  the  nonlinear diffusion problem}
\label{sec:nonlinear_diff}

We consider the application to the nonlinear diffusion problem introduced in  \cref{sec:nonlinear_diff}. We apply   \cref{alg:randomBC} with $n_{\rm train}=200$; we set $\{ p_{\mu}^{\bullet}\}_{\bullet}$ as discussed in   \cref{sec:random_training} 
and we consider the smooth sampler described in   \cref{alg:randomBC_all} for 
$N_{\rm f}=20$, 
various choices of $\alpha$ and $\bar{u}_{\rm max}$ ---  we recall that
    \cref{alg:randomBC_all} is the generalization of   \cref{alg:randomBC} to the case of three components. We further compare performance with randomized training based on the random field
\begin{equation}
\label{eq:gaussBC_nonlinear}
\begin{array}{l}
\displaystyle{
g^{\bullet}(x; \mathbf{c}) : = \sum_{i \in \texttt{I}_{\rm dir}^{\bullet} } \, \mathfrak{f} (c_i, \bar{u}_{\rm max}) \phi_i^{\rm fe, \bullet}(x),
\quad
c_i \overset{\rm iid}{\sim} \mathcal{N} 
\left(
\frac{\bar{u}_{\rm max}}{2}, \frac{\bar{u}_{\rm max}^2}{4}
\right),
}
\\[3mm]
\displaystyle{
\mathfrak{f} (c,u) = \max\{ \min \{c, u  \}, 0 \},
}
\end{array}
\end{equation}
where $\{ \texttt{I}_{\rm dir}^{\bullet} \}_{\bullet}$ denote the set of indices of the mesh on the  patch input boundaries and $\{  \phi_i^{\rm fe, \bullet} \}_{\bullet}$ are the Lagrangian bases associated with the high-fidelity discretization.
We refer to the sampling procedure in
  \cref{alg:randomBC_all}  as \emph{smooth sampling}; 
we refer to the sampling procedure based on \eqref{eq:gaussBC_nonlinear} as \emph{Gaussian sampling}.

We compute $n_{\rm test}=30$ global solutions for $n_{\rm dd}=10$ ($N_{\rm dd}=100$) components; then, we define the test datasets
$\{ \mathcal{D}^{\bullet} \}_{\bullet \in \{ \texttt{co}, \texttt{ed}, \texttt{int}\}}$ by extracting the solution in each element of $\mathbb{V}$ --- ${\rm card} (\mathcal{D}^{\bullet}) = 1920$ (resp., $120,960$) for the internal  (resp., corner, edge) component. Finally, we introduce the localized error indicators
\begin{equation}
\label{eq:error_indicator_numerical}
E_{\rm avg,rel}^{\bullet} (  \mathcal{Z}^{\bullet}  )
= 
\frac{1}{{\rm card} (   \mathcal{D}^{\bullet}    )}
\sum_{w \in   \mathcal{D}^{\bullet} }
\;
\frac{ \| w - \Pi_{\mathcal{Z}^{\bullet}} w   \|_{\bullet}}{   \|  w  \|_{\bullet}  },
\;\;
\bullet \in \{ \texttt{co}, \texttt{ed}, \texttt{int} \},
\end{equation}
which are used to assess performance.

\begin{figure}[t]
\centering

\subfloat[$\alpha=1$]{
\includegraphics[width=0.3\textwidth]
{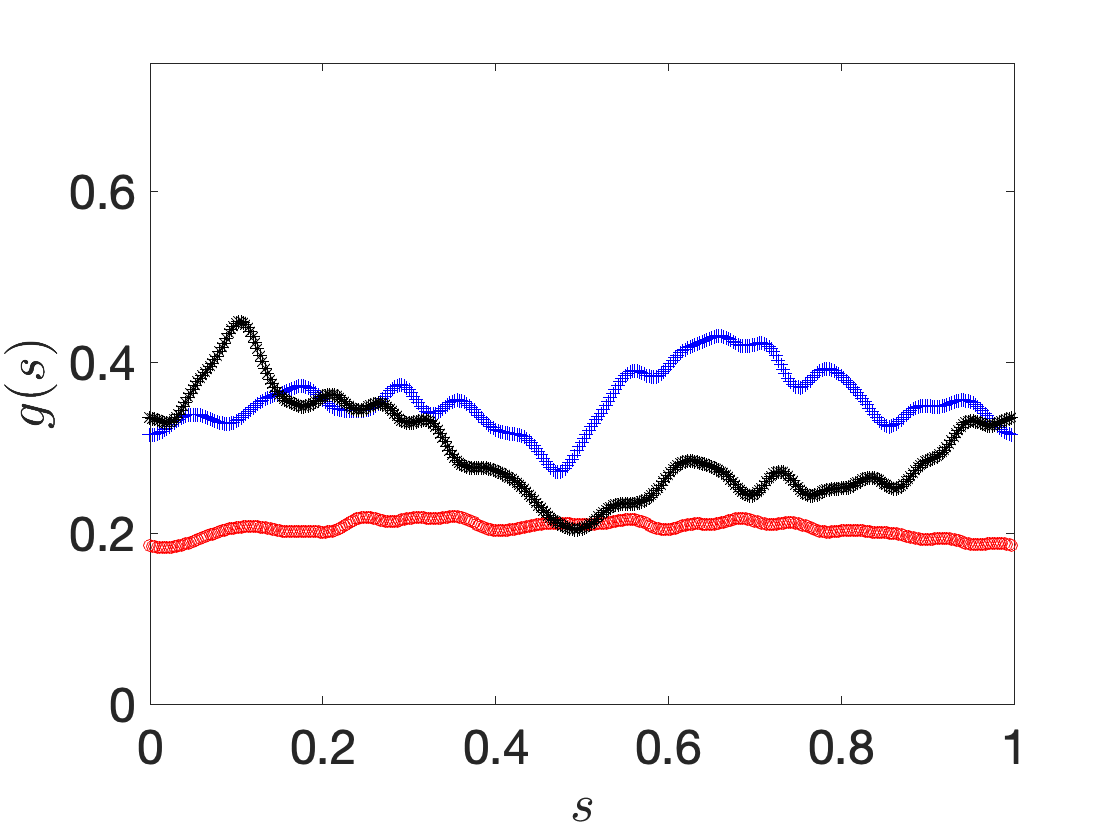}
}
~
\subfloat[$\alpha=2$]{
\includegraphics[width=0.3\textwidth]
{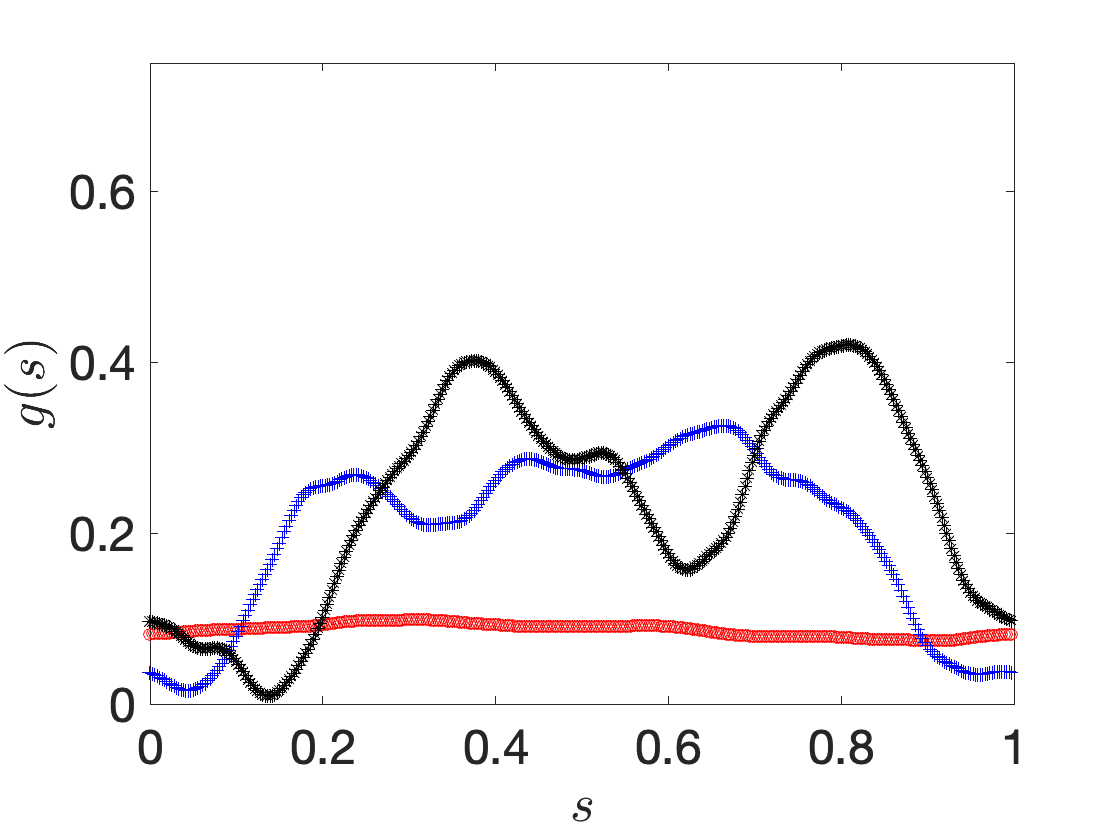}
}
~
\subfloat[$\alpha=3$]{
\includegraphics[width=0.3\textwidth]
{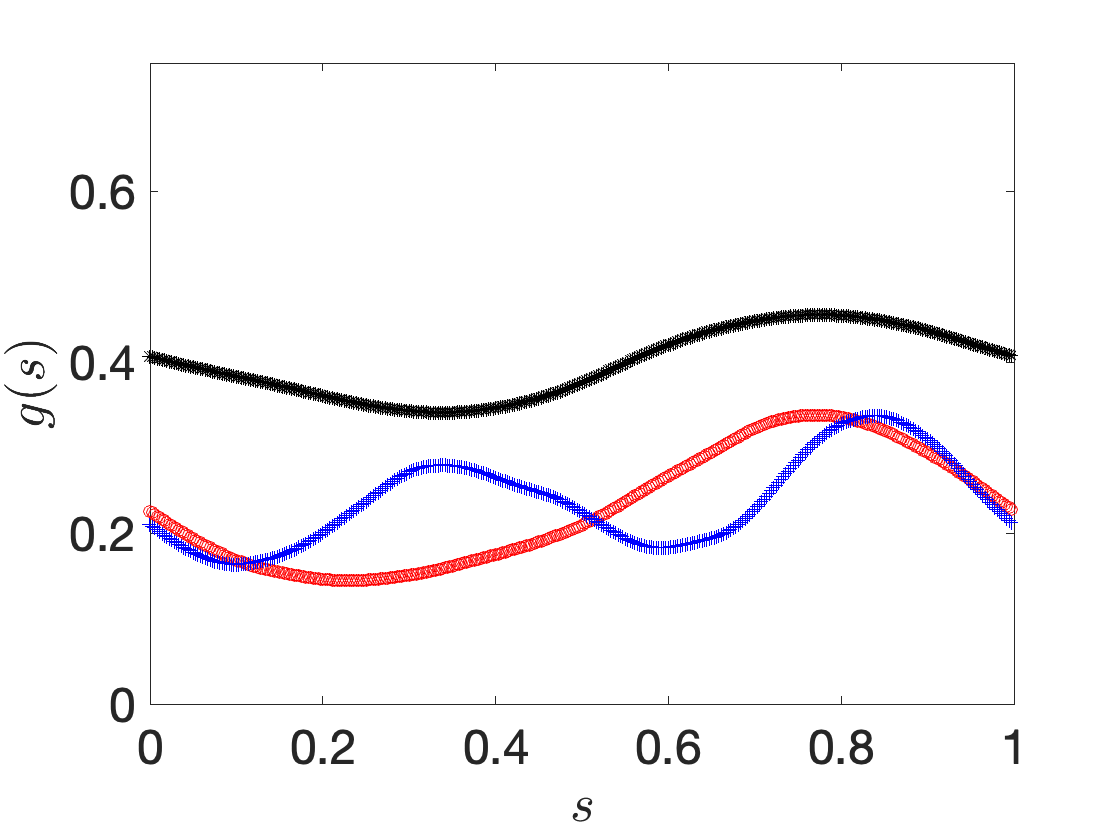}
}

\subfloat[$\alpha=1$]{
\includegraphics[width=0.3\textwidth]
{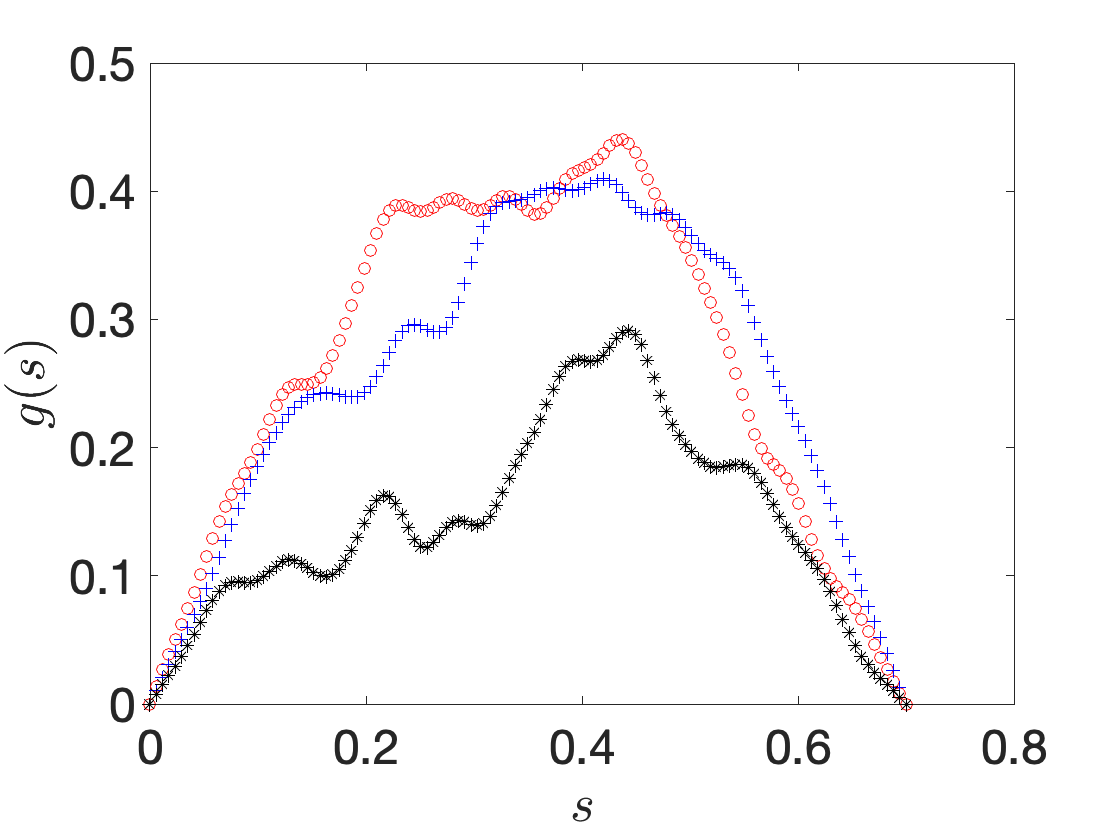}
}
~
\subfloat[$\alpha=2$]{
\includegraphics[width=0.3\textwidth]
{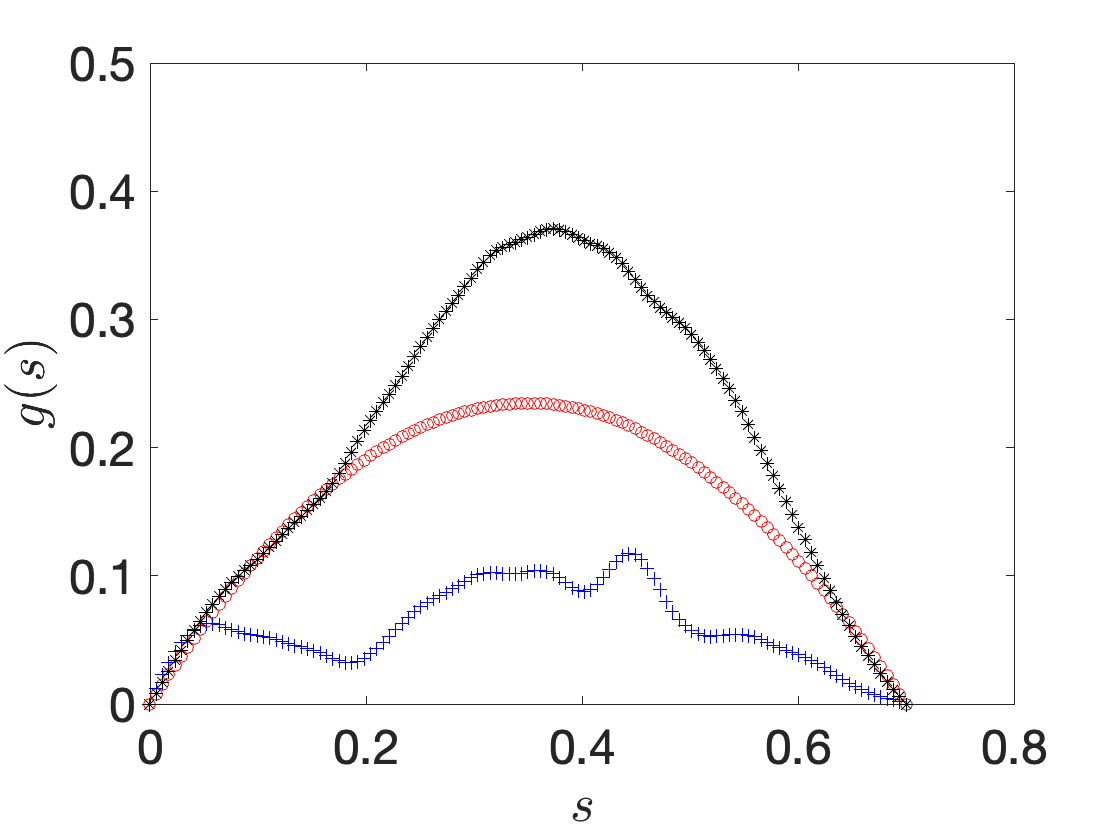}
}
~
\subfloat[$\alpha=3$]{
\includegraphics[width=0.3\textwidth]
{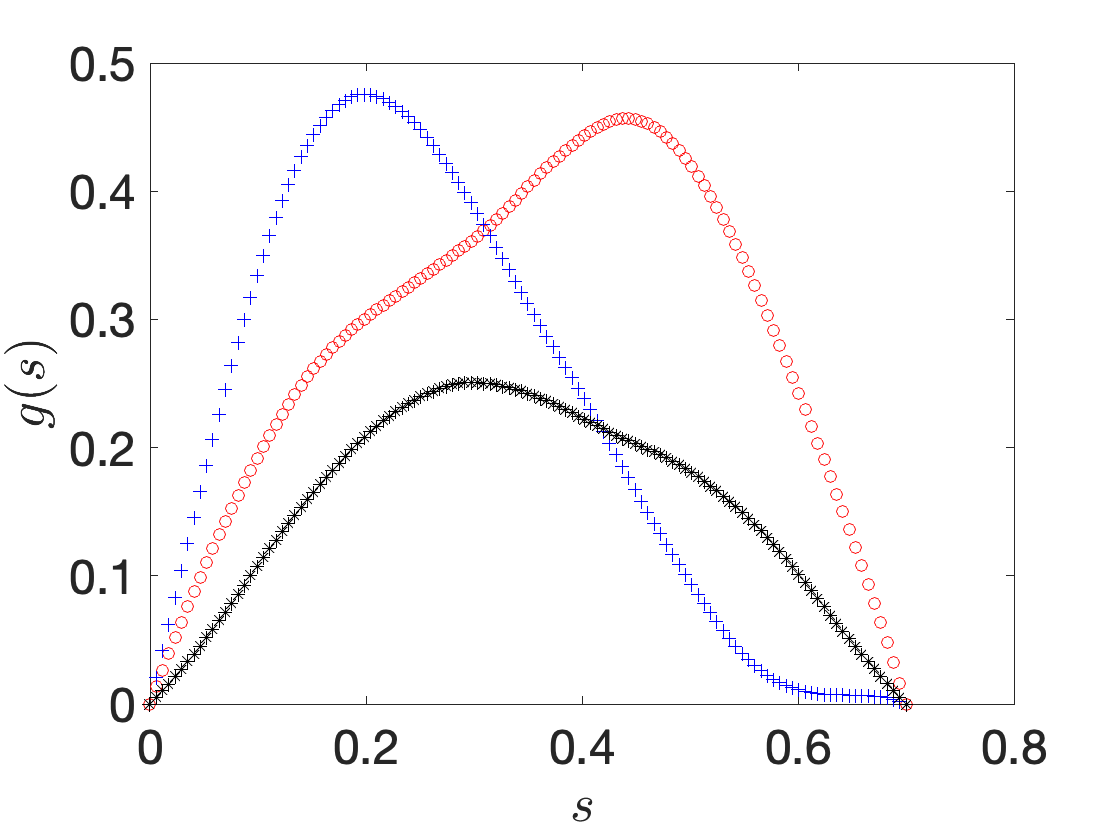}
}

 \caption{nonlinear problem. Samples of random boundary conditions for three choices of $\alpha$  ($N_{\rm f}=20$, $\bar{u}_{\rm max}=0.5$). (a)-(b)-(c) internal component.
 (d)-(e)-(f) edge component.
}
\label{fig:nonlinear_BCs}
\vspace{-15pt}
\end{figure}

 \Cref{fig:nonlinear_BCs} shows random samples of the boundary conditions on $\Gamma_{\rm in}$ for internal and edge components as provided by   \cref{alg:randomBC_all}  for various choices of $\alpha$ and $N_{\rm f}=20$ and $\bar{u}_{\rm max}=0.5$. As for the linear case, the value of $\alpha$ encodes the spatial smoothness of the samples. We further observe that  \cref{alg:randomBC_all}  automatically enforces the proper condition at the extrema $s=0$ and $s=1$ ---
$g(0)=g(1) =0$ for $\bullet \in \{ \texttt{co}, \texttt{ed}  \}$,
$g(0)=g(1)$ for $\bullet =\texttt{int} $.
 %visualization of nonlinear BCs   

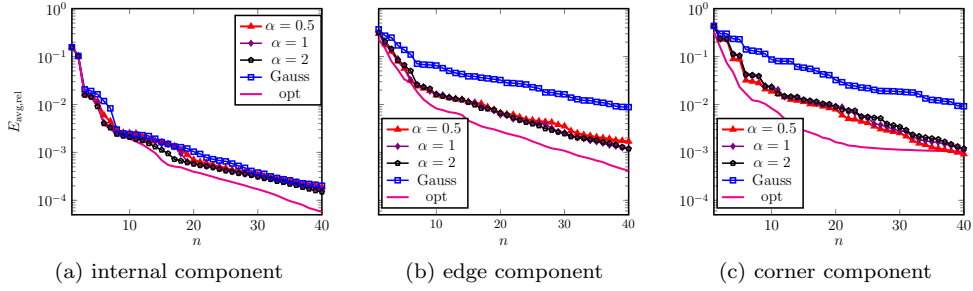
\begin{figure}[t]
\centering

\subfloat[internal component]{
\begin{tikzpicture}[scale=0.48]
\begin{semilogyaxis}[
xlabel = {\large{$n$}},
  ylabel = {\large {$E_{\rm avg,rel}$}},
legend entries = {$\alpha=0.5$, $\alpha=1$, $\alpha=2$,Gauss,opt},
  line width=1.2pt,
  mark size=2.0pt,
xmin=1,   xmax=40,
 ymin=0.00005,   ymax=1,
xlabel style = {font=\large,xshift=-0ex},
ylabel style = {font=\large,yshift=2ex},
 yticklabel style = {font=\large,xshift=0ex},
xticklabel style = {font=\large,yshift=0ex},
legend style={at={(0.65,0.75)},anchor=west,font=\large}
  ]

 \addplot[ultra thick,color=red,mark=triangle*]  table {data/loc/alpha05_int.dat};

\addplot[color=violet,mark=diamond*]  table {data/loc/alpha1_int.dat};
   
\addplot[color=black,mark=pentagon] table {data/loc/alpha2_int.dat};
  
\addplot[color=blue,mark=square] table {data/loc/gaussian_int.dat};  
 
\addplot[ultra thick,color=magenta,mark=none] table {data/loc/opt_int.dat};   
    
\end{semilogyaxis}
\end{tikzpicture} 
}
~
\subfloat[edge component]{
\begin{tikzpicture}[scale=0.48]
\begin{semilogyaxis}[
xlabel = {\large{$n$}},
legend entries = {$\alpha=0.5$, $\alpha=1$, $\alpha=2$,Gauss,opt},
  line width=1.2pt,
  mark size=2.0pt,
xmin=1,   xmax=40,
 ymin=0.00005,   ymax=1,
xlabel style = {font=\large,xshift=-0ex},
ylabel style = {font=\large,yshift=5ex},
 yticklabels={,,}
xticklabel style = {font=\large,yshift=0ex},
legend style={at={(0.01,0.25)},anchor=west,font=\large}
  ]

 \addplot[ultra thick,color=red,mark=triangle*]  table {data/loc/alpha05_edge.dat};
%\addlegendentry{model}
  
    \addplot[color=violet,mark=diamond*]  table {data/loc/alpha1_edge.dat};
   
  \addplot[color=black,mark=pentagon] table {data/loc/alpha2_edge.dat};
  
  \addplot[color=blue,mark=square] table {data/loc/gaussian_edge.dat};  
 
  \addplot[ultra thick,color=magenta,mark=none] table {data/loc/opt_edge.dat};   
    
\end{semilogyaxis}
\end{tikzpicture} 
}
~
\subfloat[corner component]{
\begin{tikzpicture}[scale=0.48]
\begin{semilogyaxis}[
xlabel = {\large{$n$}},
legend entries = {$\alpha=0.5$, $\alpha=1$, $\alpha=2$,Gauss,opt},
  line width=1.2pt,
  mark size=2.0pt,
xmin=1,   xmax=40,
 ymin=0.00005,   ymax=1,
xlabel style = {font=\large,xshift=-0ex},
ylabel style = {font=\large,yshift=2ex},
 yticklabel style = {font=\large,xshift=0ex},
xticklabel style = {font=\large,yshift=0ex},
legend style={at={(0.01,0.25)},anchor=west,font=\large}
  ]

 \addplot[ultra thick,color=red,mark=triangle*]  table {data/loc/alpha05_corner.dat};
 % \addlegendentry{model}
  
    \addplot[color=violet,mark=diamond*]  table {data/loc/alpha1_corner.dat};
   
  \addplot[color=black,mark=pentagon] table {data/loc/alpha2_corner.dat};
  
  \addplot[color=blue,mark=square] table {data/loc/gaussian_corner.dat};  
 
  \addplot[ultra thick,color=magenta,mark=none] table {data/loc/opt_corner.dat};   
    
\end{semilogyaxis}
\end{tikzpicture} 
}
 \caption{nonlinear problem. Local approximation errors \eqref{eq:error_indicator_numerical}  for three choices of $\alpha$  ($N_{\rm f}=20$, $\bar{u}_{\rm max}=0.5$), and for Gaussian sampling \eqref{eq:gaussBC_nonlinear}.
}
\label{fig:nonlinear_local_projection_error}
\vspace{-15pt}
\end{figure}

 \Cref{fig:nonlinear_local_projection_error} shows the behavior of the relative errors \eqref{eq:error_indicator_numerical} for the three components,
for smooth sampling for three choices of $\alpha$  ($N_{\rm f}=20$, $\bar{u}_{\rm max}=0.5$), and for Gaussian sampling \eqref{eq:gaussBC_nonlinear}. To provide a reference, we show also performance of the POD spaces based on the datasets
$\{ \mathcal{D}_{\rm test}^{\bullet} \}_{\bullet \in \{ \texttt{co}, \texttt{ed}, \texttt{int}}$ (``opt'') generated using  $30$ additional global simulations with $N_{\rm dd}=100$ components. We observe that smooth sampling outperforms Gaussian sampling for the boundary components: we believe that this is due to the presence of strong Dirichlet conditions on $\partial \widehat{U}^{\bullet} \setminus \widehat{\Gamma}_{\rm in}^{\bullet}$. We further observe that results weakly depend on the choice of $\alpha$.

%local approximation errors
   
\begin{figure}[t]
\centering
\subfloat[]{
\begin{tikzpicture}[scale=0.62]
\begin{semilogyaxis}[
xlabel = {\large{$n$}},
  ylabel = {\large {$E_{\rm avg,rel}$}},
legend entries = {proj ($H^1$), Galerkin ($H^1$), Galerkin ($L^2$)},
  line width=1.2pt,
  mark size=3.0pt,
xmin=1,   xmax=40,
 ymin=0.0001,   ymax=1,
xlabel style = {font=\large,xshift=-0ex},
ylabel style = {font=\large,yshift=5ex},
 yticklabel style = {font=\large,xshift=0ex},
xticklabel style = {font=\large,yshift=0ex},
legend style={at={(0.36,0.81)},anchor=west,font=\Large}
  ]

 \addplot[ultra thick,color=red,mark=triangle*]  table {data/glo/err_proj.dat};
%\addlegendentry{model}
  
    \addplot[color=violet,mark=diamond*]  table {data/glo/err_gal.dat};
   
  \addplot[color=black,mark=pentagon] table {data/glo/err_galL2.dat};
   
\end{semilogyaxis}
\end{tikzpicture}    
   }
 ~~~
 \subfloat[]{
\begin{tikzpicture}[scale=0.62]
\begin{semilogyaxis}[
xlabel = {\large{$n$}},
  ylabel = {\large {$E_{\rm avg,rel}$  ($H^1$ proj)}},
legend entries = {$\bar{u}_{\rm max}=0.5$ $\alpha=1$, $\bar{u}_{\rm max}=0.75$ $\alpha=4$},
  line width=1.2pt,
  mark size=3.0pt,
xmin=1,   xmax=40,
 ymin=0.0001,   ymax=1,
xlabel style = {font=\large,xshift=-0ex},
ylabel style = {font=\large,yshift=5ex},
 yticklabel style = {font=\large,xshift=0ex},
xticklabel style = {font=\large,yshift=0ex},
legend style={at={(0.27,0.81)},anchor=west,font=\Large}
  ]

 \addplot[ultra thick,color=red,mark=triangle*]  table {data/glo/err_proj.dat};

\addplot[color=violet,mark=diamond*]  table {data/glo/err_proj_075.dat};
   
\end{semilogyaxis}
\end{tikzpicture}    
   }  
    
 \caption{nonlinear problem. Performance of PUM CB-ROM on $n_{\rm test}=30$ global solutions for $N_{\rm dd}=100$.
  (a) Galerkin error vs projection error.
  (b) projection error for two choices of the parameters in \cref{alg:randomBC_all}. }  
 \label{fig:global_error_noadapt}
 \vspace{-15pt}
\end{figure}
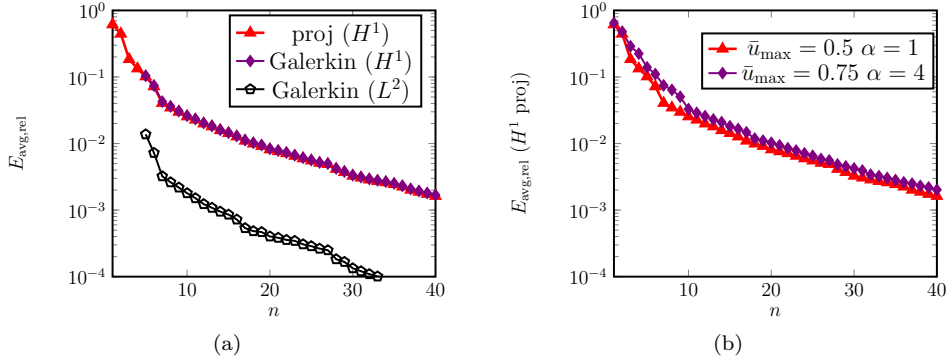

\Cref{fig:global_error_noadapt} shows the performance of the CB-ROM based on PUM.
In \cref{fig:global_error_noadapt}(a), we show the average global $L^2$   and $H^1$ relative errors over the test set of $n_{\rm test}=30$ global simulations and we also compare with the projection error.
We here consider $\bar{u}_{\rm max}=0.5$, $N_{\rm f}=20$, and $\alpha=1$. We observe that Galerkin projection is nearly optimal for all choices of $n$; we further observe exponential convergence with respect to $n$.
In  \cref{fig:global_error_noadapt}(b), we compare the $H^1$ relative projection error for $\bar{u}_{\rm max}=0.5$,  $N_{\rm f}=20$, and $\alpha=1$ with the results obtained for $\bar{u}_{\rm max}=0.75$,  $N_{\rm f}=20$, and $\alpha=4$: we observe that results weakly depend on the choice of these two hyper-parameters.
   
%global errors  

\subsection{Adaptive enrichment}
\label{sec:results_adaptive}

{We apply   \cref{alg:localized_enrichment} with error indicator $\Delta_{\mu}$ \cref{eq:residual_adaptive}} to the model problem of   \cref{sec:nonlinear_diff}. We consider $n_{\rm train}^{\rm loc}=30$, 
$n^{\rm loc}=20$, we set $\{  p_{\mu}^{\bullet}  \}_{\bullet}$ as discussed in  \cref{sec:random_training} and 
we generate random samples of boundary conditions at input ports based on (i)     \cref{alg:randomBC_all} with $N_{\rm f}=20$, $\bar{u}_{\rm max}=0.5$, $\alpha=1$ or
(ii) on iid realizations of \eqref{eq:gaussBC_nonlinear}.
We further consider $n_{\rm train}^{\rm glo}=50$, $n^{\rm glo}=10$,  \texttt{maxit}$=3$, and we generate global configurations using the strategy outlined in   \cref{sec:algorithm_adaptive}.
We assess performance based on $n_{\rm test}=20$ out-of-sample randomly-chosen configurations.
 
\begin{figure}[t]
\centering

\subfloat[]{
\includegraphics[width=0.4\textwidth]
{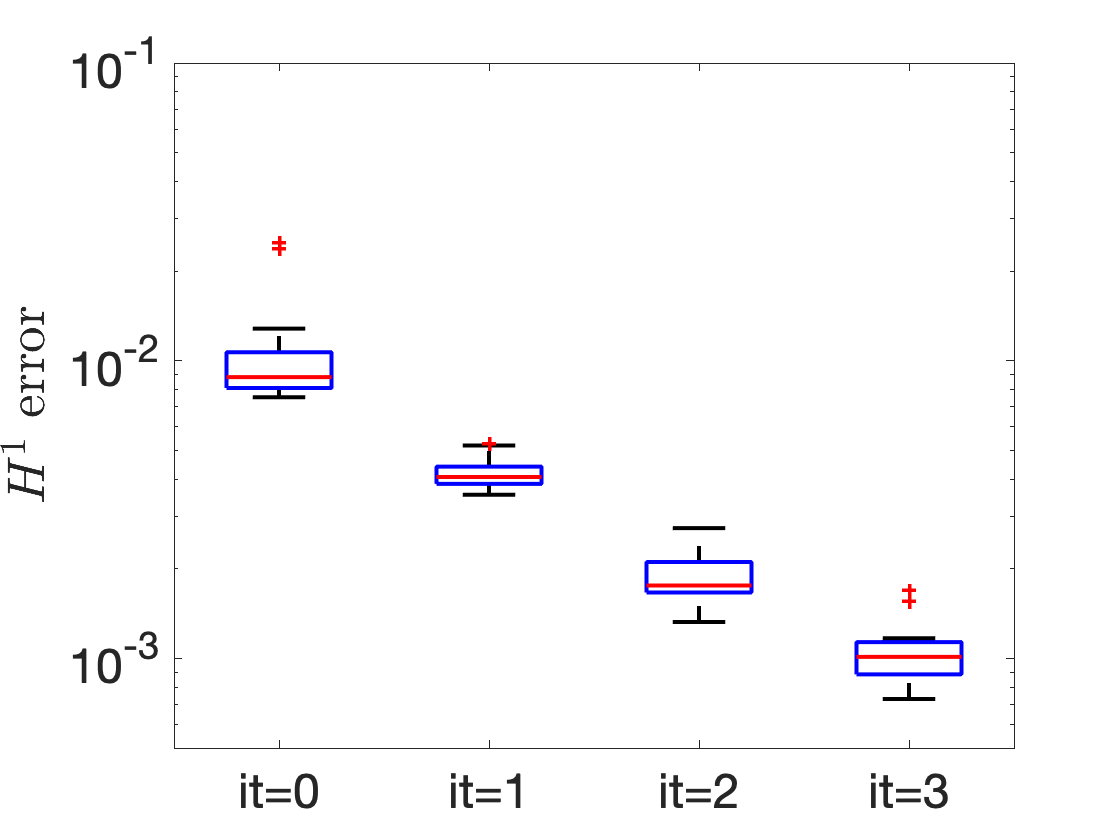}
}
 ~~
 \subfloat[]{
\includegraphics[width=0.4\textwidth]
{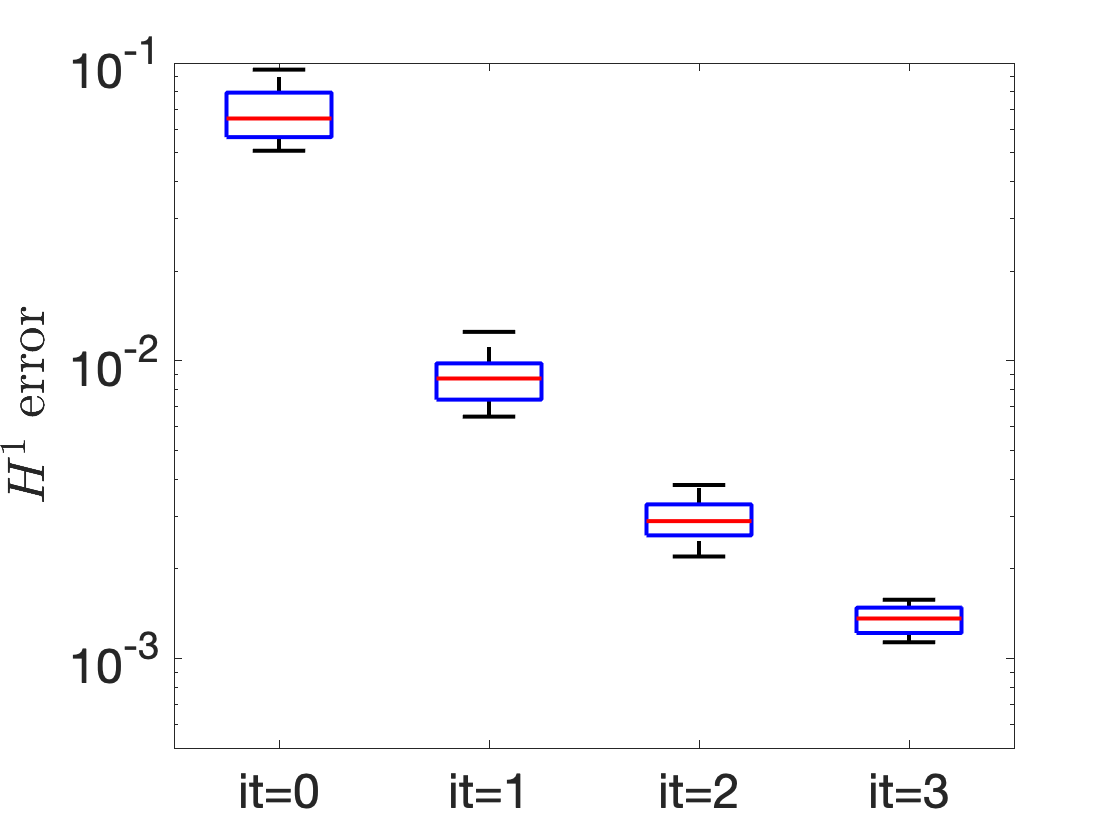}
}

\subfloat[]{
\includegraphics[width=0.4\textwidth]
{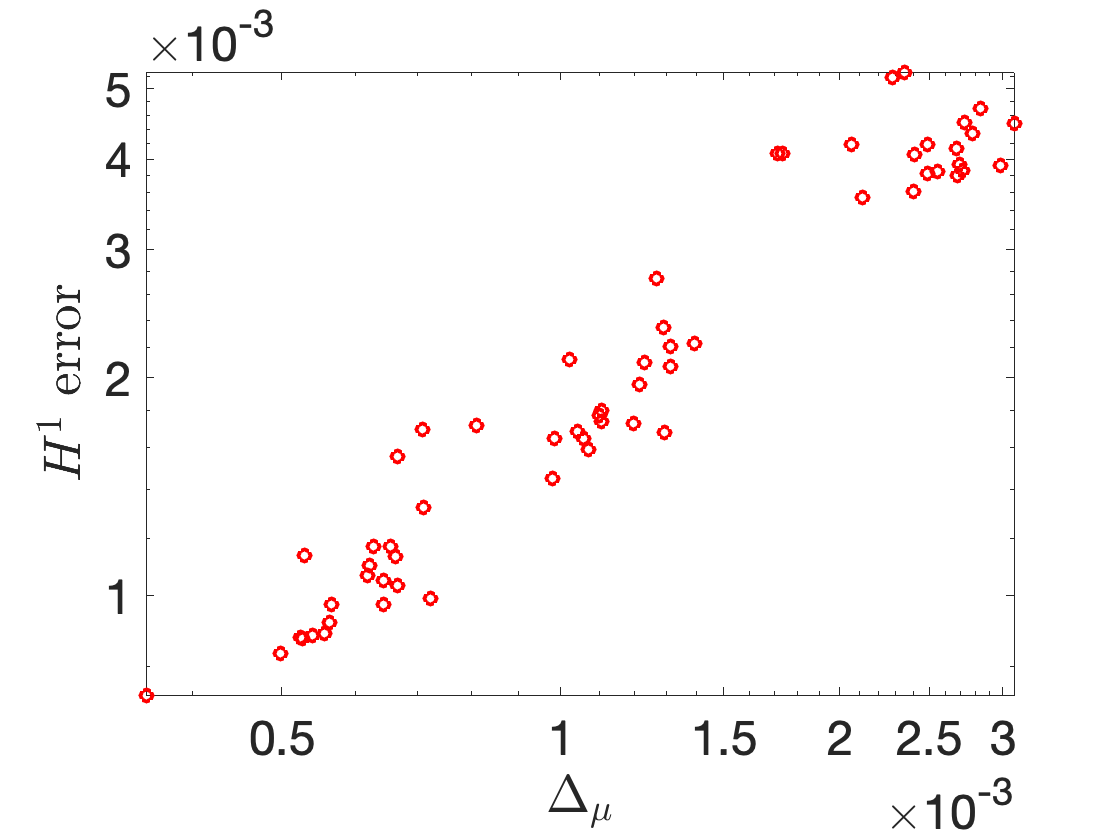}
}
~~
\subfloat[]{
\includegraphics[width=0.4\textwidth]
{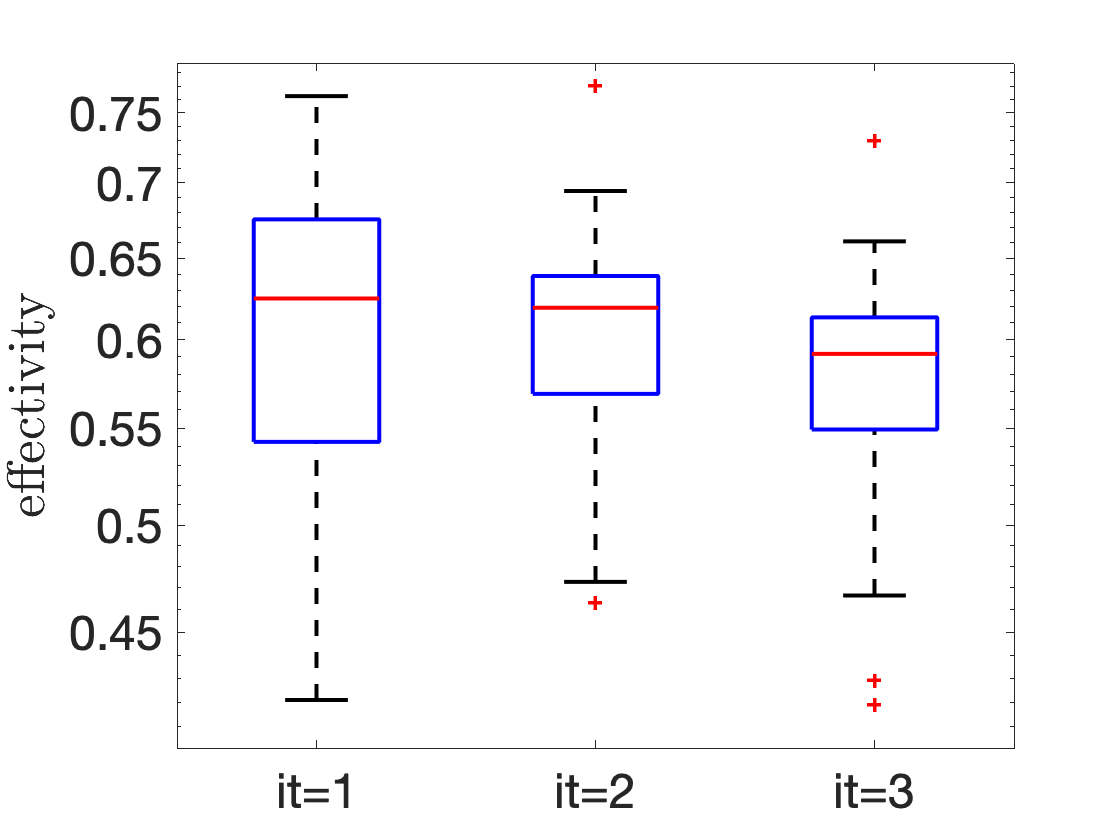}
}

\caption{nonlinear problem; adaptive enrichment. 
(a)-(b) boxplots of the relative $H^1$ error on the test set for smooth and Gaussian  sampling of localized BCs.
(b) correlation between $\Delta_{\mu}$  and relative $H^1$ error (smooth sampling).
(c) out-of-sample effectivity of the error indicator $\Delta_{\mu}/E_{\rm rel}$  \eqref{eq:residual_adaptive}  (smooth sampling).}
\label{fig:adaptive_training}
\vspace{-15pt}
\end{figure}

\Cref{fig:adaptive_training}(a) and (b) show boxplots of the relative $H^1$ error after each iteration of the training algorithm --- iteration $0$ corresponds to the performance of the CB-ROM without global enrichment.
Iteration $0$ corresponds to a reduced space of size $n=20$; iterations $it=1,2,3$ correspond to reduced spaces of size $n=20 + 10 \cdot it$.
 We observe that the enrichment iterations significantly improve performance of the CB-ROM  and reduce  the impact of the initial sampling distribution.
  \Cref{fig:adaptive_training}(c) shows the
correlation between the residual indicator \eqref{eq:residual_adaptive}
and the relative $H^1$ error on the test set for all iterations of the enrichment algorithm for smooth sampling;
  \Cref{fig:adaptive_training}(d)
shows the effectivity of the error indicator
$\eta = \Delta_{\mu}/E_{\rm rel}$ for smooth sampling.
{We observe that the residual indicator is strongly correlated with the global error.}
%: this motivates its use in the termination criterion of   \cref{alg:localized_enrichment}.
 
\section{Conclusions and perspectives}
\label{sec:conclusions}

We presented a CB-pMOR method for parameterized elliptic nonlinear PDEs.
The approach relies on the definition of several archetype components and associated local ROBs and ROMs. CB-pMOR rely on two building blocks:
(i) a localized training strategy for the construction of the local approximation  spaces, and (ii) a DD strategy for online global predictions. In this  paper, we proposed a   localized data compression procedure based on  oversampling and randomized sampling of boundary conditions of controlled smoothness, and we relied on the PUM to devise
global approximation spaces and on
 Galerkin projection to estimate the global state.
 Finally, we proposed an adaptive enrichment procedure that exploits global CB-ROM solves to improve approximation properties of the local reduced spaces.

Numerical results for a nonlinear diffusion problem show the impact of the sampling distribution on performance: given a class of nonlinear PDEs, it is thus  necessary to devise an effective sampler that is informed by the problem of interest.
The approach presented in this work (cf.   \cref{alg:randomBC}) is simple to implement, and incorporates relevant features of the problem of interest --- lower and upper bounds for the solution,  Sobolev regularity, Dirichlet boundary  conditions.
However,  it depends on several hyper-parameters that might be difficult to set \emph{a priori}.
In this respect, we numerically showed that the proposed enrichment strategy reduces the impact of the initial sampling distribution.

  In the future, we wish to  extend the approach in several directions.
  First, we wish to devise specialized hyper-reduction strategies for CB-pMOR methods based on PUM: hyper-reduction is key to reduce efficient  online memory and computational costs. 
  Second, we wish to develop rigorous a posteriori  error estimators for nonlinear PDEs for online certification.
  Third, we wish to analyze performance of randomized algorithms for nonlinear operators: this analysis is key to provide mathematical foundations for  randomized methods for nonlinear problems and also to inform the choice of the sampling distribution.
    
\appendix

\section{Proofs}\label{sec:proofs}

\subsection*{\Cref{th:residual_estimator}}
\begin{proof}
Since the form $\mathcal{G}_{\mu}$ is linear with respect to the second argument, we find that  
\begin{subequations}
\begin{equation}
\label{eq:theorem_a}
\begin{array}{rl}
\displaystyle{\mathcal{G}_{\mu}(u, v ) 
=} &
\displaystyle{
\sum_{i=1}^{N_{\rm  dd}} 
\mathcal{G}_{\mu}(u, v_i ) 
=
\sum_{i=1}^{N_{\rm  dd}} 
\int_{ \omega_i  }   {\eta}_{\mu}^{(i)} (x; u, v_i) \, dx

}
\\[3mm]
= &
\displaystyle{
\sum_{i=1}^{N_{\rm  dd}} 
\left(
\psi_{\mu}[u],  v_i
\right)_{1,\omega_i}
\leq
\sum_{i=1}^{N_{\rm dd}}   \mathcal{r}_{\mu}^{(i)} [u] \, 
\| v_i  \|_{1,\omega_i},}
\\
\end{array}
\end{equation}
where $v_i = v \phi_i$.
Note that in the second identity we used  
\eqref{eq:residual_localized_a}, while in the last bound we applied Cauchy-Schwartz  inequality.
We observe that
$$
\int_{\omega_i} (v \phi_i)^2 \, dx  \leq \int_{\omega_i} v^2 \, dx 
$$
and
$$
\begin{array}{l}
\displaystyle{
\int_{\omega_i} \| \nabla (v \phi_i) \|_2^2 \, dx
=
\int_{\omega_i} \| v \nabla  \phi_i + \phi_i \nabla v \|_2^2 \, dx
=
\int_{\omega_i} \big(
\phi_i^2  \|  \nabla v \|_2^2 +
v^2  \|  \nabla \phi_i \|_2^2 
} \\[3mm]
\displaystyle{
+
2 v \phi_i \nabla v \cdot \nabla \phi_i  \big) \, dx
\leq
\int_{\omega_i} 
  \|  \nabla v \|_2^2 +
C_i^2 v^2  
+
C_i v^2 +   C_i \|  \nabla v \|_2^2 \, dx
} \\
\end{array}
$$
The latter two bounds imply 
\begin{equation}
\label{eq:theorem_b}
\| v_i  \|_{1,\omega_i}^2 \leq
\max\{   C_i + C_i^2 + 1, 2 \}
\| v   \|_{1,\omega_i}^2.
\end{equation}
\end{subequations}

Substituting \eqref{eq:theorem_b} in \eqref{eq:theorem_a}, we obtain
$$
\big| \mathcal{G}_{\mu}(u, v )   \big|
\leq
\sum_{i=1}^{N_{\rm dd}} 
C_i^{\rm r} \,  \mathcal{r}_{\mu}^{(i)} [u] 
\,
\| v \|_{1,\omega_i}.
$$
Then, apply Cauchy-Schwartz  inequality and the estimate in 
\cite[Lemma 2]{BabMel97}, we obtain
$$
\begin{array}{rl}
\displaystyle{\big| \mathcal{G}_{\mu}(u, v )   \big| \leq }
&
\displaystyle{
\max_{i=1,\ldots,N_{\rm dd}} C_i^{\rm r} 
\sqrt{  \sum_{i=1}^{N_{\rm  dd}} \left( \mathcal{r}_{\mu}^{(i)} [u] \right)^2    }
\,
\sqrt{  \sum_{i=1}^{N_{\rm  dd}} \| v \|_{1,\omega_i}^2    }
}
\\[3mm]
\leq &
\displaystyle{
\sqrt{M} 
\max_{i=1,\ldots,N_{\rm dd}} C_i^{\rm r} 
\sqrt{  \sum_{i=1}^{N_{\rm  dd}} \left( \mathcal{r}_{\mu}^{(i)} [u] \right)^2    }
\,
\| v \|_{H^1(\Omega)} 
},
\end{array}
$$
which proves \eqref{eq:residual_estimator}.
\end{proof}

\subsection*{A posteriori error estimation}

As noted in \cref{sec:error_indicator}, the infinite-dimensional analogon of $\mathcal{G}_{\mu}$ as a map from $H^{1}_{0}(\Omega)$ to $H^{-1}(\Omega)$ is not in $C^{1}$; this can be easily seen by determining the Fr\'{e}chet derivative of $\mathcal{G}_{\mu}$. One should thus consider $\mathcal{G}_{\mu}$ as a mapping from $W^{1,p}_{0}(\Omega)$ to $W^{-1,p}(\Omega)$, $p>2$ which yields a $C^{1}$-mapping. Note that the latter is crucial to apply the Brezzi-Rappaz-Raviart (BRR) theory \cite{BrRaRa80,CalRap1997}. One may then derive a rigorous a posteriori error bound for the error in the $| \cdot |_{W^{1,p}(\Omega)}$-norm (see e.g. \cite{PouRap94, CalRap1997, SO17}). Unfortunately, the error estimator relies on the $W^{-1,p}$-norm, which is challenging to estimate even in the finite-dimensional setting. As in \cite{CalRap1997} we thus exploit that for all $z \in B(u_{\rm pum},R)$, $\mathcal{G}_{\mu}^{\prime}(z): W^{1,p}(\Omega) \rightarrow W^{-1,p}(\Omega)$ can be continuously extended as an operator in $L(H^{1}(\Omega), H^{-1}(\Omega))$. Furthermore, we require that \cref{inf-sup-apost-H1}-\cref{Lipschitz-apost-H1} are satisfied. We can then prove \cref{prop: a post est}:
\begin{proof}
Thanks to \cref{th:residual_estimator} the assumption $\tau_{\mu,p} < 1$ implies that $\widetilde{\tau}_{\mu,p} := \frac{2L_{2,p}c_{h}}{\beta_{2,p}^{2}} \|\mathcal{G}_{\mu}( \widehat{u}_{\mu},\cdot)\|_{-1,\Omega} <1$. The existence of a unique solution  $u_{\mu} \in B(\widehat{u}_{\mu},\frac{\beta_{2,p}}{L_{2,p}c_{h}})$ of \cref{eq:variational_form} and 
\begin{equation}\label{eq:aux1}
|\widehat{u}_{\mu} - u_{\mu} |_{1,\Omega} \leq \frac{\beta_{2,p}}{L_{2,p}c_{h}}(1 - \sqrt{1 - \widetilde{\tau}_{\mu,p}})
\end{equation}
then follows using the standard arguments in the BRR theory (see \cite{CaToUr09,VerPat05,CalRap1997,PouRap94} and for this particular PDE \cite{Sme13}). As the function $t(x):= 1 - \sqrt{1 - x}$ is strictly increasing on $(0,1)$, applying \cref{th:residual_estimator} to the right side of \cref{eq:aux1} concludes the proof.
\end{proof}

\subsection*{Extension of   \cref{th:exponential_convergence} to multiple configurations}
In Algorithm
\ref{alg:localized_enrichment_simplified_less},
we   generalize   \cref{alg:localized_enrichment_simplified}  to the case of multiple configurations and  
 we discuss   the proof of \emph{a priori} exponential convergence. Note that   \cref{alg:localized_enrichment_simplified_less} can still be interpreted as  a simplified version of   \cref{alg:localized_enrichment} that is more amenable for the analysis.
 We assume here that 
 $N_{\rm dd, \mu} \leq N_{\rm dd, max}$ with probability one for some $N_{\rm dd, max}$. All constants introduced below depend on 
 $N_{\rm dd, max}$ and also on the size of the training set $n_{\rm train}^{\rm glo}$.
 
\begin{algorithm}[H]                      
\caption{simplified randomized localized training with global enrichment}     
\label{alg:localized_enrichment_simplified_less}     
 
\begin{algorithmic}[1]

\State Initialize $\mathcal{Z} = \mathcal{Z}_0$.

\State Sample $n_{\rm train}^{\rm glo}$ configurations
 $\mu^k  \overset{\rm iid}{\sim} p_{\mu}^{\rm glo}$,
$\mathcal{P}_{\rm train} : =\{ \mu^k  \}_k$
\smallskip

\For{$\ell=0,\ldots,\texttt{maxit}$}

\State
Compute $\widehat{u}_{\ell,\mu}$ using the PUM-CB-ROM (cf.  \cref{sec:domain_decomposition}) for $\mu \in \mathcal{P}_{\rm train}$.
\smallskip

\State
Find $(k,\mu') = {\rm arg} 
\max_{i=1,\ldots,N_{\rm dd,\mu}, \mu \in \mathcal{P}_{\rm train}} \mathfrak{r}_{\mu}^{(i)} [\widehat{u}_{\ell,\mu}]$.
\smallskip

\State
Solve the local problem: find 
$u_{k}^{\rm loc}\in \mathcal{X}_{k,0}$ such that
$\mathcal{G}_{\mu'}(\widehat{u}_{\ell,\mu'}+u_{k}^{\rm loc} , v  ) = 0$ for all $v\in \mathcal{X}_{k,0}$.
\smallskip

\State
Define $u^{\star} = \frac{u_{k}^{\rm loc}}{\phi_k} $ and update the local space 
$\mathcal{Z}  = 
\mathcal{Z}  \cup {\rm span} \{  u^{\star}   \circ \Phi_k    \} $.
 
 \EndFor 
\end{algorithmic}
\end{algorithm}

Given $\mu\in \mathcal{P}_{\rm train}$, 
we denote by $\{ \ell_j^{\mu} \}_j \subset \{1,\ldots,\texttt{maxit} \}$ the set of indices corresponding to the iterations at which we select $\mu'=\mu$ at Line 5 of  \cref{alg:localized_enrichment_simplified_less}.
 Then, exploiting  
\cref{th:exponential_convergence}, we find
$$
\| \widehat{u}_{\ell_{j+1}^{\mu}, \mu} - u_{\mu}  \|_{a_{\mu}}
\leq
\left(
1 - \frac{1}{N_{\rm dd,\mu} c_{\rm pu,\mu}^2 }
\right)^{j/2}
\| \widehat{u}_{0, \mu} - u_{\mu}  \|_{a_{\mu}}
\leq 
C e^{-\alpha \, n_{\rm train}^{\rm glo}  j} 
$$
where $C,\alpha>0$ are constants that do not depend on the iteration count $j$.

Given $\ell\in \mathbb{N}$, there exists $\mu'\in \mathcal{P}_{\rm train}$ such  that
$\# \{ \ell_j^{\mu'} : \ell_j^{\mu'} \leq \ell  \} \geq \frac{\ell}{n_{\rm train}^{\rm glo}}$ --- {after  $\ell$ iterations of \cref{alg:localized_enrichment_simplified_less}, the parameter $\mu'$ is selected more than $\frac{\ell}{n_{\rm train}^{\rm glo}}$ times.}
This implies that 
$$
\max_{i=1,\ldots,N_{\rm dd,\mu}, \mu \in \mathcal{P}_{\rm train}} \mathfrak{r}_{\mu}^{(i)} [\widehat{u}_{\ell,\mu}]
\leq
C \, e^{-\alpha \ell}
$$
and thus
$$
\begin{array}{l}
\displaystyle{
\max_{\mu \in \mathcal{P}_{\rm train}  }
\| \widehat{u}_{\ell, \mu} - u_{\mu}  \|_{a_{\mu}}
=
\max_{\mu \in \mathcal{P}_{\rm train}  }
\|  \mathcal{G}_{\mu} \left(  \widehat{u}_{\ell, \mu} ,\cdot  \right)\|_{ \mathcal{X}_{\mu}'}
\leq
\max_{\mu \in \mathcal{P}_{\rm train}  }
c_{\rm pu, \mu}
\sqrt{
\sum_{i=1}^{N_{\rm dd, \mu}}
\;
\left(
\mathfrak{r}_{\mu}^{(i)} [\widehat{u}_{\ell,\mu}]
\right)^2
}
} \\[3mm]
\displaystyle{
\leq
C \sqrt{ N_{\rm dd, \mu}    }
\left(
\max_{\mu \in \mathcal{P}_{\rm train}  }
c_{\rm pu, \mu}
\right)
e^{-\alpha \ell}.
}\\
\end{array}
$$

\section{Extension to multiple components}\label{suppsec:mult_comp}

As discussed in the main body of the paper, in our numerical experiment,  we consider three components.
We should thus construct the 
local approximation spaces
$\mathcal{Z}^{\rm co},\mathcal{Z}^{\rm ed},\mathcal{Z}^{\rm int}$, with 
$\mathcal{Z}^{\bullet} \subset \mathcal{Y}^{\bullet}$,
such that
\begin{equation}
\label{eq:goal_data_compression_supp}
\min_{\zeta \in \mathcal{Z}^{\texttt{L}_i}}
\; \| u_{\mu} \big|_{\omega_i} - \zeta\circ \Phi_i^{-1}  \|_{1,\omega_i}
\leq \varepsilon_{\rm tol}
\;\;
{\rm for} \; i=1,\ldots,N_{\rm dd},
\;\;
\mu \in \mathcal{P}(n_{\rm dd}),
\end{equation} 
Condition \eqref{eq:goal_data_compression_supp} implies that the local spaces $\mathcal{Z}^{\rm co},\mathcal{Z}^{\rm ed},\mathcal{Z}^{\rm int}$ should approximate the manifolds
\begin{equation}
\label{eq:true_local_manifolds_supp}
\left\{
\begin{array}{l}
\displaystyle{
\mathcal{M}^{\rm int} = \left\{
 u_{\mu} \big|_{\omega_i}  \circ \Phi_i \, : \,
 \mu\in 
\mathcal{P}_{\rm glo}(n_{\rm dd}),
\;
\texttt{L}_i = \texttt{int}, \;
n_{\rm dd} \in \mathbb{N}
\right\},
} \\[3mm]
\displaystyle{
\mathcal{M}^{\rm ed} = \left\{
 u_{\mu} \big|_{\omega_i}  \circ \Phi_i \, : \,
 \mu\in 
\mathcal{P}_{\rm glo}(n_{\rm dd}),
\;
\texttt{L}_i = \texttt{ed}, \;
n_{\rm dd} \in \mathbb{N}
\right\},
} \\[3mm]
\displaystyle{
\mathcal{M}^{\rm co} = \left\{
 u_{\mu} \big|_{\omega_i}  \circ \Phi_i \, : \,
 \mu\in 
\mathcal{P}_{\rm glo}(n_{\rm dd}),
\;
\texttt{L}_i = \texttt{co}, \;
n_{\rm dd} \in \mathbb{N}
\right\}.
} \\
\end{array}
\right.
\end{equation}

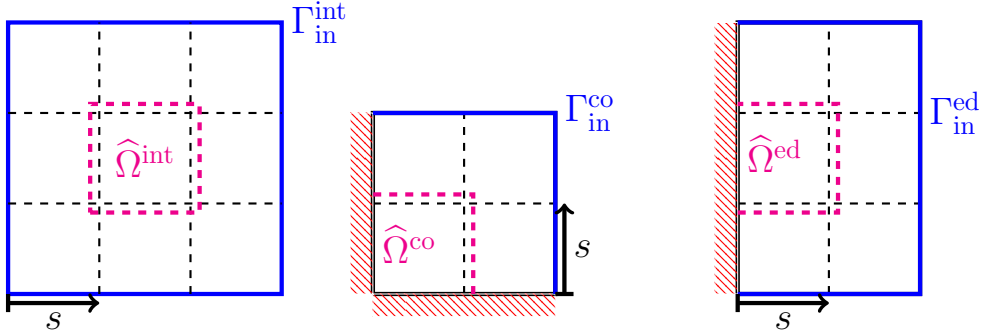
\begin{figure}[h!]
\centering
\begin{tikzpicture}[scale=1.2]
\draw [black,thick] (0,0) -- (3,0) -- (3,3) -- (0,3) --(0,0);
\draw [blue,ultra thick] (0,0) -- (3,0) -- (3,3) -- (0,3) --(0,0);
\coordinate [label={right:  {\Large {$\color{blue} \Gamma_{\rm in}^{\rm int}$}}}] (E) at (3,3) ;
  
\draw [black,thick,dashed] (0,1) -- (3,1);
\draw [black,thick,dashed] (0,2) -- (3,2);
\draw [black,thick,dashed] (1,0) -- (1,3);
\draw [black,thick,dashed] (2,0) -- (2,3);

\draw [magenta,ultra thick,dashed] (0.9,0.9) -- (2.1,0.9) -- (2.1,2.1)--(0.9,2.1) --(0.9,0.9);

\coordinate [label={center:  {\Large {$\color{magenta} \widehat{\Omega}^{\rm int}$}}}] (E) at (1.5,1.5) ;
\draw [->,black,ultra thick] (0,-0.1) -- (1,-0.1);
\draw [black,ultra thick] (0,-0.2) -- (0,0);

\coordinate [label={below:  {\Large {$s$}}}] (E) at (0.5,-0.1) ;
 
\begin{scope}[shift={(4cm,0cm)}]
\draw [black,ultra thick] (0,0) -- (2,0) -- (2,2) -- (0,2) --(0,0);

\draw [black,thick,dashed] (0,1) -- (2,1);
\draw [black,thick,dashed] (1,0) -- (1,2);

\coordinate [label={right:  {\Large {$\color{magenta} \widehat{\Omega}^{\rm co}$}}}] (E) at (0,0.55) ;

\draw[white,pattern=north west lines,pattern color=red] (0,0) -- (2,0)--(2,-0.25) --(0,-0.25) -- (0,0);

\draw[white,pattern=north west lines,pattern color=red] (0,0) -- (0,2)--(-0.25,2) --(-0.25,0) -- (0,0);

\draw [magenta,ultra thick,dashed] (1.1,0) -- (1.1,1.1) -- (0,1.1);
\draw [blue,ultra thick]  (2,0) -- (2,2) -- (0,2) ;

\coordinate [label={right:  {\Large {$\color{blue} \Gamma_{\rm in}^{\rm co}$}}}] (E) at (2,2) ;

\draw [->,black,ultra thick] (2.1,0) -- (2.1,1);
\draw [black,ultra thick] (2,0) -- (2.2,0);
\coordinate [label={right:  {\Large {$s$}}}] (E) at (2.1,0.5) ;
\end{scope}
 
\begin{scope}[shift={(8cm,0.0cm)}]
\draw [black,ultra thick] (0,0) -- (2,0) -- (2,3) -- (0,3) --(0,0);

\draw [black,thick,dashed] (0,1) -- (2,1);
\draw [black,thick,dashed] (0,2) -- (2,2);
\draw [black,thick,dashed] (1,0) -- (1,3);
%\draw [black,thick,dashed] (0,0) -- (1,2);

\coordinate [label={right:  {\Large {$\color{magenta} \widehat{\Omega}^{\rm ed}$}}}] (E) at (0,1.55) ;

\draw[white,pattern=north west lines,pattern color=red] (0,0) -- (0,3)--(-0.25,3) --(-0.25,0) -- (0,0);

\draw [blue,ultra thick]  (0,0) --(2,0) -- (2,3) -- (0,3) ;

\draw [magenta,ultra thick,dashed] (0,0.9) -- (1.1,0.9) -- (1.1,2.1) -- (0,2.1);

\coordinate [label={right:  {\Large {$\color{blue} \Gamma_{\rm in}^{\rm ed}$}}}] (E) at (2,2) ;

\draw [->,black,ultra thick] (0,-0.1) -- (1,-0.1);
\draw [black,ultra thick] (0,-0.2) -- (0,0);
\coordinate [label={below:  {\Large {$s$}}}] (E) at (0.5,-0.1) ;

  \end{scope}
\end{tikzpicture}
\caption{nonlinear diffusion. Archetype components with corresponding oversampling domain.}
\label{fig:nonlineardiff_vis_components}
\end{figure}

For the three components in our library, we consider the oversampling domains depicted in   \Cref{fig:nonlineardiff_vis_components}.
Note that $U^{\bullet}$ comprises $N_{\rm dd}^{\bullet}=9$ (resp., $N_{\rm dd}^{\bullet}=4$, $N_{\rm dd}^{\bullet}=6$) subdomains of $\Omega$ for the internal (resp., corner and edge) component: the active set of parameters is thus equal to 
$\mathcal{P}^{\bullet} = \bigotimes_{i=1}^{N_{\rm dd}^{\bullet}} \widehat{\mathcal{P}} \times \{1,\ldots,N_{\rm dd}^{\bullet}, 0\}$ where $i^{\star} = 0$ means that the source term is outside the patch. 

\subsubsection*{Random boundary conditions}

In   \cref{alg:randomBC_all}, we discuss the complete random sample generator of boundary conditions considered in the numerical experiments.
Note that, if 
 $\bullet \in  \{ \texttt{co}, \texttt{ed} \}$, 
since by construction $\frac{d^k}{d s^k}  g^{(1)}(0) = \frac{d^k}{d s^k}  g^{(1)}(1)$ for $k\in \mathbb{N}$,  
  we   define    $g^{(2)}(s )  = g^{(1)} \, ( c   s )$ with $c=0.7$\footnote{The  choice $c = 0.7$ is not crucial for the methodology.} (cf. Line 7);
then, we enforce that $g(s) = 0$ for $s\in \{ 0,1 \}$ and $g\geq 0$ (cf. Line 8); finally, in Line 9, we ensure that
${\rm Im} [g] \subset [0, \bar{u}_{\rm max}]$.

\begin{algorithm}                   
\caption{Random sample generator of boundary conditions}     
\label{alg:randomBC_all}     

\small
\begin{flushleft}
\emph{Inputs:}  $N_{\rm f}, \alpha$  (cf. 
\eqref{eq:randomBC_basic}),
$\bar{u}_{\rm max} \in (0,1]$, $\bullet \in \{\texttt{co}, \texttt{ed}, \texttt{int}\}$.
\smallskip

\emph{Output:} 
$g: [0,1]\to [0,1)$
boundary condition.
\end{flushleft}                      

 \normalsize 

\begin{algorithmic}[1]

 \State
Draw $\mathbf{c}^{\rm re}, \mathbf{c}^{\rm im} \in \mathbb{R}^{N_{\rm f}}$ s.t.
$c_k^{\rm re}, c_k^{\rm im} \overset{\rm iid}{\sim} 
 \mathcal{N}(0,1)$.
 \smallskip
 
\State
Draw $X_1,X_2,X_3 
\overset{\rm iid}{\sim} 
 {\rm Uniform}(0,\bar{u}_{\rm max}) 
 $, set $a = \min\{X_1,X_2\}$,  $b = \max\{X_1,X_2\}$.
 \smallskip

\State
Set  $\displaystyle{
g^{(1)} 
\, = \, {\texttt{Real}} 
\left[ 
\widetilde{g}(\cdot; \mathbf{c}^{\rm re}, \mathbf{c}^{\rm im}) \right]}.
$
 \smallskip

\If{ $\bullet = \texttt{int}$ }

\State $g = a + \frac{b-a}{\max g^{(1)}  - \min g^{(1)} } \left( g^{(1)} -  \min g^{(1)}  \right)$.
 \smallskip
 
\Else

\State $g^{(2)}(s )  = g^{(1)} \, ( 0.7   s )$.
 \smallskip
 
\State $g^{(3)}(s)  =   
\left(
a + \frac{b-a}{\max g^{(2)}  - \min g^{(2)} } \left( g^{(2)} -  \min g^{(2)}  \right)
\right)
 s (1 - s)$.
 \smallskip
 
\State $g =   \frac{ X_3  }{ \max  g^{(3)}  }   g^{(3)}$.

\EndIf
  
\end{algorithmic}
\end{algorithm}

\subsubsection*{Adaptive enrichment}

\Cref{alg:localized_enrichment_all} is the generalization of   \cref{alg:localized_enrichment} to the case of multiple archetype components.
Note that we  here enforce 
that the dimension of the local spaces is the same for all components to simplify the implementation of the CB-ROM;
however, the Algorithm can be modified to allow local spaces of different size.  
 
\begin{algorithm}                 
\caption{randomized localized training with global enrichment}     
\label{alg:localized_enrichment_all}     

\small
\begin{flushleft}
\emph{Inputs (localized training):}  
$n_{\rm train}^{\rm loc}$ = number of  solves,
 $n^{\rm loc}$ = size of the POD spaces,
 $\{ p_{\mu}^{\bullet}, p_{\rm bc}^{\bullet}   \}_{\bullet}$
 sampling pdfs.
 \medskip
 
 \emph{Inputs (enrichment):}  
$n_{\rm train}^{\rm glo}$ = number of global simulations per iteration,
 $n^{\rm glo}$ = number of   modes added at each iteration,
 $\texttt{maxit}$ = maximum number of outer loop iterations,
 $tol$ = tolerance for termination criterion,
 $p_{\mu}^{\rm glo}$ = global configuration sampler,
 $m_{\rm r}$ = percentage of marked components at each iteration.
 \medskip
 
\emph{Outputs:} 
$\{  \mathcal{Z}^{\bullet} \}_{\bullet \in \{ \texttt{co}, \texttt{ed}, \texttt{int} \} }$ local approximation spaces.
\end{flushleft}                      

 \normalsize 

\begin{center}
\textbf{Localized training}
\end{center}

\begin{algorithmic}[1]
\State
Apply   \cref{alg:localized_training} to obtain the local spaces
$\{  \mathcal{Z}^{\bullet} \}_{\bullet \in \{ \texttt{co}, \texttt{ed}, \texttt{int} \} }$.
\end{algorithmic} 

 \begin{center}
\textbf{Enrichment}
\end{center}

\begin{algorithmic}[1]

\State Sample $n_{\rm train}^{\rm glo}$ configurations
 $\mu^{(k)} \overset{\rm iid}{\sim} p_{\mu}^{\rm glo}$,
$\mathcal{P}_{\rm train} : =\{ \mu^{(k)}  \}_k$
\smallskip

\For{$\ell=1,\ldots,\texttt{maxit}$}

\State Initialize the datasets $  \mathcal{D}^{\bullet} = \emptyset$ for $\bullet \in \{ \texttt{co}, \texttt{ed}, \texttt{int} \} $.
\smallskip

\For{ $\mu \in \mathcal{P}_{\rm train}$}

\State
Compute $\widehat{u}_{\mu}$ using the PUM-CB-ROM (cf.   \cref{sec:domain_decomposition}).
\smallskip

\For{ $\bullet \in \{ \texttt{co}, \texttt{ed}, \texttt{int} \}$}

\State
Compute local residuals 
\eqref{eq:residual_localized}
$ \mathfrak{r}_{\mu}^{i} =  
\mathfrak{r}_{\mu}^{(i)} [  \widehat{u}_{\mu} ] $ 
for 
$i=1,\ldots,N_{\rm dd,\mu}$ s.t. $\texttt{L}_i = \bullet$.
\smallskip

\State
Mark the $m_{\rm r}$ $\%$ instantiated components of type $\bullet$
with the largest  residual,
$\{  V_{i} \}_{i \in \texttt{I}_{\rm mark, \bullet}^{\mu}   }$.

\smallskip

\State
Solve the local problems
\eqref{eq:local_solution}
in
$\{  \omega_{i} \}_{i \in \texttt{I}_{\rm mark, \bullet}^{\mu}  }$,
$u_{i,\mu}^{\bullet} = \frac{1}{\phi_i} T_{\mu}^{(i)} (  \widehat{u}_{\mu} |_{ \omega_{i}})$.
\smallskip

\State
Augment the dataset $\mathcal{D}^{\bullet} = \mathcal{D}^{\bullet} \cup \{
u_{i,\mu}^{\bullet} \circ \Phi_{i} : 
i \in \texttt{I}_{\rm mark, \bullet}^{\mu}   \} $.
\smallskip

\EndFor

\State
Compute $\Delta_{\mu}$   \eqref{eq:residual_adaptive}.

\EndFor 

\State
Update the POD spaces
$\mathcal{Z}^{\bullet} = 
\mathcal{Z}^{\bullet}
\cup
{\rm POD}\left( \{ w - \Pi_{\mathcal{Z}^{\bullet}} w : w\in \mathcal{D}^{\bullet} \}, (\cdot,\cdot)_{\bullet}, n^{\rm glo}       \right)$.

\If{ $\max_{ \mu \in \mathcal{P}_{\rm train}} \Delta_{\mu}    
 < tol$ }

\State  \texttt{BREAK}
 \smallskip
 
\EndIf

\EndFor 
\end{algorithmic}
\end{algorithm}

\footnotesize

 {\bibliographystyle{abbrv}
\bibliography{random}}

\end{document}